%% file: CochranArchive.tex
\documentstyle{amltd}
\begin{document}
\annalsline{157}{2003}
\received{August 24, 1999}
\startingpage{433}
\def\bye{\input Cochran.refs\end{document}}
 \font\tenrm=cmr10
\def\ritem#1{\item[{\rm #1}]}
\input amssym.def
\input amssym.tex
\input boxedeps.tex 
\SetepsfEPSFSpecial 
\HideDisplacementBoxes
\def\figin#1#2{
$$
 {\BoxedEPSF{#1.eps scaled
#2}%
}%
$$
\noindent}
\def\xrightarrow#1{\stackrel{#1}{\rightarrow}}
\def\joinrel{\mathrel{\mkern-4mu}}
\def\relbar{\mathrel{\smash-}}
\def\lrar{\relbar\joinrel\relbar\joinrel\rightarrow}
\def\vlrar{\relbar\joinrel\relbar\joinrel\relbar\joinrel\relbar\joinrel\rightarrow}
\def\eqref#1{(\ref{#1})}

\catcode`\@=11
\font\twelvemsb=msbm10 scaled 1100
\font\tenmsb=msbm10
\font\ninemsb=msbm10 scaled 800
\newfam\msbfam
\textfont\msbfam=\twelvemsb  \scriptfont\msbfam=\ninemsb
  \scriptscriptfont\msbfam=\ninemsb
\def\msb@{\hexnumber@\msbfam}
\def\Bbb{\relax\ifmmode\let\next\Bbb@\else
 \def\next{\errmessage{Use \string\Bbb\space only in math
mode}}\fi\next}
\def\Bbb@#1{{\Bbb@@{#1}}}
\def\Bbb@@#1{\fam\msbfam#1}
\catcode`\@=12

 \catcode`\@=11
\font\twelveeuf=eufm10 scaled 1100
\font\teneuf=eufm10
\font\nineeuf=eufm7 scaled 1100
\newfam\euffam
\textfont\euffam=\twelveeuf  \scriptfont\euffam=\teneuf
  \scriptscriptfont\euffam=\nineeuf
\def\euf@{\hexnumber@\euffam}
\def\frak{\relax\ifmmode\let\next\frak@\else
 \def\next{\errmessage{Use \string\frak\space only in math
mode}}\fi\next}
\def\frak@#1{{\frak@@{#1}}}
\def\frak@@#1{\fam\euffam#1}
\catcode`\@=12
\def\operatorname#1{\mathop{\rm #1}\nolimits}
\newcommand{\spec}{\operatorname{spec}}
\newcommand{\rk}{\operatorname{rk}_{\KK}} 
\newcommand{\Herm}{\operatorname{Herm}}
\newcommand{\GL}{\operatorname{GL}}
\newcommand{\Ker}{\operatorname{Ker}}
\newcommand{\End}{\operatorname{End}}
\newcommand{\Hom}{\operatorname{Hom}}
\newcommand{\Ext}{\operatorname{Ext}}
\newcommand{\Rep}{\operatorname{Rep}}
\newcommand{\Tors}{\operatorname{Tors}}
\newcommand{\dom}{\operatorname{dom}}
\newcommand{\tr}{\operatorname{tr}}
\newcommand{\id}{\operatorname{id}}

\newcommand{\Z}{{\Bbb Z}}
\newcommand{\N}{{\Bbb N}}
\newcommand{\C}{{\Bbb C}}
\newcommand{\Q}{{\Bbb Q}}
\newcommand{\K}{{\Bbb K}}
\newcommand{\R}{{\Bbb R}}
\newcommand{\F}{{\Bbb F}}
\newcommand{\f}{\Bbbk}
\newcommand{\G}{\Gamma}

\newcommand{\RR}{{\cal R}}
\newcommand{\UU}{{\cal U}}
\newcommand{\NN}{{\cal N}}
\newcommand{\DD}{{\cal D}}
\newcommand{\KK}{{\cal K}}
\newcommand{\FF}{{\cal F}}
\newcommand{\MM}{{\cal M}}
\newcommand{\LL}{{\cal L}}
\newcommand{\CC}{{\cal C}}
\newcommand{\BB}{{\cal B}}
\newcommand{\QQ}{{\cal Q}}
\newcommand{\PP}{{\cal P}}

\newcommand{\e}{\epsilon}
\newcommand{\hra}{\hookrightarrow}
\newcommand{\imra}{\looparrowright}
\newcommand{\sra}{\twoheadrightarrow}
\newcommand{\ira}{\rightarrowtail}
\newcommand{\sd}{\rtimes}
\newcommand{\ra}{\longrightarrow}
\newcommand{\la}{\longleftarrow}
\newcommand{\lr}{\longleftrightarrow}

\renewcommand{\SS}{{\cal S}}
\renewcommand{\AA}{{\cal A}}
\renewcommand{\l}{\ell}
\renewcommand{\a}{\alpha}
\renewcommand{\i}{\iota}
\renewcommand{\b}{\beta}

\newcommand{\ngth}{\negthickspace}
\newcommand{\torsionp}{\Z_{(p)}/\Z}
\newcommand{\defeq}{\stackrel{\mathrm{def}}{=}}
\newcommand{\2}{\Z[\pi/\pi^{(2)}]}
\newcommand{\n}{\Z[\pi/\pi^{(n)}]}
\newcommand{\sbq}{\subseteq}
\newcommand{\spq}{\supseteq}
\newcommand{\sqp}{\sqsupset}
\newcommand{\thra}{\twoheadrightarrow}
\newcommand{\gu}{{\G_n^U}}
\newcommand{\gop}{{\G_0^U}}


\font\emi= cmmi10 scaled 1700 
\font\eightmi=cmmi10
\font\titr=cmr10
\font\fiver=cmr5

\title{Knot concordance, Whitney towers\\ and {\emi L}\raise8pt\hbox{\titr 2}-signatures}  
 \def\titleheadline#1{\def\one{#1}\ifx\one\empty\else
\gdef\thetitle{{\frenchspacing%
\let\\ \relax
{#1}}}\fi}
\newif\ifshort
\def\shortname#1{\global\shorttrue\xdef
\theauthors{{\eightsc\uppercase{#1}}}}
\let\shorttitle\titleheadline
\shorttitle{ \eightsc\uppercase{Knot concordance, Whitney towers and} {\eightpoint \it L}
\hskip-3pt\raise3pt\hbox{{\fiver 2}}\eightsc\uppercase{-signatures}}

 \acknowledgements{All authors were supported by MSRI and NSF.  The third   author was also supported by a
fellowship from the Miller   foundation, UC Berkeley.}
 \twoauthors{Tim D. Cochran, Kent E. Orr,}{Peter Teichner}
\institutions{Rice University, Houston, Texas\\
{\eightpoint {\it E-mail address\/}: cochran@math.rice.edu}\\
\vglue6pt
Indiana University, Bloomington, Indiana\\
{\eightpoint {\it E-mail address\/}:  korr@indiana.edu}\\
\vglue6pt
University of California in San Diego, La Jolla, California\\
{\eightpoint {\it E-mail address\/}: 
teichner@math.ucsd.edu}}

\centerline{\bf Abstract} \vglue6pt

 We construct many examples of nonslice knots in 3-space that
cannot be distinguished from slice knots by previously known invariants.  Using Whitney towers in place of
embedded disks, we define a geometric filtration of the 3-dimensional topological knot concordance group.  The
bottom part of the filtration exhibits all classical concordance invariants, including the Casson-Gordon
invariants. As a first step, we construct an infinite sequence of new obstructions
that vanish on slice knots.  These take values in the $L$-theory of skew fields associated to certain
{\it  universal} groups. Finally, we use the dimension theory of von Neumann algebras to define
an $L^2$-signature and use this to detect the first unknown step in our
obstruction theory.

\def\sni#1{\vglue-.5pt\noindent{#1}. }
\def\ssni#1{\vglue-1pt\noindent\hskip18pt {#1}.}

\vfill
\centerline{\bf Contents}
\vfill
\sni{1} Introduction
\ssni{1.1} Some history, $(h)$-solvability and Whitney towers
\ssni{1.2} Linking forms, intersection forms, and solvable representations of\\ \phantom{zooo 2.\ }  knot groups
\ssni{1.3} $L^2$-signatures
\ssni{1.4} Paper outline and acknowledgements
\sni{2} Higher order Alexander modules and Blanchfield linking forms
\sni{3} Higher order linking forms and solvable representations of the knot group
\sni{4} Linking forms and Witt invariants as obstructions to solvability
\sni{5} $L^2$-signatures
\sni{6} Non-slice knots with vanishing Casson-Gordon invariants
\sni{7} $(n)$-surfaces, gropes and Whitney towers
\sni{8} $H_1$-bordisms 
\sni{9} Casson-Gordon invariants and solvability of knots
\vglue-.5pt\noindent \phantom{9. \ } References
 
\section{Introduction}

This paper begins a detailed investigation into the group
 of topological concordance classes of knotted circles in the $3$-sphere. Recall
that a knot $K$ is topologically  {\it  slice } if there exists a locally flat
topological embedding of the $2$-disk into $B^4$ whose restriction to the boundary
is $K$. The knots
$K_0$ and $K_1$ are topologically
 {\it  concordant} if there is a locally flat topological
 embedding of the annulus into $S^3 \times [0,1]$ whose restriction to the
boundary components gives the knots. The set of
 concordance classes of knots under the operation of connected sum
 forms an abelian group $\CC$, whose identity element is the class of slice knots. 
 
\proclaimtitle{A special case}
\specialnumber{6.4}\proclaim{Theorem}  The knot of Figure~{\rm 6.1}   has vanishing
Casson\/{\rm -}\/Gordon invariants but is not topologically slice.
\endproclaim

  In fact, we construct infinitely many such examples that cannot be
distinguished from slice knots by previously known invariants. The new slice
obstruction that detects these knots is an $L^2$-signature formed from the
dimension theory of the von Neumann algebra of a certain rationally universal
solvable group.  To construct nontrivial maps from the fundamental group of the
knot complement to this solvable group, we develop an obstruction theory and for
this purpose, we define noncommutative higher-order versions of the classical
Alexander module and Blanchfield linking form. We hope that these generalizations
are of considerable independent interest. 
\vglue4pt

We give new geometric conditions which lead to a natural filtration of the slice
condition ``there is an embedded
$2$-disk in $B^4$ whose boundary is the knot''.  More precisely, we exhibit a new
geometrically defined filtration of the knot concordance group
$\CC$ indexed on the half integers;
$$
\cdots\subset \FF_{(n.5)} \subset \FF_{(n)} \subset \cdots \subset
\FF_{(0.5)} \subset \FF_{(0)} \subset \CC,
$$ where for $h\in\frac{1}{2}\N_0$, the group $\FF_{(h)}$ consists of all {\it  $(h)$-solvable} knots.\break
$(h)$-solvability is defined using intersection forms in certain solvable covers (see Definition~\ref{def:solvable}).  The
obstruction theory mentioned above measures whether a given knot lies in the subgroups
$\FF_{(h)}$. It provides a bridge from algebra to the
 topological techniques of A. Casson and M. Freedman.  In fact,
$(h)$-solvability has an equivalent definition in terms of the geometric notions
of {\it  gropes} and {\it  Whitney towers} (see Theorems~\ref{half} and
\ref{ndashsolv} in part~1.1   of the introduction).  Moreover, the
tower of von Neumann signatures might be viewed as an algebraic mirror of
infinite constructions in topology. Another striking example of this bridge is
the following theorem, which implies that the Casson-Gordon invariants obstruct a
specific step (namely a second \pagebreak layer of Whitney disks) in the
Freedman-Cappell-Shaneson surgery theoretic program to prove that a knot is
slice.  Thus one of the most significant aspects of our work is to provide a step
toward a new and strictly 4-dimensional homology surgery theory.

\specialnumber{9.11}\proclaim{Theorem} Let $K\subset S^3$ be
$(1.5)$\/{\rm -}\/solvable. Then all previously known concordance invariants of $K$ vanish.
\endproclaim

  In addition to the Seifert form obstruction, these are the 
 invariants introduced by A. Casson and C. McA. Gordon in 1974 and further
metabelian invariants by P. Gilmer~\cite{G1}, \cite{G2}, P. Kirk and C.
Livingston~\cite{KL}, and\break C. Letsche~\cite{Let}.  More precisely,
Theorem~\ref{vanish} actually proves the vanishing of the Gilmer invariants. 
These determine the Casson-Gordon invariants and the invariants of Kirk and
Livingston.  The Letsche obstructions are handled in a separate
Theorem~\ref{L=0}. 

The first few terms of our filtration correspond
closely to the previously known concordance invariants and we show that the
filtration is nontrivial beyond these terms. Specifically, a knot lies in
$\FF_{(0)}$ if and only if it has vanishing Arf invariant, and lies in
$\FF_{(0.5)}$ if and only if  it is algebraically slice, i.e.\ if the Levine Seifert form obstructions (that classify higher
dimensional knot concordance) vanish (see Theorem~\ref{algebraically slice} together with
Remark~\ref{rem:$(n)$-solvable}). Finally, the family of examples of Theorem~\ref{Q slice}
 proves the following:

\nonumproclaim{{C}orollary}
 The quotient group $\FF_{(2)}/\FF_{(2.5)}$ has infinite rank.
\endproclaim

 In this paper we will show that this quotient group is nontrivial.  The full proof
of the corollary will appear in another paper.

The geometric relevance of our filtration is further revealed by the following
two results, which are explained and proved in Sections~\ref{sec:surfaces}
and~\ref{sec:$H_1$-bordism}.

\specialnumber{8.11} \proclaim{Theorem}  If a knot $K$ bounds a grope of height~$(h+2)$ in $D^4$ then
$K$ is $(h)$\/{\rm -}\/solvable.
\endproclaim

 \vglue-36pt
\centerline{\BoxedEPSF{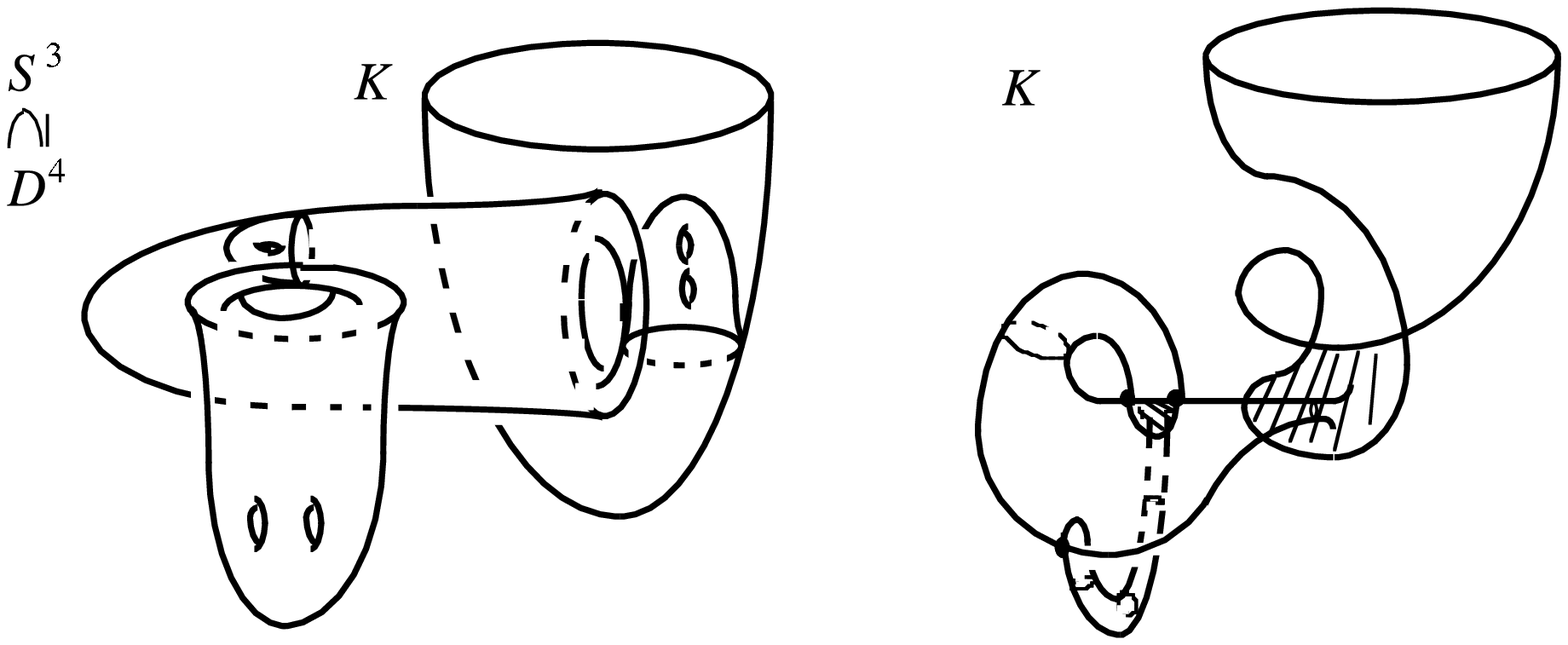 scaled 525}}
\vglue6pt
\centerline{Figure 1.1. A grope of height $2.5$ and a Whitney tower of height
$2.5$.}

\specialnumber{8.12} \proclaim{Theorem} If a knot $K$ bounds a Whitney tower of height~$(h+2)$ in
$D^4$ then $K$ is $(h)$\/{\rm -}\/solvable.
\endproclaim

We establish an infinite series of new
knot slicing obstructions lying in the $L$-theory of large skew fields, and
associated to the commutator series of the knot group.  These successively
obstruct each integral stage of our filtration~(Theorem~\ref{invariant thm}). We
also prove the desired result that the higher-order Alexander modules of an
$(h)$-solvable knot contain submodules that are self-annihilating with respect to
the corresponding higher-order linking form. We see no reason that this tower of
obstructions should break down after  three steps even though the complexity of
the computations grows.  We conjecture:

\nonumproclaim{Conjecture} For any $n \in \N_0${\rm ,} there are
$(n)$\/{\rm -}\/solvable knots that are not $(n.5)$\/{\rm -}\/solvable.  In fact $\FF_{(n)}/\FF_{(n.5)}$ has
infinite rank.
\endproclaim

For $n = 0$ this is detected by the Seifert form obstructions, for $n = 1$ this
can be established by Theorem~\ref{vanish} from examples due to Casson and Gordon,
and $n=2$ is the above corollary. Indeed, if there exists a fibered ribbon knot
whose  classical Alexander module, first-order Alexander module $\dots$ and\break
$(n-1)^{\rm st}$-order Alexander module have {\it  unique} proper submodules (analogous to
$\Z_9$ as opposed  to $\Z_3\times\Z_3$), then the conjecture is true for all $n$.
Hence our inability to establish the full
conjecture at this time seems to be merely a technical deficiency related to the
difficulty of solving equations over noncommutative fields. In
Section~\ref{sec:$H_1$-bordism} we will explain what it means for an arbitrary
link to be $(h)$-solvable. Then the following result provides plenty of
candidates for proving our conjecture in general.

\specialnumber{8.9} \proclaim{Theorem} If there exists an $(h)$\/{\rm -}\/solvable link which forms a
standard half basis of untwisted curves on a Seifert surface for a knot $K${\rm ,} then
$K$ is
$(h+1)$\/{\rm -}\/solvable.
\endproclaim

It remains open whether a $(0.5)$-solvable knot is
$(1)$-solvable and whether a $(1.5)$-solvable knot is
$(2)$-solvable but we do introduce potentially nontrivial obstructions that
generalize the Arf invariant (see Corollary~\ref{spin}).

\vglue8pt 1.1. {\it Some history{\rm ,} $(h)$\/{\rm -}\/solvability and Whitney towers}. 
In the 1960's,\break M. Kervaire and J. Levine computed the group of concordance
classes of knotted
$n$-spheres in $S^{n+2}$, $n \geq 2$, using ambient surgery techniques.
Even-dimensional knots are always slice~\cite{K}, and the odd-dimensional
concordance group can be described by a collection of computable obstructions
defined as Witt equivalence classes of linking pairings on a Seifert
surface~\cite{L1} (see also
\cite{Sto}).  One modifies the Seifert surface along middle-dimensional embedded
disks in the $(n+3)$-ball to create the slicing disk.  \pagebreak The obstructions to
embedding these middle-dimensional disks are intersection numbers that are
suitably reinterpreted as linking numbers of the bounding homology classes in the
Seifert surface. This {\it  Seifert form} obstructs slicing knotted
$1$-spheres as well.

In the mid  1970's, S. Cappell and J. L. Shaneson introduced a new strategy for
slicing knots by extending surgery theory to a theory classifying manifolds
within a homology type~\cite{CS}.  Roughly speaking, the classification of
higher dimensional knot concordance is the classification of homology circles up
to homology cobordism rel boundary. The reader should appreciate the basic fact
that a knot is a slice knot if and only if the (n+2)-manifold, $M$, obtained by
(zero-framed) surgery  on the knot is the boundary of a manifold that has the
homology of a circle and whose fundamental group is normally generated by the
meridian of the knot. More generally, for knotted
$n$-spheres in $S^{n+2}$ (n odd), here is an outline
 of the Cappell-Shaneson 
 surgery strategy. One lets $M$ bound an 
$(n+3)$-manifold $W$ with infinite cyclic fundamental group.  The
middle-dimensional homology of the universal abelian cover of $W$ admits a
$\Z[\Z]$-valued intersection form.  The Cappell-Shaneson obstruction is the
obstruction to finding a half-basis of immersed spheres whose intersection points
occur in pairs each of which admits an associated immersed Whitney disk.  As
usual, in  higher dimensions, if the obstructions vanish, these Whitney disks may
be embedded and intersections removed in pairs.  The resulting embedded spheres
are then surgically excised resulting in an homology circle, i.e.\ a slice
complement. 

These two strategies, when applied to the case $n=1$, yield the following
equivalent obstructions. (See~\cite{L1} and
\cite{CS} together with   Remark~\ref{rem:$(n)$-solvable}.2.) The theorem is
folklore except that condition~(c) is new (see Theorem~\ref{2.5}). Denote by $M$
the
$0$-framed surgery on a knot $K$. Then
$M$ is a closed $3$-manifold and $H_1(M):=H_1(M;\Z)$ is infinite cyclic. An
orientation of $M$ and a generator of
$H_1(M)$  are determined by orienting $S^3$ and $K$. 

\proclaim{Theorem} \label{algebraically slice} The following statements are equivalent\/{\rm :} 
\begin{itemize}
\ritem{(a)} {\rm (}\/The Levine condition\/{\rm )}\/  $K$ bounds a Seifert surface in $S^3$ for which
the Seifert   form contains a Lagrangian. 
\ritem{(b)} {\rm (}\/The Cappell\/{\rm -}\/Shaneson condition\/{\rm )} $M$ bounds a compact spin manifold $W$
with the following properties\/{\rm :}\/
\begin{itemize}
\ritem{1.} The inclusion induces an isomorphism $H_1(M)  
\xrightarrow{\cong} H_1(W)$.
\ritem{2.} The $\Z[\Z]$\/{\rm -}\/valued intersection form $\lambda_1$ on  
$H_2(W;\Z[\Z])$ contains a totally isotropic submodule whose image is a
Lagrangian in
$H_2(W)$.
\end{itemize}
\ritem{(c)} $K$ bounds a grope of height~$2.5$ in $D^4$.
\end{itemize}

\endproclaim

A submodule is {\it  totally isotropic} if the corresponding form vanishes on it.
A {\it  Lagrangian} is a totally isotropic direct summand of half rank.  Knots
satisfying the conditions of Theorem~\ref{algebraically slice} are the
aforementioned class of {\it  algebraically   slice} knots. In particular, slice
knots satisfy these conditions, and in higher dimensions, Levine showed
 that algebraically slice implies slice \cite{L1}.  

If the Cappell-Shaneson homology surgery machinery worked in dimension four,
algebraically slice knots would be slice as well. However, in the mid $1970's$,
Casson and Gordon discovered new slicing obstructions proving that, contrary to
the higher dimensional case, algebraically slice knotted
$1$-spheres are not necessarily slice~\cite{CG1}, \cite{CG2}. The problem is that
the Whitney disks that pair up the intersections of a spherical Lagrangian may no
longer be embedded, but may themselves have intersections, which might or might
not
 occur in pairs, and if so may have their own Whitney disks. One naturally
speculates that the Casson-Gordon invariants should obstruct a second layer of
Whitney disks in this approach. This is made precise by Theorem~\ref{vanish}
together with the following theorem (compare Definitions~\ref{def:Whitney
tower},~\ref{def:n-solv}~and~\ref{def:n.5-solv}). Moreover this theorem shows that
$(h)$-solvability  filters the Cappell-Shaneson approach to disjointly embedding
an integral homology half basis of spheres in the $4$-manifold.

\specialnumber{8.4\ \&\ 8.8} \proclaim{Theorems}  A knot is $(h)$\/{\rm -}\/solvable if and only if
$M$ bounds a compact spin manifold $W$ where the inclusion induces an isomorphism on
$H_1$ and such that there exists a Lagrangian $L\subset H_2(W;\Z)$  that has the
following additional geometric property\/{\rm :} $L$ is generated by immersed spheres
$\l_1,\dots,\l_k$ that allow a {\rm  Whitney tower} of height~$h$.
\endproclaim

We conjectured above that there is a nontrivial step from each height of the
Whitney tower to the next. However, even an infinite Whitney tower might not lead
to a slice disk. This is in contrast to finding {\it  Casson towers}, which in
addition to the Whitney disks have so called {\it  accessory disks} associated to
each double point. By Freedman's main result, any Casson tower of height four
contains a topologically embedded disk.  Thus the ultimate goal is to establish
necessary and sufficient criteria to finding Casson towers. Since a Casson tower
is in particular a Whitney tower, our obstructions also apply to Casson towers.
For example, it follows that Casson-Gordon invariants obstruct finding Casson
towers of height two in the above Cappell-Shaneson approach. Thus we provide a
proof of the heuristic argument that by Freedman's result the Casson-Gordon
invariants must obstruct the existence of Casson towers.
\vglue4pt 

We now outline the definition of $(h)$-solvability. The reader can see that it
filters the condition of finding a
 half-basis of disjointly embedded spheres by examining intersection forms with
progressively more discriminating coefficients, as indexed by the derived series. 

Let $G^{(i)}$ denote the $i^{\rm th}$ {\it  derived group} of a group $G$, inductively
defined by $G^{(0)}:=G$ and
$G^{(i+1)}:=[G^{(i)}, G^{(i)}]$. A group $G$ is {\it  $(n)$-solvable} if $G^{(n+1)}=1$  ($(0)$-solvable corresponds to abelian) and
$G$ is {\it  solvable} if such a finite $n$ exists. For a CW-complex $W$, we define
$W^{(n)}$ to be the regular covering corresponding to the subgroup
$(\pi_1(W))^{(n)}$. If
$W$ is an oriented $4$-manifold then there is an intersection form
$$
\lambda_{n}:H_2(W^{(n)}) \times H_2(W^{(n)}) 
\ra \Z[\pi_1(W)/\pi_1(W)^{(n)}].
$$  (see ~\cite[Ch.\ 5]{Wa}, and our \S \ref{sec:surfaces}  where we also
explain the   self-intersection invariant $\mu_n$). 
  For $n\in\N_0$, an {\it  $(n)$-Lagrangian} is a submodule
$L\subset H_2(W^{(n)})$ on which $\lambda_n$ and $ \mu _n$ vanish and which maps
onto a Lagrangian of $\lambda_0$.

\numbereddemo{Definition}\label{def:solvable} A knot is called {\it  $(n)$-solvable} if $M$ bounds a spin $4$-manifold $W$,
 such that the inclusion map induces an isomorphism  on first homology and such
that W admits two {\it  dual}
$(n)$-Lagrangians. This means that the form $\lambda_n$ pairs the two Lagrangians
nonsingularly and that their images together freely generate
$H_2(W)$ (see Definition~\ref{def:Lagrangian}).

A knot is called {\it  $(n.5)$-solvable}, $n\in \N_0$, if $M$ bounds a spin
$4$-manifold $W$ such that the inclusion map induces an isomorphism on first
homology and such that W admits an
$(n+1)$-Lagrangian and a dual $(n)$-Lagrangian in the above sense. We  say that
$M$ is {\it  $(h)$-solvable via $W$} which is called an {\it  $(h)$-solution for $M$} (or $K$).
\enddemo

\numbereddemo{{R}emark}\label{rem:$(n)$-solvable} It is appropriate to mention the following
facts:
\begin{itemize}
\ritem{1.} The size of an $(h)$-Lagrangian $L$ is controlled only by its image in
$H_2(W)$; in particular, if $H_2(W)=0$ then the knot $K$ is $(h)$-solvable for all
$h\in\frac{1}{2}\N$. This holds for example if $K$ is topologically slice. More
generally, if $K$ and $K'$ are topologically concordant knots, then $K$ is
$(h)$-solvable if and only if $K'$ is
$(h)$-solvable. (See Remark~\ref{rem:homology cobordism}.)
\item[2.] One easily shows $(0)$-solvable knots are exactly knots with trivial Arf
invariant. (See Remark~\ref{rem:Arf}.)  One sees that a knot is algebraically
slice if and only if it is $(0.5)$-solvable by observing that the definition
above for $n=0$ is exactly condition (b.2) of Theorem~\ref{algebraically slice}.
\item[3.] By the naturality of covering spaces and homology with   twisted
coefficients, if
$K$ is $(h)$-solvable then it is  
$(h')$-solvable for all $h'\leq h$.
\item[4.] Given an $(n.5)$-solvable or $(n)$-solvable knot with a $4$-manifold $W$ as
in Definition~\ref{def:solvable} one can do surgery on elements in
$\pi_1(W^{(n+1)})$, preserving all the conditions on $W$. In particular, if
$\pi_1(W)/\pi_1(W)^{(n+1)}$ is finitely presented then one can arrange for
$\pi_1(W)$ to be $(n)$-solvable. This motivated our choice of terminology. Moreover,
since this condition does hold for $n=0$, we see that, in the classical case of
$(0.5)$-solvable (i.e., algebraically slice) 
  knots, one can always assume that $\pi_1(W)=\Z$. This is the way that
condition~(b) in Theorem~\ref{algebraically slice} is usually formulated, namely
as the vanishing of the Cappell-Shaneson surgery obstruction in  
$ \G_0(\Z[\Z]\to\Z).$ In particular, this proves the equivalence of   conditions
(a) and (b) in Theorem~\ref{algebraically slice}. The equivalence of (b) and (c)
will be proved in Section~\ref{sec:surfaces}.
\end{itemize}

\enddemo

\vglue-16pt
  1.2. {\it Linking forms{\rm ,} intersection forms{\rm ,} and solvable representations of
knot groups.}
The Casson-Gordon invariants exploit the observation that linking of
$1$-dimensional objects in a $3$-manifold may be computed via the intersection
theory of a homologically simple
$4$-manifold that it bounds. Thus, $2$-dimensional intersection pairings for the
$4$-manifold are subtly related to the fundamental group of the bounding
$3$-manifold. Casson and Gordon utilize the $\Q/\Z$, or {\it  torsion linking
pairing}, on prime power cyclic knot covers to access intersection data in 
metabelian covers of $4$-manifolds.  A secondary obstruction theory results,
with vanishing criteria determined by first order choices. 
 
Our obstructions are Witt classes of intersection forms on the homology of
higher-order solvable covers, obtained from a sequence of new higher-order
linking pairings (see Section~\ref{sec:Alexander}). We define what we call {\it  rationally universal
$n$-solvable knot groups}, constructed from universal torsion modules,
which play roles analogous to $\Q/\Z$ in the torsion linking pairing on a
rational homology sphere, and to
$\Q(t)/\Q[t^{\pm 1}]$ in the classical Blanchfield pairing of a knot. 
Representations of the knot group into these groups are parametrized by elements
of the higher-order
Alexander modules. The key point is that if $K$   is slice  (or merely
$(n)$-solvable), then some predictable fraction of
these representations extends to the complement of the slice disk (or the
$(n)$-solution W).
The Witt classes of the intersection forms of these $4$-manifolds then
constitute invariants 
that vanish for slice knots (or merely $(n.5)$-solvable knots). 

For any {\it  fixed} knot and any {\it  fixed} $(n)$-solution $W$ one can show that a
signature vanishes by using certain solvable quotients
of $\pi_1(W)$, and not using the universal groups.  However a general
obstruction theory requires the introduction of these universal groups just as the
study of torsion linking pairings on all rational homology $3$-spheres requires the
introduction of $\Q/\Z$.
 
We first define the rationally universal solvable groups. The
metabelian group is a rational analogue of the group  used by
Letsche~\cite{Let}. Let
$ \G_0:=\Z$ and let $\KK_0$ be the quotient field of  
$\Z\G_0$. Consider a PID $\RR_0$ that lies in between $\Z\G_0$ and $\KK_0$. For
example, a good choice is $\Q[\mu^{\pm 1}]$ where $\mu$ generates $\G_0$.  Note
that
$\KK_0=\Q(\mu)$. For any choice of $\RR_0$, the abelian group $\KK_0/
\RR_0$ is a bimodule over $ \G_0$ via left (resp.\ right)   multiplication. We
choose the right multiplication to define the semi-direct product
\vglue-12pt
$$
 \G_1:=(\KK_0/ \RR_0)\sd \G_0 .
$$ This is our {\it  rationally universal} metabelian (or  
$(1)$-solvable) group for knots in $S^3$. Inductively, we obtain {\it  rationally
universal}
$(n+1)$-solvable groups by setting
$$
 \G_{n+1}:=(\KK_n/\RR_n) \sd  \G_n
$$ for certain PID's $\RR_n$ lying in between $\Z\G_n$ and its quotient field
$\KK_n$. To define the latter
 we show in Section~\ref{sec:Alexander} that the ring $\Z\G_n$ satisfies the so-called {\it  Ore condition} which is
necessary and sufficient to construct the (skew) quotient field $\KK_n$ exactly as in the commutative case.

Now let $M$ be the 0-framed surgery on a knot in $S^3$. We begin with a fixed
representation into $ \G_{0}$ that is normally just the abelianization isomorphism
$\pi_1(M)^{ab}
\cong  \G_0$. Consider $\AA_0:= H_1(M;
 \RR_0)$, the ordinary (rational) Alexander module. Denote its dual by
$$
\AA^{\#}_0:= \Hom_{ \RR_0}(\AA_0, \KK_0/  \RR_0).
$$ Then the Blanchfield form
$$ B\l_0: \AA_0\times\AA_0 \ra \KK_0/  \RR_0
$$ is nonsingular in the sense that it provides an isomorphism
$\AA_0\cong\AA^{\#}_0$. Using basic properties of the semi-direct product, we
show in Section~\ref{sec:Alexander} that there is a one-to-one-correspondence
$$
\AA^{\#}_0 \lr
\Rep^{\ast}_{ \G_0}(\pi_1 (M), \G_1).
$$ Here $\Rep^{\ast}_{ \G_n}(G, \G_{n+1})$ denotes the set of representations of
$G$ into $ \G_{n+1}$ that agree with some fixed representation into $
\G_n$, modulo conjugation by elements in the subgroup $\KK_n/\RR_n$. Hence
when $a_0\in\AA_0$ the Blanchfield form $B\l_0$ defines  an action of
$\pi_1(M)$ on 
$ \RR_1$ and we may define the next Alexander module
$\AA_1=\AA_1(a_0):=H_1(M; \RR_1)$. We prove that a nonsingular Blanchfield form
$$ B\l_1:\AA_1\xrightarrow{\cong}\AA^{\#}_1:=\Hom_{ \RR_1}(\AA_1, \KK_1/ \RR_1)
$$ exists and induces a one-to-one correspondence
$$
\AA_1 \lr \Rep^{\ast}_{ \G_1}(\pi_1(M), \G_2).
$$ Iterating this procedure leads to the $(n-1)$-st Alexander module
$$
\AA_{n-1}=\AA_{n-1} (a_0,a_1,\dots,a_{n-2}):=H_1(M;\RR_{n-1})
$$ together with the $(n-1)$-st Blanchfield form
$B\l_{n-1}:\AA_{n-1}\xrightarrow{\cong}\AA_{n-1}^{\#}$ and a one-to-one correspondence
$$
\AA_{n-1}\lr \Rep^{\ast}_{ \G_{n-1}}(\pi_1(M), \G_n).
$$

We show in Section~\ref{sec:Witt} that for an $(n)$-solvable knot there exist
choices
 $(a_0,a_1,\dots,a_{n-1})$ that correspond to a representation
$\phi_n:\pi_1(M)\to
\G_n$ which extends to a spin 4-manifold $W$ whose boundary is $M$. We then
observe that the intersection form on
$H_2(W; \KK_n)$ is nonsingular and represents an element $B_n=B_n(M,\phi_n)$ of
$L^0(
\KK_n)$ which is well-defined (independent of $W$) modulo the image of $L^0(\Z
\G_n)$. Here $L^0(R)$, $R$ a ring with involution, denotes the Witt group of
nonsingular hermitian forms on finitely generated free $R$-modules, modulo
metabolic forms.  

We can now formulate our obstruction theory for
$(h)$-solvable knots. A more general version, Theorem~\ref{invariant thm}, is
stated and proved in Section~\ref{sec:Witt}.

\proclaimtitle{A special case}
\specialnumber{4.6} \proclaim{Theorem} \label{algebra}   Let $K$ be a knot in $S^3$
with $0$\/{\rm -}\/surgery
$M$.
\begin{itemize}
\item[{\rm (0):}]  If $K$ is $(0)$-solvable then there is a well-defined obstruction
$B_0\in L^0( \KK_0)/\break i(L^0(\Z \G_0))$.
\item[{\rm (0.5):}]  If $K$ is $(0.5)$\/{\rm -}\/solvable then $B_0=0$.
\item[{\rm (1):}]   If $K$ is $(1)${\rm -}solvable then there exists a submodule
$P_0\subset\AA_0$ such that
$P^{\perp}_0=P_0$ and such that for each $p_0 \in P_0$ there is an obstruction
$B_1=B_1(p_0)\in L^0( \KK_1)/i(L^0(\Z \G_1))$.
\item[{\rm (1.5):}]  If $K$ is $(1.5)$\/{\rm -}\/solvable then there is a
$P_0$ as above such that for all $p_0\in P_0$ the obstruction $B_1$ vanishes.\\
$\vdots$
\item[{\rm (n):}]   If $K$ is $(n)$\/{\rm -}\/solvable then there exists
$P_0$ as above such that for all $p_0\in P_0$ there exists
$P_1=P_1(p_0)\subset\AA_1(p_0)$  with  
$P_1^{\perp}=P_1$ and such that for all $p_1\in P_1$ there exists
$P_2=P_2(p_0,p_1)\subset\AA_2(p_0,p_1)$ with
$P_2=P^{\perp}_2$ and such that $\dots$ there exists
$P_{n-1}=P_{n-1}(p_0,\dots,p_{n-2})$ with
$P_{n-1}=P^{\perp}_{n-1}${\rm ,} and such that any
$p_{n-1}\in P_{n-1}$  corresponds to a representation
$\phi_n(p_0,\dots,p_{n-1}):\pi_1(M)\to \G_n$ that extends to some bounding
$4$\/{\rm -}\/manifold and thus induces a class
$B_n=B_n(p_0,\dots,p_{n-1})  
\in L^0( \KK_n)/i(L^0(\Z \G_{n}))$. 
\item[{\rm (n.5):}]  If $K$ is $(n.5)$\/{\rm -}\/solvable then there is an inductive sequence 
$$ P_0,P_1(p_0), \dots, P_{n-1}(p_0,\dots,p_{n-2})
$$ as above such that $B_n=0$ for  all $p_{n-1} \in P_{n-1}$.
\end{itemize}

\endproclaim

 Note that the above obstructions depend only on the
$3$-manifold $M$. In a slightly imprecise way one can reformulate the integral
steps in the theorem as follows. (The imprecision only comes from   the fact that
we translate the conditions $P_i^{\perp}=P_i$ into   talking about ``one-half''
of the representations in question.) We try to count those representations of
$\pi_1(M)$ into
$ \G_n$   that extend to $\pi_1(W)$ for some $4$-manifold
$W$.
$$
\BoxedEPSF{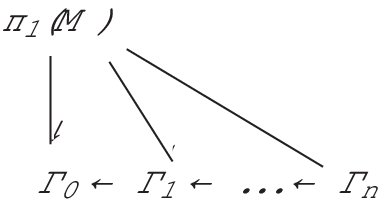 scaled 1000}
$$ 
If the knot $K$ is $(0)$-solvable, i.e.\ the Arf invariant   vanishes, then the
abelianization $\pi_1(M) \to  \G_0$ extends to a $4$-dimensional spin manifold
$W$. Then
 $B_0$ is defined. For $(0.5)$-solvable (or algebraically slice) knots this
invariant vanishes, giving
$P_0\subset\AA_0
\cong\Rep_{ \G_0}(\pi_1(M), \G_1)$. The corresponding representations to $
\G_1$ {\it  may not} extend over $W$. But if the knot $K$ is $(1)$-solvable via a
$4$-manifold $W$, then one-half of the   representations to $ \G_1$ {\it  do
extend} to $\pi_1(W)$.

For each such extension $p_0$ we form the next Alexander module $\AA_1(p_0)$,
which para\-me\-trizes representations into $ \G_2$, fixed over $ \G_1$, and
consider
$B_1\in L^0( \KK_1)$ (which depends on $p_0$). If $K$ is $(1.5)$-solvable, this
invariant vanishes and gives $P_1\subset\AA_1$. Again the corresponding
representations to $ \G_2$ might not extend to this $4$-manifold $W$. But if
$K$ is
$(2)$-solvable, then one {\it  quarter} of the representations to $ \G_2$ extend
to a
$(2)$-solution $W$. Continuing in this way, we get the following meta-statement:
\begin{quote}
 {\it  If $K$ is $(n)$-solvable via $W$ then $\frac{1}{2^n}$   of all
representations into $\G_n$ extend from $\pi_1(M)$ to $\pi_1(W)$.} 
\end{quote}

To be more precise, the following rather striking statement follows from
Lemma~\ref{2comp} and Proposition~\ref{rank coef}: For any slice knot for which
the degree of the Alexander polynomial is greater than
$2$ let $W$ be the complement of a slice disk for $K$. Then, for any n, at least
one $\G_n$-representation extends from $\pi_1(M)$ to $\pi_1(W)$. Moreover, this
representation is nontrivial in the sense that it does {\it  not} factor through
$\G_{n-1}$.

\vglue8pt 1.3. {\it $L^2$\/{\rm -}\/signatures.}
There remains the issue of detecting nontrivial classes in the $L$-theory of the
quotient fields $\KK$ of $\Z\G$. Our numerical invariants arise from
$L^2$-homology and von Neumann algebras (see Section~\ref{sec:example}).  We
construct an
$L^2$-signature
$$
\sigma_{\G}^{(2)} : L^0(\KK) \to \R
$$ by factoring through $L^0( \UU\G)$, where
$ \UU\G$ is the algebra of (unbounded) operators affiliated to the von Neumann
algebra
$ \NN\G$ of the group $\G$.  We show in Section~\ref{sec:L2} that this invariant
can be easily calculated in a large number of examples. The {\it  reduced}
$L^2$-signature, i.e.\ the difference of
$\sigma_{\G}^{(2)}$ and the ordinary signature,  turns out to be exactly what we
need to detect our obstructions $B_n$ from Theorem~\ref{invariant thm}.
 The fact that it does not depend on the choice of an
$(n)$-solution can be proved in three essentially different ways. Firstly, one
can show~\cite{Ma}, \cite{R} that the reduced $L^2$-signature of a $4k$-manifold with boundary
$M$ equals the reduced von Neumann 
$\eta$-invariant of the signature operator\break (associated to the regular
$\G$-cover of the $(4k-1)$-manifold $M$). This so-called von Neumann
$\rho$-invariant was introduced by J. Cheeger and M. Gromov \cite{ChG} who showed
in particular that it does not depend on a Riemannian metric on $M$ since it is a
difference of
$\eta$-invariants.  It follows that the reduced $L^2$-signature does not depend
on a bounding
$4$-manifold (which might not even exist) and can thus be viewed as a function of
$(M,\phi:\pi_1(M)\ra\G)$.
 
In the presence of a bounding $4$-manifold, the well-definedness of the invariant
can be deduced from Atiyah's $L^2$-index theorem
\cite{A}. This is even true in the topological category (see
Section~\ref{sec:L2}). There we also explain the third point of view, namely that
for groups $\G$ for which the analytic assembly map is onto, the reduced
$L^2$-signature actually vanishes on the image of $L^0(\Z\G)$ and thus clearly is
well-defined on our obstructions $B_n$ from Theorem~\ref{invariant thm}. By a
recent result of N. Higson and G. Kasparov \cite{HK} this applies in particular
to all torsion-free amenable groups (including our rationally universal solvable
groups). This last point of view is the strongest in the sense that it shows that
in order to define our obstructions one can equally well work with
$(n)$-solutions $W$ that are finite Poincar\'{e} $4$-complexes (rather than
topological
$4$-manifolds).

It seems that the invariants of Casson-Gordon should also be interpretable in
terms of
$\rho$-invariants (or signature defects) associated to {\it  finite-dimensional}
unitary representations of finite-index subgroups of $\pi_1(M)$
\cite{CG1}, \cite[p.~661]{KL}, \cite{Let}. J. Levine, M. Farber and W. Neumann have
also investigated finite dimensional $\rho$-invariants as applied to knot
concordance
\cite{L3}, \cite{N}, \cite{FL}. More recently C. Letsche used such
$\rho$-invariants together with a universal metabelian group to construct
concordance invariants
\cite{Let}. 

Since the invariants we employ are {\it  von Neumann} $\rho$-invariants,  they are
associated to the regular representation of our rationally universal solvable
groups on an {\it  infinite} dimensional Hilbert space.  These groups have to
allow homomorphisms from arbitrary knot (and slice) complement fundamental
groups, hence they naturally have to be huge and thus might not allow any
interesting finite dimensional representations at all.

The following is the result of applying Theorem~\ref{invariant thm} (just at the
level of obstructions to $(1.5)$-solvability) and the $L^2$-signature to the case
of genus one knots in homology spheres which should be compared to \cite[Th.\
4]{G2}. The proof, which will appear in another paper, is not difficult. It uses
the fact that in the simplest case of an $L^2$-signature for knots, namely where
one uses the abelianization homomorphism $\pi_1(M) \to \Z$, the real number
$\sigma_{\Z}^{(2)}(M)$ equals the integral over the
circle of the Levine signature function.

\proclaimtitle{[COT]}
\proclaim{Theorem} \label{genusone} 
Suppose $K$ is a $(1.5)$\/{\rm -}\/solvable
knot with a genus~one Seifert surface $F$. Suppose that the classical Alexander
polynomial of $K$ is nontrivial. Then there exists a homologically essential
simple closed curve $J$ on $F${\rm ,} with self\/{\rm -}\/linking zero, such that the integral over the
circle of the Levine signature function of $J$ {\rm (}\/viewed as a knot\/{\rm )} vanishes.
\endproclaim

1.4. {\it Paper outline and acknowledgments}. The paper is organized as follows:
Section~\ref{sec:PTFA} provides the necessary algebra to define the higher-order
Alexander modules and Blanchfield linking forms. In Section~\ref{sec:Alexander}
we construct our rationally universal solvable groups and investigate the
relationship between representations into them and higher-order Blanchfield
forms. We define our knot slicing obstruction theory in  Section~\ref{sec:Witt}.
Section~\ref{sec:L2} contains the proof that the
$L^2$-signature may be used to detect the $L$-theory classes of our obstructions.
In Section~\ref{sec:example}, we construct knots with vanishing Casson-Gordon
invariants that are  not topologically slice, proving our main Theorem~\ref{Q
slice}. Section~\ref{sec:surfaces} reviews intersection theory and defines
Whitney towers and gropes.  Section~\ref{sec:$H_1$-bordism} defines
$(h)$-solvability, and proves our theorems relating this filtration to gropes and
Whitney towers.  In Section~\ref{sec:CG} we prove Theorem~\ref{vanish}, showing
that Casson-Gordon invariants obstruct a second stage of Whitney disks. 
   
The authors are happy to thank Jim Davis and Ian Hambleton for interesting
conversations.  Wolfgang L\"uck, Holger Reich, Thomas Schick and Hans Wenzl answered
numerous questions on Section~\ref{sec:L2}. The heuristic argument concerning
Casson-Gordon invariants  and Casson-towers appears to be well-known. For the second
author, this argument was first explained by Shmuel Weinberger in $1985$ and he thanks
him for this insight. Moreover, we thank the Mathematical Sciences Research Institute in
Berkeley for providing both space and financial support during the 1996--97 academic
year, and the best possible environment for this project to take flight.

\section{Higher order Alexander modules and Blanchfield linking
forms}\label{sec:PTFA}

In this section we show that the classical Alexander module and Blanchfield
linking form associated to the infinite cyclic cover of the knot complement can
be extended to torsion modules and linking forms associated to {\it  any}
poly-torsion-free abelian covering space. We refer to these as {\it  higher-order
Alexander modules} and {\it  higher-order linking forms}. A forthcoming paper will
discuss these higher-order modules
 from the more  traditional viewpoint of Seifert\break surfaces [C].  

Consider a tower of regular covering spaces
$$ M_n\to M_{n-1}\to\dots\to M_1\to M_0=M
$$ such that each $M_{i+1}\to M_i$ has a torsion-free abelian group of deck
translations and each $M_i\to M$ is a regular cover. Then the group of deck
translations $\G$ of $M_n\to M$ is a {\it  poly-torsion-free abelian group\/} (see
below) and it is easy to see that such towers correspond precisely to certain
normal series for such a group. In this section we use such towers to generalize
the Alexander module. We will show that if
$\b_1(M)=1$ then $H_1(M_n;\Z)$ is a torsion $\Z\G$-module. 

\numbereddemo{Definition}\label{def:PTFA} A group $\G$ is {\it  poly\/{\rm -}\/torsion\/{\rm -}\/free abelian}
(PTFA)  if it admits a normal series $\left<1\right>=G_0\triangleleft
G_1\triangleleft\dots\triangleleft G_n=\G$ such that the factors $G_{i+1}/G_i$
are torsion-free abelian. (In the literature only a subnormal series is required.)
\enddemo

\numbereddemo{Example}\label{ex:PTFA} If $G$ is the fundamental group of a (classical) knot
exterior then $G/G^{(n)}$ is PTFA since the quotients of successive terms in the
derived series $G^{(i)}/G^{(i+1)}$ are torsion-free abelian  \cite{Str}. The
corresponding covering space is obtained by taking iterated universal abelian
covers.
\enddemo

\numbereddemo{{R}emark}\label{rem:PTFAext} If $A\triangleleft G$ is torsion-free abelian and
$G/A$ is PTFA then $G$ is PTFA. Any subgroup of a PTFA group is a PTFA group
(Lemma~2.4, p.~421 of
\cite{P}). Clearly any PTFA group is torsion-free and solvable (although the
converse is false!). The class of PTFA groups is quite large --- it contains all
torsion-free nilpotent groups ~\cite[Cor.~1.8]{Str}.
\enddemo

For us there are two especially important properties of PTFA groups, which we
state as propositions. These should be viewed as natural generalizations of
well-known properties of the free abelian group. The first is an algebraic
generalization of the fact that any infinite cyclic cover of a $2$-complex with
vanishing $H_2$ also has vanishing $H_2$.  It holds, more generally, for any locally indicable group $\G$.

\proclaimtitle{[Str, p.~305]}
\proclaim{Proposition}\label{injective}  Suppose $\G$ is a
{\rm PTFA} group and
$R$ is a commutative ring. Any map between projective right
$R\G$\/{\rm -}\/modules whose image under the functor
$-\otimes_{R\G}R$ is injective{\rm ,} is itself injective.
\endproclaim 

The second important property is that $\Z\G$ has a  (skew) quotient field. Recall
that if $A$ is a {\it  commutative\/} ring and $S$ is a subset closed under
multiplication, one can construct the {\it  ring of fractions\/} $AS^{-1}$ of
elements
$as^{-1}$ which add and multiply like normal fractions. If
$S=A-\{0\}$ and $A$ has no zero divisors, then $AS^{-1}$ is called the {\it  quotient field\/} of $A$. However, if $A$ is {\it  noncommutative\/} then
$AS^{-1}$ does not always exist (and $AS^{-1}$ is not {\it a priori} isomorphic to
$S^{-1}\! A$). It is known that if $S$ is a {\it  right divisor set\/} then
$AS^{-1}$ exists (\cite[p.~146]{P} or ~\cite[p.~52]{Ste}). If $A$ has no zero
divisors and
$S=A-\{0\}$ is a right divisor set  then $A$ is called an {\it  Ore domain\/}. In
this case $AS^{-1}$ is a skew field, called the {\it  classical right ring of
quotients\/} of
$A$. We will often refer to this merely as the {\it  quotient field } of $A$ . A
good reference for noncommutative rings of fractions is Chapter~2 of
\cite{Ste}. In this paper we will always use {\it  right\/} rings of fractions.
The following holds more generally for any torsion-free amenable group.

\proclaim{Proposition}\label{qG} If $\G$ is {\rm PTFA} then $\Q\G$ is a right {\rm (}\/and left\/{\rm )} Ore
domain\/{\rm ;} i.e.\ $\Q\G$ embeds in its classical right ring of quotients $\KK${\rm ,} which
is a skew field.
\endproclaim

\demo{Proof}  For the fact (due to A.A. Bovdi) that $\Z\G$ has no
 zero divisors see \cite[pp.~591-592]{P} or \cite[p.~315]{Str}.  As we have
remarked, any PTFA group is solvable. It is a result of J.~Lewin~\cite{Le} that
for solvable groups such that $\Q\G$ has no zero divisors, $\Q\G$ is an Ore
domain   (see
 Lemma~3.6 iii, p.~611 of \cite{P}).
\enddemo

If $\RR$ is an integral domain then a right $\RR$-module
$\AA$ is said to be a {\it  torsion module\/} if, for each $a\in\AA$, there exists
some nonzero $r\in \RR$ such that $ar=0$. If $\RR$ is an Ore domain then $\AA$
is a torsion module if and only if $\AA\otimes_\RR\KK=0$ where $\KK$ is the
quotient field of $\RR$. 
\cite[II Cor.~3.3]{Ste}. In general,  the set of torsion elements of $A$ is
a submodule.

\numbereddemo{{R}emark}\label{facts} We shall need the following elementary facts about the
right skew field of quotients $\KK$. It is naturally a right
$\KK$-module and is a $\Z\G$-bimodule.
\begin{description}
\item[Fact {\rm 1:}]  $\KK$ is flat as a left $\Z\G$-module; i.e.\ $\cdot\otimes_{\Z\G}\KK$ is
exact~\cite[Prop.~II.3.5]{Ste}.
\item[Fact {\rm 2:}] Every module over $\KK$ is a free
module~\cite[Prop.~I.2.3]{Ste} and  such modules have a well defined rank
$\rk$ which is additive on short exact sequences \cite[p.~48]{Co1}.
\end{description}
\enddemo

{\it Homology of {\rm PTFA} covering spaces.}
Suppose $X$ has the homotopy type of a connected CW-complex,
$\G$ is a group and $\phi:\pi_1(X)\ra\G$ is a homomorphism. Let $X_\G$ denote the
{\it  regular $\G$-cover of $X$ associated to $\phi$} (by pulling back the
universal cover of
$B\G$). Note that if $\pi=\mbox{image}(\phi)$ then $X_\G$ is a disjoint union of
$\G/\pi$ copies of the connected cover
$X_\pi$ (where $\pi_1(X_\pi)\cong\Ker (\phi$)). Fixing a certain convention
(which will become clear in Section~\ref{sec:example}), $X_\G$ becomes a right
$\G$-set. For simplicity, the following are stated for the ring $\Z$, but also
hold for $\Q$ and
$\C$. Let $\MM$ be a
$\Z\G$-bimodule (for us usually $\Z\G$, $\KK$, or a ring $\RR$ such that
$\Z\G\subset
\RR\subset\KK$, or $\KK/\RR$). The following are often called the equivariant
homology and cohomology of X.

\numbereddemo{Definition}\label{def:XR} Given $X$, $\phi$, $\MM$ as above, let
$$ H_{\ast}(X;\MM)\equiv H_{\ast}(C_{\ast}(X_\G;\Z)\otimes_{\Z\G}
\MM)
$$ as a right $\Z\G$ module, and $H^{\ast}(X;\MM)\equiv
H_{\ast}\left(\Hom_{\Z\G}(C_{\ast}(X_\G;\Z),\MM)\right)$ as a left
$\Z\G$-module.
\enddemo

But these are well-known to be  isomorphic (respectively) to the homology (and
cohomology) of $X$ with coefficient system induced by
$\phi$ (see Theorems~VI 3.4 and 3.4$^*$ of \cite{W}).

\numbereddemo{{R}emark}\label{postlore} 
 1.   Note that $H_{\ast}(X;\Z\G)$ as in Definition~\ref{def:XR} is merely
$H_{\ast}(X_\G;\Z)$ as a right $\Z\G$-module. Thus if $\MM$ is flat as a left
$\Z\G$-module then $H_{\ast}(X;\MM)\cong H_{\ast}(X_\G;\Z)\otimes_{\Z\G}\MM$.
Hence the homology groups we discuss have an interpretation as homology of
$\G$-covering spaces. However the cohomology
$H^{\ast}(X;\Z\G)$ does not have such a direct interpretation, although it can be
interpreted as cohomology of $X_\G$ with compact supports (see, for
instance,~\cite[p.~5--6]{Hi}.)
\vglue4pt 
2.   Recall that if $X$ is a compact, oriented $n$-manifold then by Poincar\'e
duality
$H_p(X;\MM)$ is isomorphic to $H^{n-p}(X,\partial X;\MM)$ which is made into a
right
$\Z\G$-module using the involution on this ring \cite{Wa}.

\vglue4pt 3. We also have a universal coefficient spectral sequence (UCSS) as in
\cite[Th.~2.3]{L2}. If $R$ and $S$ are rings with unit, $C$ a free right
chain complex over $R$ and
$\MM$ an $(R-S)$ bimodule, there is a convergent spectral sequence
$$ E_2^{p,q}\cong\Ext^q_R(H_p(C),\MM)\Longrightarrow H^{\ast}(C;\MM)
$$ of left $S$-modules (with differential $d^r$ of degree
$(1-r,r)$). Note in particular that the spectral sequence collapses when
$R=S=\KK$ is the (skew) field of quotients since
$\Ext^n_{\Z\G}(M,\KK)\cong\Ext^n_\KK(M\otimes_{\Z\G}\KK,\KK)$ by change of rings
\cite[Prop. 12.2]{HS}, and the latter is zero if $n\geq 1$ since all
$\KK$-modules are free. Hence
$$ H^n(X;\KK)\cong\Hom_\KK(H_n(X;\KK),\KK).
$$ More generally it collapses when $R=S$ is a (noncommutative) principal ideal domain.
\enddemo

Suppose that $\G$ is a PTFA group and $\KK$ is its (skew) field of quotients. We
now investigate $H_0$, $H_1$ and
$H_2$ of spaces with coefficients in $\Z\G$ or $\KK$. First we show that
$H_0(X;\Z\G)$ is a torsion module.

\proclaim{Proposition}\label{acyclic} Given $X$, $\phi$ as in Definition~{\rm \ref{def:XR},}
suppose a ring homomorphism
$\psi:\Z\G\ra\RR$ defines $\RR$ as a $\Z\G$\/{\rm -}\/bimodule. Suppose some element of the
augmentation ideal of $\Z[\pi_1(X)]$ is invertible {\rm (}\/under $\psi\circ\phi${\rm )} in
$\RR$. Then
$H_0(X;\RR)=0$. In particular{\rm ,} if
$\phi:\pi_1(X)\ra\G$ is a nontrivial coefficient system then
$H_0(X;\KK)=0$.
\endproclaim

\demo{Proof}  By \cite[p. 275]{W} and \cite[p.34]{Br}, $H_0(X;\RR)$ is isomorphic
to the cofixed set $\RR/\RR I$ where $I$ is the augmentation ideal of
$\Z\pi_1(X)$ acting via
$\psi\circ\phi$.
\enddemo

The following proposition is enlightening, although in low dimensions its use can
be avoided by short ad hoc arguments. Here $\Q$ is a $\Z\G$ module via the
composition
$\Z\G
\xrightarrow{\e} \Z \to \Q$ where $\e$ is the augmentation homomorphism.

\proclaim{Proposition}\label{n-conn}  {\rm a)} If $C_{\ast}$ is a nonnegative $\Q\G$ chain complex which is
finitely generated and free in dimensions $0\le i\le n$ such that
$H_i(C_{\ast}\otimes_{\Q\G}\Q)=0$ for $0\le i\le n${\rm ,} then
$H_i(C_{\ast}\otimes_{\Q\G}\KK)=0$ for $0\le i\le n$. 
\vglue4pt {\rm b)} If $f:Y\to X$ is a continuous map{\rm ,} between {\rm CW} complexes with finite
$n$\/{\rm -}\/skeleton which is $n$\/{\rm -}\/connected on rational homology{\rm ,} and
$\phi:\pi_1(X)\ra\G$ is a coefficient system{\rm ,} then $f$ is
$n$\/{\rm -}\/connected on homology with $\KK$\/{\rm -}\/coefficients.
\endproclaim

\demo{Proof}   Let $\e:\Q\G\to\Q$ be the augmentation and $\e(C_{\ast})$ denote
$C_{\ast}\otimes_{\Q\G}\Q$. Since $\e(C_{\ast})$ is acyclic up to dimension~$n$,
there is a $``$partial'' chain homotopy
$$
\{h_i:\e(C_{\ast})_i\to\e(C_{\ast})_{i+1}\mid0\le i\le n\}
$$ between the identity and the zero chain homomorphisms. By this we mean that
$\partial h_i+h_{i-1}\partial =\id$ for $0\le i\le n$.

Since $C_i \xrightarrow{\e} \e(C_i)$ is surjective, for any basis element
$\sigma$ of
$C_i$ we can choose an element, denoted $\widetilde h_i(\sigma)$, such that
$\e\circ\widetilde h_i(\sigma)=h_i(\e(\sigma))$. Since
$C_{\ast}$ is free, in this manner $h$ can be lifted to a partial chain homotopy
$\{\widetilde h_i\mid0\le i\le n\}$ on
$C_{\ast}$ between some ``partial'' chain map $\{f_i\mid0\le i\le n\}$ and the
zero map. Moreover $\e(f_i)$ is the identity map on $\e(C_{\ast})_i$, and in
particular, is injective. Thus, by Proposition~\ref{injective}, $f_i$ is
injective for each $i$. Consequently, $\widetilde h_i\otimes\id$ is a partial
chain homotopy on
$C_{\ast}\otimes_{\Q\G}\KK$ between the zero map and the partial chain map
$\{f_i\otimes\id\}$, such that $f_i\otimes\id$ is injective (since $\KK$ is flat
over
$\Q\G$). Any monomorphism between finitely generated, free $\KK$-modules of the
same rank is necessarily an isomorphism. Therefore a partial chain map exists
which is an inverse to $f\otimes\id$. It follows that
$C_{\ast}\otimes_{\Q\G}\KK$ is acyclic up to and including dimension~$n$.

The second statement follows from applying this to the relative cellular chain
complex associated to the mapping cylinder of $f$.
\enddemo

\proclaim{Proposition}\label{torsion} Suppose $X$ is a {\rm CW-}\/complex such that $\pi_1(X)$ is
finitely generated{\rm ,} and
\mbox{$\phi:\pi_1(X)\ra\G$} is a nontrivial coefficient system. Then
$$
\rk H_1(X;\KK)\le \b_1(X)-1.
$$ In particular{\rm ,} if $\b_1(X)=1$ then
$H_1(X;\KK)=0${\rm ;} that is{\rm ,} $H_1(X;\Z\G)$ is a $\Z\G$ torsion module.
\endproclaim

\demo{Proof}  Let $Y$ be a wedge of
$\b_1(X)$ circles. Choose $f:Y\to X$ which is $1$-connected on rational homology.
Applying Proposition~\ref{n-conn}, one sees that
$f_{\ast}:H_1(Y;\KK)\ra H_1(X;\KK)$ is surjective. We claim that $\phi\circ f_\ast$ is nontrivial on $\pi_1(Y)$.  Suppose not.  Let $G$ denote the image of $\phi$.  Note that
if $\{x_i\}$ generates $\pi_1(Y)$ then $\{\phi\circ  f_\ast(x_i)\}$ generates $G/G^{(1)}\otimes \Q$, which, under our
supposition, would imply that the {\it nontrivial} PTFA group $G$ had a finite abelianization.  But one sees from Definition
2.1 that the abelianization of a PTFA group has a quotient $(G_n/G_{n-1}$ in the language of 2.1) that is a nontrivial
torsion-free abelian  group and therefore must contain an element of infinite order.  This contradiction implies $\phi\circ
f_\ast$ is nontrivial.
 Finally Lemma~\ref{2comp}
below shows that
$$
\rk H_1(Y;\KK)=\b_1(Y)-1.
$$ The claimed inequality follows. If $H_1(X;\KK)=0$ then
$H_1(X;\Z\G)$ is a $\Z\G$ torsion module by Remark~\ref{facts}.1 and
\cite[II Cor.~3.3]{Ste}.
\enddemo

\proclaim{Lemma}\label{2comp}\hskip-8pt Suppose $Y$ is a finite connected
$2$\/{\rm -}\/complex with $H_2(Y;\Z)= 0$ and $\phi:\pi_1(Y)\ra\G$ is nontrivial. Then
$H_2(Y;\Z\G)= H_2(Y;\KK)=0$ and
$\rk H_1(Y;\KK)=\b_1(Y)-1$.
\endproclaim

\demo{Proof}  Let 
$$
\CC= \left( 0 \ra C_2 \xrightarrow{\partial_2} C_1 
\xrightarrow{\partial_1} C_0 \ra 0\right)
$$ be the free $\Z\G$ chain complex for the cellular decomposition of $Y_\G$ (the
$\G$ cover of $Y$) obtained by lifting cells of $Y$. Since $H_2(Y;\Z)=0$,
$\partial_2\otimes_{\Z\G}\id$ is injective, which implies, by
Proposition~\ref{injective}, that $\partial_2$ itself is injective. Thus
$H_2(Y;\Z\G)=0$ by Remark~\ref{postlore}.1 and $H_2(Y;\KK)=0$ by
Remark~\ref{facts}.1. Since $\phi$ is nontrivial, Proposition~\ref{acyclic}
implies that
$H_0(Y;\KK)=0$. Since the $C_i$ are finitely generated free modules, the Euler
characteristic of
$\CC\otimes\KK$ equals the Euler characteristic of
$\CC\otimes\Q$ (by Remark~\ref{facts}.2) and the result follows.
\enddemo

It is interesting to note that \ref{torsion} and \ref{2comp} are false without the finiteness assumptions (see
Section 3 of~\cite{C}.)

Thus we have shown that the definition of the classical Alexander module, i.e.\ the torsion module
 associated to the first homology of the infinite cyclic cover of the knot
complement, can be extended to {\it  higher\/{\rm -}\/order Alexander modules} which are
$\Z\G$ torsion modules
$\AA = H_1(M;\Z\G)$  associated to {\it  arbitrary\/} PTFA covering spaces.
Indeed, by Proposition~\ref{torsion}, this is true for any
$3$-manifold with $\b_1(M)=1$, such as zero surgery on the knot or a prime-power
cyclic cover of $S^3-K$. In this paper we will work with the zero surgery.
 
Furthermore, we will now show that the Blanchfield linking form associated to the
infinite cyclic cover generalizes to linking forms on these higher-order
Alexander modules. Under some mild restrictions, we can get a nonsingular
linking form in the sense of A.~Ranicki. Recall from \cite[p. 181--223]{Ra2} that
$(\AA,\lambda)$ is a  {\it  symmetric linking form\/} if $\AA$ is a torsion
$\RR$-module of (projective) homological dimension~$1$ (i.e.\ $\AA$ admits a
finitely-generated projective resolution of length $1$) and
$$
\lambda :\AA\ra \Hom_\RR(\AA,\KK/\RR)\equiv \AA^\#
$$ is an $\RR$-module map such that $\lambda (x)(y)=\overline{\lambda (y)(x)}$
(here
$\KK$ is the field of fractions of $\RR$ and $\AA^\#$ is made into a right
$\RR$-module using the involution of $\RR$.). The linking form is {\it  nonsingular\/} if
$\lambda$ is an isomorphism. If $\RR$ is an integral domain then $\RR$ is a {\it  {\rm (}\/right\/{\rm )} principal ideal
domain\/} (PID) if every right ideal is principal.

\proclaim{Theorem}\label{higher Bforms} Suppose $M$ is a closed{\rm ,} oriented{\rm ,} connected
$3$\/{\rm -}\/manifold with $\b_1(M)=1$ and
$\phi:\pi_1(M)\ra\G$ a nontrivial {\rm PTFA} coefficient system. Suppose $\RR$ is a
ring such that $\Z\G\sbq \RR\sbq\KK$. Then there is a symmetric linking form 
$$ B\l:H_1(M;\RR)\ra H_1(M;\RR)^\#
$$ defined on the {\rm  higher-order Alexander module} $\AA:=H_1(M;\RR)$.
 If either $\RR$ is a {\rm PID,} or some element of the augmentation ideal of
$\Z\pi_1(M)$ is sent {\rm (}\/under~$\phi$\/{\rm )} to an invertible element of $\RR${\rm ,} then
$B\l$ is nonsingular.
\endproclaim

\demo{Proof}  Note that $\AA$ is a torsion $\RR$-module by
Proposition~\ref{torsion}, since
$\KK$ is also the quotient field of the Ore domain $\RR$. Define $B\l$ as the
composition of the Poincar\'e duality isomorphism to $H^2(M;\RR)$, the inverse of
the Bockstein to
$H^1(M;\KK/\RR)$, and the usual Kronecker evaluation map to $\AA^\#$. The
Bockstein 
$$  B:H^1(M;\KK/\RR)\ra H^2(M;\RR)
$$  associated to the short exact sequence
$$  0\ra \RR\ra\KK\ra\KK/\RR\ra0
$$  is an isomorphism since
$H^2(M;\KK)\cong H_1(M;\KK)=0$ by Proposition~\ref{torsion}, and
$H^1(M;\KK)=0$ by Proposition~\ref{torsion} and Remark~\ref{postlore}.3. Under
the second hypothesis on $\RR$, the Kronecker evaluation map 
$$ H^1(M;\KK/\RR)\ra\Hom_\RR\left(H_1(M;\RR),\KK/\RR\right)
$$ is an isomorphism by the UCSS since $H_0(M;\RR)=0$ (see
Remark~\ref{postlore}.3  and Proposition~\ref{acyclic}).  If
$\RR$ is a PID then $\KK/\RR$ is an injective $\RR$-module since it is clearly
divisible~\cite[I Prop.~6.10]{Ste}. Thus
$$
\Ext^i_\RR\left(H_0(M;\RR),\KK/\RR\right)=0
$$ for $i>0$ and therefore the Kronecker map is an isomorphism.

We need to show that $\AA$ has homological dimension one and is finitely
generated. This is immediate if $\RR$ is a PID \cite[p.\ 22]{Ste}. Since, in
this paper we shall only need this special case we omit the proof of the general
case. 

We also need to show that $B\l$ is $``$conjugate symmetric". The diagram below
commutes up to a sign (see, for example,~\cite[p.~410]{M}), where $B'$ is the
homology Bockstein
$$
\begin{array}{ccc}
{H_2(M; \KK/\RR)}  &\hskip-24pt\stackrel{B'}{\longrightarrow}&{H_1(M;
\RR)}
\\
{\cong}\Big\downarrow{\rm P.D.} &&{\cong}\Big\downarrow{\rm P.D.}  \\
{H^1(M; \KK/\RR)} &\hskip-24pt\stackrel{B}{\longrightarrow}&{H^2(M;\RR)}\\
\Big\downarrow{\kappa}&&\\
{\Hom_{\RR}(H_1(M; \RR), \KK/\RR)}
\end{array}
$$ and the two vertical homomorphisms are Poincar\'e duality. Thus our map
$B\l$ agrees with that obtained by going counter-clockwise around the square and
thus agrees with the Blanchfield form defined by J.~Duval in a noncommutative
setting~\cite[p.~623--624]{D}. The argument given there for symmetry is in
sufficient generality to cover the present situation and the reader is referred
to it. 
\enddemo 

The implications of the following for the higher-order Alexander polynomials of
slice knots will be discussed in a forthcoming paper. This is the noncommutative
analogue of the result that the classical Alexander polynomial of a slice knot
factors as a product
$f(t)f(t^{-1})$.

\proclaim{Lemma}\label{Bl iso} If $\AA$ is a generalized Alexander module {\rm (}\/as in
Theorem~{\rm \ref{higher Bforms})} which admits a submodule $P$ such that $P=P^\bot${\rm ,}
then the map $h:P\ra(\AA/P)^\#${\rm ,} given by $p
\mapsto B\l(p,\cdot)${\rm ,} is an isomorphism.
\endproclaim

\demo{Proof}  Since the Blanchfield form is nonsingular by Theorem~\ref{higher
Bforms}, $h(p)$ is actually a monomorphism if $p\neq0$ and so $h$ is certainly
injective. Since $B\l:\AA\to\AA^\#$ was shown to be an isomorphism, it is easy to
see that $h$ is onto when
$P=P^\bot$.
\enddemo

\vglue-12pt

 \section{Higher-order linking forms and solvable representations\\ of the knot
group}\label{sec:Alexander}
\vglue-4pt

We now define and restrict our attention to certain families
$\G_0$, $\G_1,\dots,\G_n$ of PTFA groups that are constructed as semi-direct
products, inductively, beginning with
$\G_0\equiv\Z$, and defining $\G_n=A_{n-1}\sd\G_{n-1}$ for certain ``universal''
torsion $\Z\G_{n-1}$ modules $A_{n-1}$. We then show that if coefficient systems
$\phi_i:\pi_1(M)\ra\G_i$, $i<n$, are defined, giving rise to the higher-order
Alexander modules $\AA_0,\dots,\AA_{n-1}$,  then any nonzero choice
$x_{n-1}\in\AA_{n-1}$ corresponds to a nontrivial extension of $\phi_{n-1}$ to
$\phi_n:\pi_1(M)\ra\G_n$. This coefficient system is then used to define the
$n^{\rm{th}}$ Alexander module
$\AA_n(x_{n-1})$. Thus, if the ordinary Alexander module $\AA_0$ of a knot $K$ is
nontrivial, then there exist nontrivial $\G_1$-coefficient systems. This allows
for the definition of $\AA_1$, and if this module is nontrivial there exist
nontrivial
$\G_2$-coefficient systems. In this way, higher Alexander modules and actual
coefficient systems are constructed inductively from choices of elements of the
lesser modules. Naively stated: if $H_1$ of a covering space
$\widetilde M$ of $M$ is not zero then $\widetilde M$ itself possesses a
nontrivial abelian cover.

We close the section with a crucial result concerning when such coefficient
systems extend to bounding $4$-manifolds.

\vglue8pt {\it Families of universal {\rm PTFA} groups}.
We now inductively define families 
\mbox{$\{\G_n\mid n\geq 0\}$} of PTFA groups. These groups $\G_n$ are
``universal'' in the sense that the fundamental
group of any knot complement with nontrivial classical Alexander polynomial admits
nontrivial
$\G_n$-representations, a nontrivial fraction of which extend to the fundamental group of the complement of
a slice disk for the knot. These are the groups we shall use to construct our
knot slicing obstructions. Our approach elaborates work of Letsche who first
used an analogue of the group
$\G^U_1$~\cite{Let}.
 
Let $\G_0=\Z$, generated by $\mu$. Let $\KK_0=\Q(\mu)$ be the quotient field of
$\Q\G_0$ with the involution defined by
$\mu\to\mu^{-1}$. Choose a ring $\RR_0$ such that 
$\Q\G_0\subset \RR_0\subset\KK_0$. Note that $\KK_0/\RR_0$ is a
$\Z\G_0$-bimodule. Choose the right multiplication and define $\G_1$ as the
semidirect product $\KK_0/\RR_0\rtimes\G_0$. Note that if, for example,
$\RR_0=\Q[\mu^{\pm1}]=\Q\G_0$ then $\KK_0/\RR_0$ is a torsion
$\Q\G_0$ module that is, in fact, a direct limit of all cyclic torsion
$\Q\G_0$ modules.

In general, assuming $\G_{n-1}$ is defined (a PTFA group), let
$\KK_{n-1}$ be the quotient field of $\Q\G_{n-1}$ (by Proposition~\ref{qG}).
Choose any ring $\RR_{n-1}$ such that $\Q\G_{n-1}\subset
\RR_{n-1}\subset\KK_{n-1}$. Consider
$\KK_{n-1}/\RR_{n-1}$ as a right $\Z\G_{n-1}$-module and define
$\G_n$ as the semi-direct product $\G_n\equiv
(\KK_{n-1}/\RR_{n-1})\rtimes\G_{n-1}$. Since $\KK_{n-1}/\RR_{n-1}$ is a
$\Q$-module, it is torsion-free abelian. Thus $\G_n$ is PTFA by
Remark~\ref{rem:PTFAext}. We have the epimorphisms $\G_n
\stackrel{\pi}{\thra}\G_{n-1}$, and canonical splittings
$s_n:\G_{n-1}\to\G_n$. The family of groups depends on the choices for $\RR_i$.
The larger $\RR_i$ is, the more elements of $\Z\G_i$ will be invertible in
$\RR_i$; hence the more (torsion) elements of $H_1(M;\Z\G_i)$ will die in
$H_1(M;\RR_i)$; hence the more information will be potentially lost.  However,
in Proposition~\ref{acyclic} and Theorem~\ref{higher Bforms} we saw that it is
useful to have $\RR_i\neq\Q\G_i$ if
$i>0$ because it often ensures nonsingularity in higher-order Blanchfield forms.

For the final result of this section, concerning when coefficient systems extend
to bounding
$4$-manifolds, we find it necessary to introduce a rather severe (and hopefully
unnecessary) simplification: we take our Alexander modules (as in \ref{higher
Bforms}) to have coefficients in certain principal ideal domains
$\RR_0\dots,\RR_{n-1}$ where $\Z\G_i\sbq \RR_i\sbq\KK_i$. In some cases this can
have the unfortunate effect of completely killing
$H_1(M;\Z\G_i)$, which means that no interesting higher modules can be
constructed by the procedure below. However for most knots this does not happen.
Because of the importance, in this paper, of the family of groups corresponding
to these particular
$\RR_i$, we give it
 a specific notation:

\numbereddemo{Definition}\label{def:univ gr}  The family of {\it  rationally universal groups\/}
$\{\gu\}$ is defined inductively as above with $\gop=\Z$, $\RR^U_0=\Q[\mu^{\pm
1}]$, for $n\geq 0$,  
$$ S_n=\Q[\gu,\gu]-\{0\},\qquad \RR^U_n=(\Q\gu)S^{-1}_n
$$ and 
$$
\G^U_{n+1}=\KK_n/\RR^U_n\rtimes\gu.
$$ Here $\KK_n$ is the right ring of quotients of $\gu$.
\enddemo

Observe that this is quite a drastic localization. To form
$\RR_n$ we have inverted all the nonzero elements of the rational group ring of
the commutator subgroup of $\G^U_n$.

Note that $[\gu,\gu]$ is PTFA by Remark~\ref{rem:PTFAext} so that
$\Q[\gu,\gu]$ is an Ore domain. Therefore $S_n$ (above) is a right divisor set of
$\Q\gu$ by Chapter 13, Lemma~3.5 
of \cite[p.~609]{P}. One easily shows that $\gu$ is
$(n)$-solvable.

We will now show that the rings $\RR^U_n$ of Definition~\ref{def:univ gr} are in fact {\it  skew Laurent polynomial rings\/} which are
(noncommutative) principal right (and left) ideal domains by~\cite[2.1.1
p.~49]{Co2} generalizing the case $n=0$ where $\RR^U_0=\Q[\mu^{\pm1}]$. If
$\K$ is a skew field, $\a$ is an automorphism of
$\K$ and $\mu$ is an indeterminate, the {\it  skew {\rm (}\/Laurent\/{\rm )} polynomial ring in
$\mu$ over $\K$ associated with $\a$\/}, denoted $\K[\mu^{\pm1}]$, is the ring
consisting of all expressions
$$ f=\mu^{-m}a_{-m}+\dots+a_0+\mu a_1+\mu^2a_2+\dots+\mu^na_n
$$ where $a_i\in\K$, under ``coordinate-wise addition'' and multiplication
defined by the usual multiplication for polynomials and the rule
$a\mu=\mu\a(a)$~\cite[p.~54]{Co1}. The form above for any element $f$ is unique
\cite[p.~49]{Co2}. Note also that (for $a_{-m}$ and $a_n$ nonzero), the
nonnegative function $\deg f=n+m$ is additive under multiplication of
polynomials.

Now if $\G$ is a PTFA group and $G$ is a normal subgroup such that $\G/G\cong\Z$ is
generated by $\mu\in\G$, there is an automorphism of
$G$ given by $a\mapsto\mu^{-1}a\mu\equiv a^\mu$. It is a rather tedious
calculation to show that the abelianizations of our
$\G^U_i$ {\it  are in fact} $\Z$. Thus the G which is relevant for these cases is
actually the commutator subgroup. Since this fact is not crucial, we do not
include it.  In any case, there are other situations where
 one needs the extra generality of the following argument. Continuing, this
automorphism extends to a ring automorphism of
${\Q}G$ and hence, to one of $\K$, the quotient field of $\Q G$. Let $S=\Q
G-\{0\}$ and
$\RR=(\Q\G)S^{-1}$.

\proclaim{Proposition}\label{extendiso} The embedding
$g:\Q G\ra\K$ extends to an isomorphism
$\RR\ra\K[\mu^{\pm1}]$.
\endproclaim

\demo{Proof}  As an additive group, $\Q\G$ is isomorphic to
$\bigoplus^\infty_{i=-\infty}\Q G$ since the cosets of $G$ partition $\G$. But
$\K[\mu^{\pm1}]$, as a group, is isomorphic to a countable direct sum of copies of
$\K$. Therefore $g$ extends in an obvious way to an additive group homomorphism
$g:\Q\G\ra\K[\mu^{\pm1}]$ such that $g(\mu^ia_i)=\mu^ig(a_i)$ for
$a_i\in\Q G$. Since the automorphism
$a\mapsto a^\mu$ defining $\K[\mu^{\pm1}]$ agrees with conjugation in $\G$, this
map is a ring homomorphism. Clearly the nonzero elements of $\Q G$ are sent to
invertible elements. Moreover, any element of $\K[\mu^{\pm1}]$ is of the form
$\left(\Sigma\mu^ig(a_i)\right)s^{-1}$ where
$a_i\in\Q G$ and $s\in S$. This establishes that
$(\Q\G)S^{-1}\cong\K[\mu^{\pm1}]$, \cite[p.~50]{Ste}.
\enddemo

\proclaim{{C}orollary}\label{Knmu} For each $n\geq 0$ the rings $\RR^U_n$ of
Definition~{\rm \ref{def:univ gr}} are left and right principal ideal domains{\rm ,} denoted
$\K_n[\mu^{\pm1}]${\rm ,} where $\K_n$ is the right ring of quotients of
$\Z[\gu,\gu]$.
\endproclaim

\numbereddemo{{R}emark}\label{ex:Gnu} Suppose $\gu$ is one of the rationally universal groups
defined by Definition~\ref{def:univ gr}. Then, if $\phi$ is nontrivial on
$[\pi_1(X),\pi_1(X)]$, Proposition~\ref{acyclic} applies and $H_0(X;\RR^U_n)=0$ if
$n>0$. However, beware: $H_0(X,\RR^U_0)$ is certainly {\it  not\/} zero.
\enddemo

Suppose $M$ is a closed $3$-manifold with $\b_1(M)=1$. A choice of a generator of
$H_1(M;\Z)$ modulo torsion induces an epimorphism
$\phi_0:\pi_1(M,m_0)\to\G_0=\Z$. In case $M$ is an oriented knot complement this
choice is usually made using the knot orientation. Let $\AA_0\equiv H_1(M;\RR_0)$
be the {\it  rational Alexander module\/}, and suppose (inductively) that we are
given
$\phi_{n-1}:\pi_1(M)\ra\G_{n-1}$. Then we can define the {\it  higher-order
Alexander module\/}
$\AA_{n-1}\equiv H_1(M;\RR_{n-1})$, using the $\Z\G_{n-1}$ local coefficients
induced by
$\phi_{n-1}$. Varying $\phi_{n-1}$ by an inner automorphism of
$\G_{n-1}$ changes $H_1(M;\RR_{n-1})$ by an isomorphism induced by the
conjugating element. Let
${\rm Rep}_{\G_{n-1}}\left(\pi_1(M),\G_n\right)$ denote the set of homomorphisms
from
$\pi_1(M,m_0)$ to $\G_n$ which agree with $\phi_{n-1}$ after composition with the
projection
$\G_n\ra\G_{n-1}$. 

Recall that $\KK_{n-1}/\RR_{n-1}$ is a right $\Z\G_{n-1}$ module and hence
becomes a right $\Z\pi_1(M)$ module via
$\phi_{n-1}$. By a universal property of semi-direct products
\cite[VI Prop.~5.3]{HS}, there is a one-to-one correspondence between
${\rm Rep}_{\G_{n-1}}\left(\pi_1(M),\G_n\right)$ and the set of derivations
$d:\pi_1(M)\ra\KK_{n-1}/\RR_{n-1}$. One   checks that varying by a principal
derivation corresponds to varying the representation by a
$\KK_{n-1}/\RR_{n-1}$-conjugation (i.e.\ composing with an inner automorphism of
$\G_n$ given by conjugation with an element of the subgroup
$\KK_{n-1}/\RR_{n-1}$). Thus if we let
${\rm Rep}^{\ast}_{\G_{n-1}}\left(\pi_1(M),\G_n\right)$ denote the
representations modulo
$\KK_{n-1}/\RR_{n-1}$-conjugations, it follows that this set is in bijection with
$H^1(M;\KK_{n-1}/\RR_{n-1})$ (by the well-known identification of the latter
with derivations modulo principal derivations \cite[p.~195]{HS}). Moreover this
bijection is natural with respect to continuous maps. This establishes part (a)
of Theorem~\ref{characterizing char} below. Moreover, any choice
$x_{n-1}\in\AA_{n-1}$ will (together with $\phi_{n-1}$) induce $\phi_n$ under the
correspondence (from the proof of Theorem~\ref{higher Bforms})
$$
\AA_{n-1}\equiv H_1(M;\RR_{n-1})\cong H^2(M; \RR_{n-1})
\cong H^1(M;\KK_{n-1}/\RR_{n-1}).
$$ We will refer to this as {\it  the coefficient system corresponding to
$x_{n-1}$ {\rm (}\/and}
$\phi_{n-1}$). This coefficient system is well-defined up to conjugation. It is 
also sometimes convenient to think of (the image of) this element $x_{n-1}$ as
living in
$\AA_{n-1}^{\#}=\Hom_{\RR_{n-1}}(\AA_{n-1},\KK_{n-1}/\RR_{n-1})$
 under the Kronecker map. This image is called the {\it  character induced by
$x_{n-1}$}.
 Indeed it is important to note at this point that
$$
\phi_n:\pi_1(M)\ra\G_n=\KK_{n-1}/\RR_{n-1}\rtimes\G_{n-1}
$$ induces a map from $\pi_1(M_{n-1})$, the $\G_{n-1}$ cover defined by
$\phi_{n-1}$, to $\KK_{n-1}/\RR_{n-1}$, and that the abelianization of this map
$H_1(M_{n-1})\ra\KK_{n-1}/\RR_{n-1}$ is precisely the character induced by
$x_{n-1}$ as above. This is true by construction. Finally, given $\phi_n$ we can
define the $n^{\rm th}$ Alexander module $\AA_n\equiv H_1(M;\RR_n)$. Hence
$$
\AA_n=\AA_n(x_0,x_1,\dots,x_{n-1})
$$ is a function of the choices $x_i\in\AA_i$.

Of course, $\AA_0$ is  $H_1$ of the $\G_0$ cover of
$M$. Given $x_0\in\AA_0$, a ``$\KK_0/\RR_0$-cover'' of the $\G_0$ cover is
induced and
$\AA_1$ is $H_1$ of this composite $\G_1$-cover modulo $S_1$-torsion where
$S_1$ is the set of elements of $\Z\G_1$ which have inverses in
$\RR_1$. Generally $\AA_n$ is $H_1$ of the $\G_n$-cover of $M$, modulo
$S_n$-torsion. In summary we have the following:

\proclaim{Theorem}\label{characterizing char} Suppose $\{\G_n\mid n\geq0\}$ are as in
the beginning of Section~{\rm \ref{sec:Alexander}
 (}\/but not necessarily as in Definition~{\rm \ref{def:univ gr}).} Suppose
$M$ is a compact manifold and
$\phi_{n-1}:\pi_1(M)\ra\G_{n-1}$ is given.
\begin{itemize}
\item[{\rm (a)}] There is a bijection
$f:H^1(M;\KK_{n-1}/\RR_{n-1})\longleftrightarrow{\rm Rep}^{\ast}_{\G_{n-1}}
\left(\pi_1(M),\G_n\right)$ which is natural with respect to continuous maps{\rm ;}
\item[{\rm (b)}]  If
$M$ is a closed oriented $3$\/{\rm -}\/manifold with
$\b_1(M)=1$ then the isomorphism
$H_1(M;\RR_{n-1})\cong H^1(M;\KK_{n-1}/\RR_{n-1})$ with $f$ gives a natural
bijection
$\tilde
f:\AA_{n-1}\longleftrightarrow{\rm Rep}^{\ast}_{\G_{n-1}}\left(\pi_1(M),\G_n\right)$.
\item[{\rm (c)}] If $x\in\AA_{n-1}$ then the character induced by $x$ is given by
$y\mapsto B\l_{n-1}(x,y)$.
\end{itemize}
\endproclaim

{\it Extension of characters and coefficient systems to bounding
$4$\/{\rm -}\/manifolds.}
Suppose $M=\partial W$ and $\phi:\pi_1(M)\ra\G_n$ is given. When does $\phi$
extend over $\pi_1(W)$? In general this is an extremely difficult question
because of our relative ignorance about the types of groups which may occur as
$\pi_1(W)$. This problem has obstructed the generalization of the invariants of
Casson and Gordon and no doubt blocked many other assaults (for example, see
\cite[Cor.~5.3]{KL}, \cite{N}, \cite{L3}, \cite{Let}). Our success in this
regard is the crucial element in defining concordance invariants. If $M$ is the
zero surgery on a slice knot (or more generally an $(n)$-solvable knot) and $W$
is the $4$-manifold which exhibits this (i.e.\ the complement of the slice disk in
the first case) then, under some restrictions on the family $\G_i$, $i\le n$, we
will show that (loosely speaking)
$1/2^n$ of the possible representations from $\pi_1(M)$ to $\G_n$ extend to
$\pi_1(W)$. In particular as long as the generalized Alexander modules $\AA_0$,
$\AA_1,\dots,\AA_{n-1}$ are nonzero there exist nontrivial representations
$\phi$ which do extend. This allows for the construction of an invariant in
$L^0(\KK_n)/i_{\ast}\left(L^0(\Z\G_n)\right)$, which is discussed in \pagebreak
Section~\ref{sec:Witt}.

For the following, let $\G_{n-1}$ be an {\it  arbitrary\/}
$(n-1)$-solvable PTFA group and suppose
$\G_n=\KK_{n-1}/\RR_{n-1}\rtimes\G_{n-1}$ as in Section~\ref{sec:Alexander}. We
need not assume that $\G_{n-1}$ is constructed as in Section~\ref{sec:Alexander}.
We proceed inductively by assuming
$\phi_{n-1}:\pi_1(M)\ra\G_n\ra\G_{n-1}$ already extends to
$\pi_1(W)$.

\proclaim{Theorem}\label{extend}\hskip-8pt  Suppose $M\!=\!\partial W$ with
$\b_1(M)\!=\!1$ and $\phi_n:\pi_1(M)\ra\G_n$\break is given{\rm ,} where
$\G_n=\KK_{n-1}/\RR_{n-1}\rtimes\G_{n-1}$ as in Section~{\rm \ref{sec:Alexander} (}\/but
$\G_{n-1}$ is\break allowed to be an {\rm arbitrary} {\rm PTFA} group\/{\rm ).} Assume that the
nontrivial map
$\phi_{n-1}:\pi_1(M)\ra\G_{n-1}$ extends to a map
$\psi_{n-1}:\pi_1(W)\ra\G_{n-1}$ and that $\phi_n$ is a representative of a class
in
${\rm Rep}^{\ast}_{\G_{n-1}}\left(\pi_1(M),\G_n\right)$  corresponding to $x\in
H_1(M;\RR_{n-1})$. Let
$$ P_{n-1}\equiv\mbox{$\Ker$}\{j_{\ast}:H_1(M;\RR_{n-1})\ra H_1(W;\RR_{n-1})\}.
$$ Then
{\rm 1.}  If $\RR_{n-1}$ is a {\rm PID (}\/or if $j_{\ast}$ is surjective\/{\rm ),} then
$\phi_n$ extends to $\pi_1(W)$ if and only if
$x\in P^\bot_{n-1}$. {\rm (}\/Recall that $P^\bot_{n-1}=\{x\in
H_1(M;\RR_{n-1})|B\l_{n-1}(x,p)=0$\ \ $\forall p\in P_{n-1}\}${\rm .)}
\vglue3pt {\rm 2.} If $M$ is $(n)$\/{\rm -}\/solvable via $W$ then $\phi_n$ extends if and only if
$x\in P_{n-1}$.
\endproclaim

The reader will note that Theorem~\ref{extend}.2 applies to any slice knot. The
difficulty with using Theorem~\ref{extend}.2 is that in applications, often $W$
is unknown and one cannot insure that $P_{n-1}$ is nontrivial. In
Theorem~\ref{self annihil} we shall show that if the hypotheses of {\it  both\/}
\ref{extend}.1 and .2 are satisfied then
$P_{n-1}=P^\bot_{n-1}$. This then is a useful condition which says that $1/2$ (in
a loose sense) of these $\phi_n$ extend. The astute reader will note that
Theorem~\ref{self annihil} is a logical consequence of Theorems~\ref{extend}.1 and
 2. For this reason and because, with our current knowledge,
Theorem~\ref{extend}.2 is useless without
\ref{extend}.1, we shall postpone the proof of Theorem~\ref{extend}.2 until after
Theorem~\ref{self annihil}.

\demo{Proof of Theorem~{\rm \ref{extend}.1}} If $\RR_{n-1}$ is a PID then
$\KK_{n-1}/\RR_{n-1}$ is an injective
$\RR_{n-1}$-module (since it is clearly divisible)~\cite[I~Prop.~6.10]{Ste}. Since
$j_{\ast} : H_1(M;\RR_{n-1})/P_{n-1}\to H_1(W;\RR_{n-1})$ is a monomorphism,
$j^\# : H_1(W)^\#\to\left(H_1(M)/P_{n-1}\right)^\#$ is surjective. Therefore the
``character''
$B\l_{n-1}(x,\cdot)$ extends to $H_1(W;\RR_{n-1})$ if and only if it annihilates
$P_{n-1}$, i.e.\ if $x\in P_{n-1}^\bot$.

Since the bijection between
$H^1\left(\pi_1(M);\KK_{n-1}/\RR_{n-1}\right)$ and
${\rm Rep}^{\ast}_{\G_{n-1}} (\pi_1(M)$, $\G_n)$ is functorial and since
the Kronecker map
$$ H^1(\pi_1(W);\KK_{n-1}/\RR_{n-1})\ra\Hom(\pi_1(W);\KK_{n-1}/\RR_{n-1})
$$ is an isomorphism (as in the proof of Theorem~\ref{higher Bforms}), the
extension of the ``character''
$B\l_{n-1}(x,\cdot)$ is equivalent to an extension of
$\phi_n$ on the $\pi_1$ level.

A similar argument works if $j_{\ast}$ is surjective since this implies that
$j^\#$ is an \pagebreak isomorphism.
\enddemo

\section{Linking forms and Witt invariants\\ as obstructions to
solvability}\label{sec:Witt}
\vglue8pt 

In this section we introduce knot invariants that we prove are defined for
$(n)$-solvable knots and 
 vanish for $(n.5)$-solvable knots. This allows us to state our main theorem
concerning the
 existence of higher-order obstructions to a knot's being slice. These invariants
lie in Witt groups of hermitian forms and 
 are closely related to the Witt classes of our higher-order linking forms via a
localization
 sequence in $L$-theory. In this section we also ask what can be said about a
higher-order linking form $B\l$ on $M$ as in Theorem~\ref{higher Bforms}
 if $M$ is the boundary of a certain type of
$4$-manifold over which the coefficient system extends. A consequence of our
answer to this question is that for $(n)$-solvable knots certain large families
of the higher-order linking forms $B\l_0$,...,$B\l_{n-1}$ are hyperbolic.

Suppose $M$ is equipped with a nontrivial PTFA coefficient system
$\phi:\pi_1(M)\ra\G$ that extends to
$\pi_1(W)$ where $M$ is the boundary of $W$ and
$j_{\ast} : H_1(M;\Q)\to H_1(W;\Q)\cong\Q$ is an isomorphism. Then, since
$H_*(M;\KK)=0$ by Proposition~\ref{torsion} and Remark~\ref{postlore}.3, the
chain complex of the induced cover of $W$ with coefficients in $\KK$ is a {\it  $4$\/{\rm -}\/dimensional symmetric
 Poincar{\rm \'{\it e}} complex} over $\KK$, called the {\it  symmetric
chain complex} of
$W$, and hence represents an element B in $L^0(\KK)$, the cobordism classes of
such complexes
 \cite[pp.1--24]{Ra2}. Since all $\KK$-modules are free, this complex is known to
be cobordant to
 one given by the intersection form $\lambda$ on $H_2(W;\KK)$ (which is
nonsingular by the above remarks and is discussed in detail in
Section~\ref{sec:surfaces}) \cite[Lemma 4.4 (ii)]{Da}. Moreover, in this case
$L^0(\KK)$ is known to be isomorphic to the usual Witt group of nonsingular 
hermitian forms on finitely-generated $\KK$ modules. Let
$W^{\prime}$
 be another such 4-manifold and $B^{\prime}$ the corresponding class. Let V be
the closed 4-manifold obtained by taking the union of $W$ and $-W^{\prime}$
 along M and consider the the symmetric complex of V with $\Z\G$ coefficients
(the {\it  symmetric signature} of V). Let A denote the image of this element of
$L^0(\Z\G)$ under the map
$i_{\ast} : L^0(\Z \G) \to L^0( \KK)$. Thus A is the symmetric signature of V
with $\KK$ coefficients which, as above, is equal to that obtained from
 the intersection form on $H_2(V;\KK)$. Since, by a Mayer-Vietoris sequence, the
latter is the difference of $B$ and $B^{\prime}$, we see that $B =  B(M,\phi)$ is well defined (independent of
$W$) modulo the image of $i_{\ast}$.

In Section~\ref{sec:L2} we discuss $L^{(2)}$-signature invariants which can detect the nontriviality
of $B(M, \phi)$.  Specifically, a homomorphism $\sigma_{\Gamma}^{(2)} : L^0(\KK) \to \Bbb  R$ is defined
which is equal to the ordinary signature $\sigma$ on the image of $L^{(0)}(\Bbb  Z \G)$.   This
$L^2$-signature has additivity properties similar to $\sigma$.  Then, given $(M, \phi)$ as above, one can
define the {\it  reduced $L^2$-signature} (von Neumann $\rho$-invariant) $\rho_{\G}(M, \phi) =
\sigma_{\G}^{(2)}(B) -
\sigma(W)$, a real number independent of $W$.

In this section the groups $\G$ are general PTFA groups unless specified
otherwise. All
$3$- and $4$-manifolds are compact, connected and oriented.  Recall from Section~1
that
$W^{(n)}$ denotes the regular cover of $W$ corresponding to the $n^{\rm th}$ derived
subgroup
$\pi_1(W)^{(n)}$.

\numbereddemo{Definition}\label{Qsolv} The  manifold $M$ is {\it  rationally $(n)$-solvable via
$W$} if it is the boundary of a compact $4$-manifold $W$ such that the inclusion
induces an isomorphism on $H_1(\ ;\Q)$ and such that $W$ admits a {\it  rational
$(n)$-Lagrangian\/} with {\it  rational $(n)$-duals}; that is, there exist classes
$\{\l_1,\dots,\l_m\}$ and $\{d_1,\dots,d_m\}$ in $H_2(W^{(n)};\Q)$ such that
$\lambda_n(\l_i,\l_j)=0$ and
$\lambda_n(\l_i,d_j)=\delta_{ij}$, and where the class images (under the covering map)
together form a basis of $H_2(W;\Q)$. $M$ is {\it  rationally 
$(n.5)$-solvable\/} if {\it  in addition} there exist classes
$\{\l^{\prime}_1,\dots,\l^{\prime}_m\}$ of $H_2(W^{(n+1)};\Q)$ which map to
$\l_i$ as above and such that
$\lambda_{n+1}(\l^{\prime}_i,\l^{\prime}_j)=0$ . It follows that $\sigma(W) = 0$.  Note that if M
is
$(h)$-solvable (see  Sections~1 and~8) then M is rationally
$(h)$-solvable.  
\enddemo

\proclaim{Theorem}\label{B=0} Suppose $\G$ is an $(n)$\/{\rm -}\/solvable group. If $M$ is
rationally
$(n.5)$\/{\rm -}\/solvable via a
$4$\/{\rm -}\/manifold $W$ over which the coefficient system $\phi$ extends{\rm ,} then
$B(M,\phi) = 0$ and $\rho_{\Gamma}^{(2)}(M, \phi) = 0$.
\endproclaim

\demo{Proof}  Let $$L=\{\l_1,\dots,\l_m\}\subset
H_2(W^{(n+1)};\Q)=H_2\left(W;\Q\!\left[\pi_1(W)/\pi_1(W)^{(n+1)}\right]\right)$$
be the rational $(n+1)$-Lagrangian. Since $\G$ is $(n)$-solvable,
$\psi:\pi_1(W)\ra\G$ factors through the quotient
$\pi_1(W)/\pi_1(W)^{(n+1)}$. Using this we can let $L'$ be the submodule
generated by the image of $L$ in
$H_2(W;\KK)$. By naturality, the intersection form with $\KK$ coefficients,
$\lambda$, vanishes on $L'$. Since all $\KK$-modules are free,
$L'$ is a free summand of $H_2(W;\KK)$. It suffices therefore to show that $\rk
L'$ is one-half that of
$H_2(W;\KK)$. The latter has rank equal to $2m$ by the first part of
Proposition~\ref{rank coef} below. We are given that the image of $L$ in
$H_2(W;\Q)$ is linearly independent. By the flatness of $\KK$, it is sufficient
(and necessary) to show that
 $\{\l_1,\dots,\l_m\}$ is linearly independent in
$H_2(W;\Q\G)$. Now apply the second part of the proposition below with $n=3$,
noting that by assumption $H_3(W;\Q) \cong H^1(W,M;\Q)=0$.  Thus $B(M, \phi) = [\lambda] = 0$ and
hence $\sigma_{\G}^{(2)}(B(M, \phi)) = 0.$  Since $\sigma(W) = 0$, $\rho_{\G}^{(2)}(M,\phi) = 0$ as well.
\enddemo

\proclaim{Proposition}\label{rank coef} \hskip-8pt Suppose $W$ is a compact{\rm ,} connected{\rm ,} oriented
$4$\/{\rm -}\/manifold with connected boundary $M$ such that
$ H_1(M;\Q)\ra H_1(W;\Q)$ is an isomorphism. Suppose $\phi:\pi_1(W)\ra\G$ is a
nontrivial {\rm PTFA} coefficient system. Then
$$
\rk H_2(W;\KK)\le\b_2(W)
$$ with equality if $\b_1(W)=1$.

Now suppose $W$ is a connected {\rm (}\/possibly infinite\/{\rm )}
$n$\/{\rm -}\/complex such that $H_n(W;\Q)=0$ and there exist
$(n-1)$\/{\rm -}\/dimensional manifolds $S_i${\rm ,} continuous maps
$f_i:S_i\to W$ and lifts $\tilde f_i:S_i\to W_\G$ such that
$\{[f_i]\mid i\in I\}$ is linearly independent in
$H_{n-1}(W;\Q)$. Then $\{[\tilde f_i]\mid i\in I\}$ is $\Q\G$ linearly
independent in
$H_{n-1}(W;\Q\G)$.
\endproclaim 

\demo{Proof}  Note that any compact topological 4-manifold has the homotopy type
of a finite simplicial complex \cite[Th.~4.1]{KS}.
Choose a finite, $3$-dimensional CW structure for $W$ and let $C_{\ast}(W)$
denote the cellular
chain complex of W with $\Q$ coefficients. Let $b_i=\rk H_i(W;\KK)$.
Then $b_0=0$ by Proposition~\ref{acyclic} and by Proposition~\ref{torsion},
$$ 
b_1\le\b_1(W)-1.
$$
Since
$H_3(W;\Q)\cong H^1(W,\partial W;\Q)=0$, the boundary homomorphism
$\partial:C_3(W)\to C_2(W)$ is injective. Let
$C_{\ast}(W_\G)$ be the corresponding $\Q\G$ chain complex free on the cells of
$W$. By Proposition~\ref{injective},
$\partial :C_3(W_\G)\to C_2(W_\G)$ is injective so that 
$H_3(W;\Q\G)\cong H_3(C_{\ast}(W_\G))=0$. Hence
$b_3=0$ by Remark~\ref{facts}.1. Finally, as noted in the proof of
Lemma~\ref{2comp},
$\chi(W;\Q)=\chi(W;\KK)$ so that we get $b_2\le\b_2(W)$, with equality if
$\b_1(W)=1$ (by Proposition~\ref{torsion}). This completes the proof of the first
part of the proposition.

Let $X$ be the one point union of the $S_i$ (using some base paths), and define
$f:X\to W$ and $\tilde f:X\to W_\G$ to restrict to the given maps on the $S_i$.
After taking mapping cylinders, we may assume $C_{\ast}(X)$ is an
$(n-1)$-dimensional subcomplex of the $n$-dimensional
$C_{\ast}(W)$, and similarly for $C_{\ast}(W_\G)$ and the subcomplex
$C_{\ast}(X_\G)$ where $X_\G$ is the induced cover of $X$. Then $C_i(W)$ is
naturally identified with
$C_i(W_\G)\otimes_{\Q\G}\Q$ and 
$$ p_\# : C_{\ast}(W_\G) \to C_{\ast}(W)
$$ coincides with the obvious homomorphism defined using the augmentation. The
hypothesis is that $f_{\ast}$ is injective on $H_{n-1}(\  \  ;\Q)$.

Since $\phi\circ f_{\ast}$ is trivial on $\pi_1$, $X_\G$ is the trivial cover
consisting of $\G$ copies of~$X$. Thus
$H_{n-1}(X_\G;\Q)$ is a free $\Q\G$-module on $\{\tilde f_i\}$, and hence to
establish the result, we must show that $\tilde f_{\ast}$ is injective on
$H_{n-1}$. Note that, as in the proof of the first part of the proposition,
$H_n(W_\G;\Q)=0$ (finiteness is not needed), and it follows that $\tilde
f_{\ast}$ is injective on
$H_{n-1}$ if and only if $H_n(W_\G,X_\G;\Q)$ is zero. The latter is equivalent to
the injectivity of
$$
\partial^{{\rm rel}}:C_n(W_\G,X_\G)\to C_{n-1}(W_\G,X_\G).
$$ By Proposition~\ref{injective} it suffices to see that
$$
\partial^{{\rm rel}}\otimes\id:C_n(W_\G,X_\G)\otimes_{\Q\G}\Q\to
C_{n-1}(W_\G,X_\G)\otimes_{\Q\G}\Q \pagebreak
$$ is injective. The last statement is equivalent to the vanishing of
$H_n(C_{\ast}(W_\G,X_\G)\break\otimes_{\Q\G}\Q)$. But
$C_{\ast}(W_\G,X_\G)\otimes_{\Q\G}\Q$ can be identified with
$C_{\ast}(W,X;\Q)$. Since $f_{\ast}$ is injective on
$H_{n-1}$ by hypothesis and since $H_n(W;\Q)=0$, it follows that $H_n(W,X;\Q)$
vanishes.
\enddemo

Now we can show that if $K$ is a slice knot (even in a rational homology ball)
and a chosen coefficient system extends to the $4$-manifold {\it  and}
coefficients are taken in a PID $\RR$, then the induced higher-order linking form
on the higher-order Alexander module is hyperbolic. (In fact, under the more
general conditions of Theorem~\ref{higher Bforms}, we can show  that these forms
are {\it  stably} hyperbolic, but this more general result is not needed here).
The consequences of this for the higher-order Alexander polynomials will be
discussed in a later paper (compare
\cite[Cor.~5.3]{KL}).

\proclaim{Theorem}\label{self annihil} Suppose $M$ is rationally
$(n)$\/{\rm -}\/solvable via $W${\rm ,} $\b_1(M)=1$ and
$\phi:\pi_1(M)\ra\G$ is a nontrivial coefficient system that extends to
$\pi_1(W)$ and
$\G$ is an $(n-1)$\/{\rm -}\/solvable {\rm PTFA} group. If $\RR$ is a {\rm PID} such that
$\Q\G\subset \RR\subset\KK$ then the linking form
$B\l(M,\phi)$ {\rm (}\/as defined in Theorem~{\rm \ref{higher Bforms})} is hyperbolic{\rm ,} and in
fact the kernel of
$j_{\ast}:H_1(M;\RR)\ra H_1(W;\RR)$ is self\/{\rm -}\/annihilating.
\endproclaim

\demo{Proof}  Let $P=\mbox{$\Ker$}\{j_{\ast}:H_1(M;\RR)\ra H_1(W;\RR)\}$. Since
all finitely generated modules over a principal ideal domain are homological
dimension at most~1, it suffices to show $P=P^\bot$ with respect to
$B\l$~\cite[p.~253]{Ra2}. We now need the following crucial lemma.

\proclaim{Lemma}\label{Nopid} Assume the hypotheses of Theorem~{\rm \ref{self annihil}}
except that here we do not need that $\RR$ is a {\rm PID}. Then
$TH_2(W,M;\RR)\xrightarrow{\partial}H_1(M, \RR)\xrightarrow{j_{\ast}} H_1(W;\RR)$
is exact. {\rm (}\/Recall that $TH_2$ denotes the $\RR$\/{\rm -}\/torsion submodule{\rm .)}  Moreover\break $H_2(W,
\RR)$
 is the direct sum
of its torsion submodule and a free module.
\endproclaim

\demo{Proof}  Let $\{\l_i,d_i\mid i=1,\dots,m\}$ denote the classes in
$$ H_2(W;\Q[\pi_1(W)/\pi_1(W)^{(n)}])
$$ generating the rational $(n)$-Lagrangian and its duals. Since $\G$ is
$(n-1)$-solvable, the coefficient system $\psi:\pi_1(W)\ra\G$ (extending $\phi$)
descends to $\overline\psi:\pi_1(W)/\pi_1(W)^{(n)}\ra\G$. Let
$\{\l'_i,d'_i\}$ denote the images of $\{\l_i,d_i\}$ in
$H_2(W;\RR)$. By naturality, these are still dual and the intersection form
$\lambda$, with coefficients in $\RR$, vanishes on the span of
$\{\l'_i\}$. Consider $\RR^m\oplus \RR^m$, the free module on
$\{\l'_i,d'_i\}$, and the composition 
$$
\RR^m\oplus \RR^m\xrightarrow{i_{\ast}} H_2(W;\RR)\xrightarrow{\lambda}
H_2(W;\RR)^*
\xrightarrow{i^{\ast}}(\RR^m\oplus
\RR^m)^*.
$$ This map is represented by a block diagonal matrix of the form 
$$
\left(\begin{array}{cc} 0  &I\\  I  &X\end{array}\right)
$$ for some $m\times m$ matrix $X$. This matrix has
$$
\left(\begin{array}{cc} -X  &I\\  I  &0\end{array}\right)
$$ as its inverse implying that $i^{\ast}$ is a (split) epimorphism and
$i_{\ast}$ is a monomorphism. Since the $\KK$-rank of 
$$ H_2(W;\RR)\otimes_\RR \KK  \cong H_2(W;\KK)
$$
 is $2m$ by Proposition~\ref{rank coef}, the cokernel $C$ of $i_{\ast}$ is a
torsion module, and thus
$\Hom_\RR(C,\RR)=0$. Hence, applying the functor $\Hom_\RR(-,\RR)$ to the map
$i_{\ast}$, we see that its Hom-dual $i^{\ast}$ is injective. Therefore
$i^{\ast}$ is an isomorphism. It follows that $\lambda$ is surjective (and hence $H_2(W, \RR)$ is a direct sum
of a free module of rank $2m$ and its torsion module.)  Now consider the commutative diagram below for
(co)-homology with
$\RR$-coefficients.
$$
\BoxedEPSF{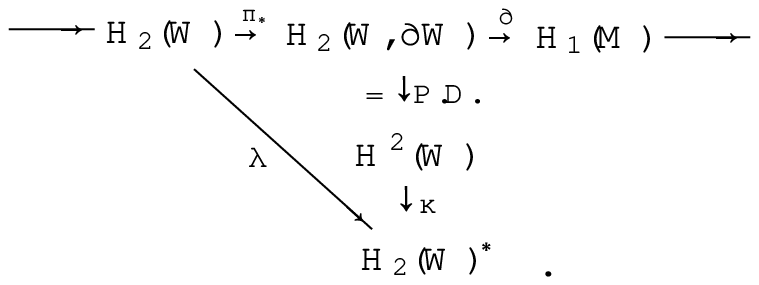 scaled 1000}
$$ 
Note that  $\kappa$ is a split surjection between modules of the same rank over
$\KK$ and thus the kernel of $\kappa\circ\operatorname{P.D.}$ is torsion. Now,
given $p\in P$, choose $x$ such that
$\partial x=p$. Let $y$ be an element of the set $\lambda^{-1}(\kappa\circ
\operatorname{P.D.}(x))$. Then $\partial(x-\pi_{\ast}(y))=p$ and
$x-\pi_{\ast}(y)$ is torsion since it lies in the  kernel of
$\kappa\circ\operatorname{P.D.}$. This concludes the proof of the lemma.
\enddemo

Continuing the proof of Theorem~\ref{self annihil}, consider the following
diagram, commutative up to sign, where coefficients are in $\RR$ unless otherwise
specified:
\begin{equation}\label{diag:extend}
\begin{array}{cccccc} 
{TH_2(W,\partial W)} &\stackrel{\partial}{\lrar}&{H_1(M)} &\stackrel{j_{\ast}}{\lrar}&{H_1(W)}\\
 {\cong}\Big\downarrow {\rm P.D.}&&\cong\Big\downarrow {\rm P.D.}\\
{TH^2(W)} &\stackrel{j^{\ast}}{\lrar} &{H^2(M)}\\
{ \cong}\Big\downarrow{B^{-1}}&&
 { \cong}\Big\downarrow{B^{-1}}\\
{H^1(W;\KK/\RR)}&\stackrel{j^{\ast}}{\lrar}& {H^1(M;\KK/\RR)}\\
\Big\downarrow{\kappa}&&
\Big\downarrow{\kappa}\\
{H_1(W)^{\#}} &\stackrel{j^{\#}}{\lrar}&{H_1(M)^{\#}}&&&.
\end{array}   \pagebreak
\end{equation}  

The vertical homomorphisms above are Poincar\'e Duality, inverse of the Bockstein
$B$ and the Kronecker evaluation map $\kappa$. These compositions are denoted
$\b_{{\rm rel}}$ and
$B\l$ respectively. To see that this ``linking form''
$\b_{{\rm rel}}$ exists, examine the sequence
$$ H^1(W;\KK)\ra H^1(W;\KK/\RR)\xrightarrow{B} H^2(W;\RR)\ra H^2(W;\KK)
$$ and note that $H^1(W;\KK)=0$, $H^2(W;\KK)$ is $\RR$-torsion free and all
homology with coefficients in $\KK/\RR$ is $\RR$-torsion. It follows that $B$ is
an isomorphism onto $TH^2(W;\RR)$.  If $x\in P$ then $x=\partial y$ using
Lemma~\ref{Nopid}. Thus
$B\l(x)=j^\#(\b_{{\rm rel}}(y))$ and hence, for any $p\in P$,
$B\l(x)(p)=\b_{{\rm rel}}(y)(j_{\ast}(p))=0$ so that $x\in P^\bot$. Hence $P\subset
P^\bot$. $\phantom{\sum^1}$

Finally, we will use the fact that $\RR$ is a PID to show that
$P^\bot\subset P$. Consider the monomorphism
$H_1(M;\RR)/P\xrightarrow{j_{\ast}} H_1(W;\RR)$. Clearly, $\KK/\RR$ is  a divisible
$\RR$-module which implies it is injective since $\RR$ is a PID \cite[I
6.10]{Ste}. Therefore the map
$$ j^\#:H_1(W;\RR)^\#\ra (H_1(M;\RR)/P)^\#
$$ is onto. Now, given $x\in P^\bot$, it follows that
$B\l(x)(p)=0$ for all $p\in P$ so $B\l(x)$ lifts to an element of
$(H_1(M)/P)^\#$. Thus
$B\l(x)$ lies in the image of
$j^\#$.  Moreover, the Kronecker map
$$
\kappa:H^1(W;\KK/\RR)\ra H_1(W;\RR)^\#
$$  is an isomorphism since $\RR$ is a PID (see proof of Theorem~\ref{higher
Bforms}). By Diagram~\ref{diag:extend}, $x$ lies in the image of
$\partial$ and so $x\in P$. Hence $P^\bot\subset P$.\hfill\qed
\enddemo

\demo{Proof  of Theorem~{\rm \ref{extend}.2}} Note that Lemma~\ref{Nopid} holds.
Consider Diagram~\ref{diag:extend}. The Kronecker maps may no longer be
isomorphisms so ignore them. If $x\in P_{n-1}$ then $x=\partial y$ as above so
the image of $x$ in
$H^1(M;\KK/\RR)$ is in the image of $j^{\ast}$. Recall that extensions as in
Theorem~\ref{characterizing char} correspond fundamentally to these cohomology
classes and the proof is finished as in the second paragraph of the proof of
Theorem~\ref{extend}.1.
\enddemo

We can now prove our main theorem by combining Theorems~\ref{B=0}, \ref{self
annihil},
\ref{characterizing char}, and \ref{extend}. This can be applied to the zero
surgery on a knot $K$ in a rational homology sphere, or to a prime-power 
cyclic cover of such a manifold.

If $\b_1(M)=1$ then, modulo torsion, $H_1(M)\cong H_1(W)\cong\Z$, and the
inclusion induces multiplication by some nonzero integer, whose absolute value
we call the {\it  multiplicity\/}. Note that if $M=\partial W$ with
$j_*:H_1(M;\Q)\ra H_1(W;\Q)\cong\Q$ an isomorphism, then there are precisely two
{\it  epimorphisms}
$\psi_0:\pi_1(W)\ra\gop=\Z$. Let
$\phi_0=\psi_0\circ j_*$. This map $\phi_0:\pi_1(M)\ra\gop$ (up to sign) is
canonically associated \pagebreak to $M$ and the multiplicity of $M\ra W$,  and extends to
$\pi_1(W)$ by definition. If $j_*$ is an isomorphism on integral homology, as is
the case for a slice knot in an integral homology $4$-ball, the multiplicity is
$1$ and $\phi_0$ is the canonical epimorphism.

\proclaim{Theorem}\label{invariant thm} Let $\gop${\rm ,}
$\G^U_1,\dots,\G^U_n$ be the family of universal groups of
Definition~{\rm \ref{def:univ gr}}. Suppose $M$ is a closed{\rm ,} oriented{\rm ,} $3$\/{\rm -}\/manifold
with $\b_1(M)=1$. Then
\begin{itemize}
\item[$(0)${\rm :}] If $M$ is rationally $(0)$\/{\rm -}\/solvable via $W_0$ then either of the
two maps
$\phi_0:\pi_1(M)\ra\gop$ extends to $\pi_1(W_0)$ inducing a class
$B_0=B(M,\phi_0)$ in
$L^0(\KK_0)$ {\rm (}\/modulo the image of $L_0(\Bbb  Z \G_0^U))$. Moreover
$\phi_0$ induces an Alexander module $\AA_0$ and a nonsingular Blanchfield
linking form $B\l_0$.
\item[$(0.5)${\rm :}] If $M$ is rationally $(0.5)$\/{\rm -}\/solvable via
$W_{0.5}$ then{\rm ,} in addition to the above holding for
$W_{0.5}${\rm ,} $B_0=0 = \rho_{\G_0}^{(2)}(M,\phi_0)$.
\item[$(1)${\rm :}] If $M$ is rationally $(1)$\/{\rm -}\/solvable via $W_1$ then{\rm ,} in addition to
the above holding for $W_1${\rm ,}
$P_0\equiv\Ker\{j_{\ast}:\AA_0\ra\AA_0(W_1)\}$ is self\/{\rm -}\/annihilating for $B\l_0$
and for any $p_0\in P_0$ a coefficient system $\phi_1(p_0):\pi_1(M)\ra\G^U_1$ is
induced which extends to $\pi_1(W_1)$ and induces a class
$B_1(p_0)=B(M,\phi_1)$ in
$L^0(\KK_1)$. Here $\AA_0(W_1)=H_1(W_1;\RR^U_0)$. Moreover
$\phi_1$ induces the generalized Alexander module $\AA_1(p_0)$ and nonsingular
linking form $B\l_1(p_0)$.
\item[$(1.5)${\rm :}] If $M$ is rationally $(1.5)$\/{\rm -}\/solvable via
$W_{1.5}$ then{\rm ,} in addition to the above holding for
$W_{1.5}${\rm ,} $B_1(p_0)=0 = \rho_{\G_1}^{(2)}(M,\phi_1)$.\\
\vdots
\item[$(n)${\rm :}] If $M$ is rationally $(n)$\/{\rm -}\/solvable via $W_n$ then{\rm ,} in addition to
the above holding for $W_n${\rm ,}
$P_{n-1}\equiv\Ker\{j_{\ast}:\AA_{n-1}\ra\AA_{n-1}(W_n)\}$ is self-annihilating
with respect to
$$
B\l_{n-1}(p_0,\dots,p_{n-2})
$$
 and for any $p_{n-1}\in
P_{n-1}(p_0,\dots,p_{n-2})$ a coefficient system
$$
\phi_n(p_0,\dots,p_{n-1}):\pi_1(M)\ra\gu
$$ is induced which extends to $\pi_1(W_n)$ and induces a class
$$ B_n(p_0,\dots,p_{n-1})=B(M,\phi_n)\in L^0(\KK_n)$$
modulo the image of $L^0(\Z\gu).$ 
\item[$(n.5)${\rm :}]  If $M$ is rationally $(n.5)$\/{\rm -}\/solvable via
$W_{n.5}$ then{\rm ,} in addition to the above holding for
$W_{n.5}${\rm ,} $$B_n(p_0,\dots,p_{n-1}) = 0 = \rho_{\G_n}^{(2)}(M, \phi_n).$$
\end{itemize} In particular{\rm ,} if a knot $K$ is $(n.5)$\/{\rm -}\/solvable then{\rm ,} for any
choices
$(p_0,p_1,\dots,\break p_{n-1})${\rm ,} there exist self\/{\rm -}\/annihilating submodules
$P_i\subset\AA_i(p_0,\dots,p_{i-1})${\rm ,} $0\le i\break <n${\rm ,} and an induced coefficient
system
$\phi_n(p_0,\dots,p_{n-1}):\pi_1(M)\ra\gu$
 {\rm (}\/up to conjugation\/{\rm )} such that $B_n(M,\phi_n)$ and $\rho_{\G_n}^{(2)}(M, \phi_n)$ are defined and equal to
zero. Here
$M$ is zero surgery on $K$.
\endproclaim

\numbereddemo{{R}emark}\label{post inv rem} 1.  If one chooses $p_{n-1}=0$ in Theorem~\ref{invariant thm}~$(n)$ then
$\phi_n$ is ``trivial'' in the sense that it factors through $\phi_{n-1}$ via the
splitting
$\G_{n-1}\ra\G_n$. It follows that $B_n$, $\AA_n$ and
$B\l_n$ are all just $B_{n-1}$, $\AA_{n-1}$ and
$B\l_{n-1}$ ``tensored up to $\Z\G_n$'' in the appropriate fashion. Therefore
there is no additional information and further conclusions are trivial
consequences of previous stages. Consequently if
$\AA_{n-1}(p_0,\dots,p_{n-2})=0$  then no new information can be gleaned just as,
if the classical Alexander module $\AA_0$ is  trivial then Casson-Gordon's
invariants give no information. Indeed M.~Freedman has shown that a knot with
$\AA_0=0$ is topologically slice \cite{F}.  On the other hand, if
$\AA_{n-1}(p_0,\dots,p_{n-2})$ is {\it  nontrivial\/} then the nonsingularity of
$B\l_{n-1}$ guarantees that
$P_{n-1}$ is nontrivial. In fact one can show that if 
${\rm dim}_{\Q} \AA_0>2$, then $\AA_{n-1}$ is {\it  always nontrivial}. This will be
shown in a subsequent paper [C].

\vglue4pt 2. Actually a slightly stronger theorem is true.  One need not use the full
$\G^U_i$ but, once M is fixed, can replace these by a family of universal groups
(defined inductively as semi-direct products) where
 $\KK_{i-1}/\RR_{i-1}$ is replaced by the image of the smallest direct summand of
$H_1(W;\RR_{i-1})$ which contains the image of $\AA_{i-1}$. This leads to a
family of much smaller groups $\G_i$, depending only on the $\AA_i(M)$,
 which are still of the type from Section~\ref{sec:Alexander}. Although we shall
not here formalize this subtlety further, we will use it to advantage in
Section~\ref{sec:example}.
\enddemo

\demo{Proof}  Note that all maps on the  fundamental group will be nontrivial since
$\phi_0$ is. By induction, assume Theorem~\ref{invariant thm}~$(n-1).5$ holds
true. We shall establish \ref{invariant thm}~$(n)$. Suppose $M$ is rationally
$(n)$-solvable via
$W_n$. Then
$M$ is rationally $(h)$-solvable via $W_n$ for any $h<n$. By the induction
hypothesis,
$\phi_0$ extends to $\pi_1(W_n)$, and for any
$p_0\in P_0\equiv\Ker\{j_{\ast}:\AA_0\ra\AA_0(W_n)\}$, 
$\phi_1(p_0)$ is induced which extends to
$\pi_1(W_n)$ (and for any such extension), and $\dots$ for any
$$ p_{n-2}\in
P_{n-2}(p_0,\dots,p_{n-3})\equiv\Ker\{j_{\ast}:\AA_{n-2}\ra\AA_{n-2}(W_n))\}
$$
$\phi_{n-1}(p_0,\dots,p_{n-2}):\pi_1(M)\ra\G^U_{n-1}$ is induced which extends to
$\pi_1(W_n)$ and (for any such extension) induces $\AA_{n-1}$, $B\l_{n-1}$. By
Theorem~\ref{self annihil} we see that $P_{n-1}$ is self-annihilating for
$B\l_{n-1}$. This is the first condition of Theorem~\ref{invariant thm}~$(n)$.
Now choose
$p_{n-1}\in P_{n-1}$. By Theorem~\ref{characterizing char} a coefficient system
$\phi_n:\pi_1(M)\ra\gu$ is induced \pagebreak which extends to $\pi_1(W_n)$ by
Theorem~\ref{extend}.1 or 2. Then
$\phi_n$ induces $B(M,\phi_n)$ as in Section~\ref{sec:Witt} and with  $\AA_n$ and
$B\l_n$ as in Theorem~\ref{higher Bforms}. This establishes \ref{invariant
thm}~$(n)$. To establish \ref{invariant thm}~$(n.5)$, merely apply the above to
$W_{n.5}$ and then apply Theorem~\ref{B=0}.
\enddemo

\proclaim{Theorem}\label{Gammas} Let $\G_0${\rm ,} $\G_1,\dots,\G_n$ be a family of rationally
universal groups as in Definition~{\rm \ref{def:univ gr}) (}$\RR_i${\rm ,} $S_i$ variable\/{\rm )}.
Then Theorem~{\rm \ref{invariant thm}} holds with the following changes. Omit all
conclusions about
$P_i$ being self\/{\rm -}\/annihilating. Replace the conclusion that $B\l_i$ is
nonsingular with the conclusion that $B\l_i$ is nonsingular if $\RR_i$
satisfies the hypothesis of Theorem~{\rm \ref{higher Bforms}.}
\endproclaim

\demo{Proof}  Follow the proof of Theorem~\ref{invariant thm}. Apply
Theorem~\ref{extend}.2 instead of Theorem~\ref{extend}.1.
\enddemo

{\it Bordism invariants generalizing the Arf invariant.} 
The following result could lead to examples of $(n-0.5)$-solvable knots that are
not
$(n)$-solvable, but calculations have not been made.

\proclaim{{C}orollary}\label{spin} Under the hypotheses of Theorem~{\rm \ref{invariant thm},}
suppose
$K$ is rationally $(n)$\/{\rm -}\/solvable {\rm (}\/respectively $(n)$\/{\rm -}\/solvable\/{\rm ).} Then there exists
a submodule $P_0\sbq\AA_0$ which is self\/{\rm -}\/annihilating for $B\l_0$ and for any
$p_0\in P_0$ a coefficient system $\phi_1:\pi_1(M)\ra\G^U_1$ is induced
$\dots$ such that there exists a submodule
$P_{n-1}(p_0,\dots,p_{n-2})\sbq\AA_{n-1}(p_0,\dots,p_{n-2})$ which is
self\/{\rm -}\/annihilating for
$B\l_{n-1}$ and for any $p_{n-1}\in P_{n-1}$ a coefficient system
$\phi_n:\pi_1(M)\ra\gu$ is induced such that the element $(M,\phi_n)$ of
$\Omega_3(B\gu)$ {\rm (}\/respectively $\Omega^{\rm Spin}_3(B\gu)${\rm )} is zero.
\endproclaim

\demo{Proof}   This is a direct corollary of Theorem~\ref{invariant thm}~$(n)$.
\enddemo

Note that the obstruction in the case $n=0$ (the spin case) of
Corollary~\ref{spin} is well-known to be the Arf invariant of $K$. Note also that
there is a somewhat stronger version along the lines of Remark~\ref{post inv
rem}.2.

\section{$L^2$-signatures} \label{sec:L2}  

Given a PTFA group $\G$ and the
quotient field $\KK$ of $\Z\G$, the purpose of this section is to construct a
homomorphism
$$ L^0(\KK) \ra \R
$$ which detects the slice obstructions from Theorem~\ref{invariant thm}. It
turns out that such a homomorphism can be found by completing $\Z\G$, or better
$\C\G$, to the\break von Neumann algebra $\NN\G$ and also completing $\KK$ to the
algebra $\UU\G$ of unbounded operators affiliated to $\NN\G$. Then one can use
the dimension theory of von Neumann algebras to define an $L^2$-signature
$$
\sigma_\G: L^0(\UU\G) \ra \R
$$ for any group $\G$. It agrees with the ordinary signature on the image of
$L^0(\Z\G)$ if the analytic assembly map for $\G$ is onto. This property was
recently established by N. Higson and A. Kasparov \cite{HK} for all torsion-free
amenable groups, a class of groups which contains our PTFA groups.

The idea that the $L^2$-signature can be applied to concordance questions
originated after discussions with Holger Reich on the extension of the von
Neumann dimension to
$\UU\G$. His Ph.D. thesis \cite{Re} was very helpful for writing this section.
Also discussions with Wolfgang L\"uck and Thomas Schick were very useful for
understanding von Neumann algebras in necessary detail. L\"uck's paper
\cite{Lu2} is an excellent survey on the use of von Neumann algebras in topology
and geometry.

We claim no originality in the following section, except for the observation that
this beautiful theory does relate to classical knot concordance. The section is
written for nonexperts in von Neumann algebras.

Let $\G$ be a countable discrete group and  consider the Hilbert space $\l^2 \G$
of square-summable sequences of group elements with complex coefficients.  The
group $ \G$ acts by left- and right- multiplication on
$\l^2 \G$.  These operators are obviously isometries and we can consider the
embedding
$$
\C \G \hookrightarrow {\cal B}(\l^2 \G)
$$ corresponding to (sums of) {\it  right} multiplications into the space of
bounded operators on $\l^2 \G$.

\numbereddemo{Definition}\label{C algebra} The (reduced) $C^{\ast}$-algebra $C^{\ast} \G$ is
the completion of $\C \G$ with respect to the operator norm on
${\cal B}(\l^2 \G)$.  The von Neumann algebra $ \NN \G$ is the completion of
$\C \G$ with respect to pointwise convergence in
${\cal B}(\l^2 \G)$. In particular we have $\C \G \subseteq C^{\ast} \G
\subseteq
 \NN \G$.
\enddemo

If follows from von Neumann's double commutant theorem that $ \NN \G$ is equal to
the set of bounded operators which commute with the {\it  left}
$ \G$-action on $\l^2 \G$:
$$
 \NN  \G={\cal B}(\l^2 \G)^ \G.
$$ From this description, the {\it  standard $\G$-trace} 
$$
\tr_ \G:
 \NN \G \longrightarrow \C
$$
 is defined by $ \tr_ \G(a):=\langle (e)a, e\rangle_{\l^2 \G}$ where $e\in  \G
\subseteq
\l^2 \G$ is the unit element. It follows from the left invariance of $a\in
\mathcal N
 \G$ that for all (isometries)
$g\in  \G$
$$
\langle (e)(ag),e\rangle=\langle (e)a,g^{-1}\rangle=\langle
g((e)a),e\rangle=\langle (ge)a,e\rangle=\langle (e)(ga),e\rangle.
$$
 By the linearity and continuity of the $\G$-trace this implies the usual trace
property
$$
\tr_ \G(ab)=\tr_ \G(ba) \quad\hbox{  for all } a, b \in
\mathcal   N  \G.
$$ Moreover, it also implies that the $\G$-trace is {\it  faithful}, i.e.\ if
$\tr_
\G(a^{\ast}a)=0$ then $(e)a=0$ which implies
$a=0$ by left invariance and continuity. This $\G$-trace extends to $(n\times
n)$-matrices by sending a matrix to the sum of the $\G$-traces of the diagonal
entries.  For example, if
$$ p\in M_n ( \NN \G) =
{\cal B}((\l^2 \G)^n)^ \G
$$ is an orthogonal $ \G$-equivariant projection onto a subspace $V\subseteq (\l^2
\G)^n$ then we may define the {\it  $ \G$-dimension} of $V$ by
$$
\dim_ \G V:=
\tr_ \G(p) \in [0,\infty).
$$ The trace property shows that this actually only depends on the $
\G$-isometry class of the Hilbert space $V$. Thinking of K-theory as equivalence
classes of projections, by implication, we have a map
$$
\tr_\G:K_0({\mathcal N}  \G) \ra \R.
$$ If $\NN\G$ is a factor, i.e.\ has center $\C$, then this map is actually an
isomorphism. This is the case if and only if for all nontrivial $ \gamma\in\G$
the number of elements conjugate to $\gamma$ is infinite.
\vglue2pt
To define the $L^2$-signature, consider a hermitian $(n\times n)$-matrix over
$\NN\G$, $h\in
\Herm_n( \NN \G)$, as a bounded, hermitian
$ \G$-equivariant operator on the Hilbert space
$(\l^2 \G)^n$.  Its spectrum $\spec(h)$ is a (compact) subset of the real line
and for any bounded measurable function $f$ on $\spec(h)$ we may define the
bounded
$ \G$-equivariant operator $f(h) \in M_n( \NN \G)$ by functional calculus (see
e.g.\
\cite{Pe}).  In particular, consider the characteristic functions $p_+,
p_-:\R\to\R$ of
$(0, +\infty)$, respectively $(-\infty, 0)$. 

\numbereddemo{Definition}\label{L2sig} The {\it signature map} $\sigma_\G: \Herm_n( \NN \G) \to
K_0(\NN\G)$ is defined by sending $h$ to the formal difference $p_+(h)-p_-(h)$ of
projections in
$M_n( \NN \G)$. The 
$L^2$-{\it signature} of $h\in \Herm_n( \NN \G)$ is defined to   be the real number

\hfil  ${\displaystyle 
\sigma^{(2)}_\G(h):= \tr_\G(p_+(h)) -
\tr_\G(p_-(h)).
}$\hfill
\vglue2pt
\enddemo 

As the crucial example, we consider the case $ \G = \langle t\rangle
\cong
\Z$.  Then
$\C \G$ consists of Laurent polynomials $\C[t, t^{-1}]$ which embed naturally
into the space of complex-valued continuous functions on the circle $S^1$. 
Indeed, Fourier   transformation gives an isomorphism of Hilbert-spaces $\l^2
\Z  \cong L^2(S^1)$ and pointwise multiplication by a function induces the
isomorphism
$C^{\ast} \G \cong C(S^1)$. This is a consequence of the Stone-Weierstrass 
theorem on the density of polynomials in the space of all continuous functions in
the   supremum norm. Completion in the topology of pointwise convergence then
leads to the von Neumann algebra $ \NN \G$ which turns out to be the space
$L^\infty(S^1)$  of complex-valued, bounded,  measurable functions on the circle,
defined almost everywhere.  Finally, the standard $\G$-trace is just given by
integration.

Now consider $h\in \Herm_n(C(S^1))$ and think of it as a continuous map from
$S^1$ to
$\Herm_n(\C)$.  The ordinary signature $\sigma_0: \Herm_n(\C)\longrightarrow
\Z$ counts the number of positive Eigenvalues minus the number of negative
Eigenvalues.

\numbereddemo{Definition}\label{def:signature} The {\it  twisted signature} of $h$ is the step
function
$\sigma(h): S^1\longrightarrow \Z$ which assigns to each
$s\in S^1$   the signature $\sigma_0(h(s))$. Moreover, the real number
$\sigma^{(2)}(h)$ is defined to be the integral of this function $\sigma(h)$ over
the circle (normalized to have total measure~1).
\enddemo

Thus $ \sigma^{(2)}(h)$ is the average of all the twisted signatures.  It is
clear that
$\sigma(h)$ makes sense almost everywhere for $h\in
\Herm_n(L^\infty(S^1))$ and therefore $\sigma^{(2)}(h)$ is well defined even in
this case. As an example, consider the following element in
$\Herm_2(\C[t^{\pm 1}])$:
$$ h:= \left( \begin{array}{cc}\label{neq0} t+t^{-1} -2 & 1 \\ 1 & t+t^{-1} -2
\end{array} \right).
$$
 Notice that $\sigma(h)$ is a step function with jumps, at most, at the zeroes of
the ``Alexander polynomial" $\det(h) \in
\C[t, t^{-1}]$. We have $\det(h)= (t+t^{-1} -2)^2 -1= (t+t^{-1} -1)(t+t^{-1} -3)$
which has roots on  $S^1$ exactly for the two primitive $6^{\rm th}$ roots of unity.  So
we only need to calculate $\sigma(h)$ at two points on the circle  which
interlace with these two roots, e.g. at $\pm 1$. Clearly $h(1)$ is hyperbolic and
one easily checks that the ordinary signature of $h(-1)$ is
$-2$. One therefore gets
$$
\sigma^{(2)}(h) = (1/3) \cdot 0 + (2/3) \cdot (-2)= -4/3 \neq 0.
$$ 

\proclaim{Lemma}\label{agree} \hskip-8pt The average $ \sigma^{(2)}(h)$ equals the
$L^2$-signature $
\sigma^{(2)}_\G(h)$ for $\G=~\Z$.
\endproclaim

\demo{Proof}  Notice that $\tr_ \G(p_+(h))$ is the $ \G$-dimension of the
``positive Eigen\-space'' of $h$.  In the functional calculus one approximates
(say) $p_+$ by a sequence of real polynomials $p_i$ into which any operator can
easily be   substituted. Then one takes the pointwise limit to define
$$ p_+(h):= \lim_i p_i(h) \quad \hbox{  for } h\in \mathcal N \Z.
$$  
For example, if $h$ is a finite dimensional matrix, then one checks that $p_+$ is just
the projection onto the $(+1)$-Eigenspace of $h$.
This implies that for a point $s\in S^1$, a fancy way to count the number of positive
Eigenvalues of
$h(s)\in\Herm_n(\C)$ is to take the ordinary trace of $p_+(h(s)):= \lim_i
p_i(h(s))$. But now one clearly sees that the integral of the function
$\sigma(h)$ which associates for each $s\in S^1$ the difference
$p_+(h(s))-p_-(h(s))$ is almost everywhere the same as
$\sigma^{(2)}_ \G(h)$.
\enddemo

\numbereddemo{{R}emark}\label{rem:center valued} The above proof actually shows that without
the integration the twisted signature $\sigma(h) \in L^\infty(S^1;\Z)$ as in
Definition~\ref{def:signature} equals the more general {\it  center\/{\rm -}\/valued}
$L^2$-signature map coming from the center-valued trace on $\mathcal N\Z$. Since
$\Z$ is commutative, this trace is just the identity on
$\NN\Z=L^\infty(S^1,\C)$. Moreover, the elements $p_{\pm}(h) \in M_n(\NN\G)$ in
the definition of the
$L^2$-signature are mapped by the corresponding trace on matrices to
$L^\infty(S^1,\R)$. In fact, the proof above shows that they are equal to step
functions almost everywhere. \hfill \qed 
\enddemo

Since the functional calculus used above extends to self-adjoint {\it  unbounded}
operators \cite{Pe}, we can extend the $L^2$-signature to a super-ring $\UU\G$ of
$\NN\G$.  Here $\UU\G$ is the algebra of operators {\it  affiliated} to $ \NN
\G$.  This is the set of (unbounded) operators $a=(a, \dom(a))$ on $\l^2 \G$
which satisfy the following conditions:

\begin{itemize}
\item[(i)]  $a$ is densely defined, i.e.\ $\dom(a)$ is dense in
$\l^2 \G$.
\item[(ii)] $a$ is closed; i.e., its graph is closed in
$\l^2 \G
\times \l^2 \G$.
\item[(iii)] $a$ is affiliated to $ \NN \G$; i.e., for every bounded operator 
$b$  which commutes with all of
$ \NN \G$ we have $ba \subseteq ab$.  This means
$\dom(ba) \subseteq \dom(ab)$ and on the smaller subset these operators agree.
\end{itemize}
 It takes some work to show that $ \UU \G$ is indeed a ring with
involution, for example to define addition and multiplication one has to close
the operators.  Then the various associativity and distributivity laws become
actual theorems.  This was worked out by Murray and von Neumann \cite{MvN} (see
also
\cite{Re}).  It turns out however, that $ \UU \G$ can also be obtained as the
Ore-localization of $\NN \G$ with respect to all nonzero divisors.  From this
point of view the theorem is that the latter set is an Ore-domain.  As an example,
$ \UU\Z$ is the set of all measurable functions on
$S^1$, i.e.\ not necessarily bounded functions. This means that in this example the
$\G$-trace, or integral, does not extend to a map on $\UU\G$.  However, one can
extend the $L^2$-signature to a map
$$
\sigma_\G^{(2)} : \Herm_n(\UU\G)\ra\R
$$  as follows. Observe that a hermitian matrix $h$ with entries in $\UU\G$ can
be viewed as an (unbounded) self-adjoint operator on $\l^2(\G)^n$. Since the two
projections
$p_+$ and $p_-$ onto the positive respectively negative spectrum are {\it  bounded} functions, it follows that the corresponding projections $p_+(h)$ and
$p_-(h)$, obtained via functional calculus, are bounded and thus lie in
$M_n(\NN\G)$. Therefore their
$\G$-traces can be defined as before.

\proclaim{Lemma}\label{lem:independence} The $L^2$\/{\rm -}\/signature only depends on the
$ \G$\/{\rm -}\/isometry class of $h\in\Herm_n(\UU\G)${\rm ;} i.e.{\rm ,} it is unchanged under
$h \mapsto a^{\ast} h a$ for $a\in \GL_n( \UU \G)$.
\endproclaim

\demo{{P}roof}   We first give an argument for the ring $\NN\G$. Consider the
Hilbert space $H:=(\l^2  \G)^n$ with the bounded $
\G$-equivariant operators $h$ and $a$. Let $p_0$ be the characteristic function of
$\{0\}\subset\R$, i.e.\ $p_0(h)$ is the projection onto the kernel of $h$. Then
$p_0(h),p_+(h),p_-(h)$ are commuting projections which sum up to the identity and
therefore their images give an orthogonal decomposition of Hilbert spaces
$$ H=H_0 \perp H_+ \perp H_-.
$$ For a vector $v$ in one of the three summands above, one has by  definition
that
$$
 \langle h(v),v \rangle =0, >0 \hbox{  respectively }<0.
$$ It follows that depending on whether $v$ is in $a^{-1}H_0$,  
$a^{-1}H_+$, respectively $a^{-1}H_-$, one has
$$
\langle (a^{\ast}ha)(v),v \rangle = \langle h(av),av \rangle =0, >0
,  \hbox{ 
respectively, } <0.
$$
 Therefore, the three orthogonal projections
$$ a^{-1}H_\dagger \ra p_\dagger(a^{\ast}ha)H \quad\hbox{  for }  
\dagger\in\{0,+,-\}
$$ are monomorphisms and thus
$$
\dim_{ \G}H_\dagger=\dim_{ \G}a^{-1}H_\dagger\leq \dim_{ \G} 
p_\dagger(a^{\ast}ha)H\quad\hbox{  for } \dagger\in\{0,+,-\}.
$$ But the three dimensions on both sides must sum up to the total dimension
$n$ of
$H$ and therefore the inequalities are actually  equalities.

To extend this argument to the ring $\UU\G$ of unbounded operators affiliated to
$\NN\G$ one has to use some of their properties (see \cite[\S 11]{Re}). Namely,
define a subspace $L$ of the Hilbert space $H$ to be {\it  essentially dense} if
it contains a sequence of closed affiliated subspaces whose $\G$-dimension tends
to one. Here a closed subspace is called affiliated if the corresponding
projection is affiliated to
$\NN\G$.  Then the above proof applies to $\UU\G$ because for all $a\in\UU\G$,
\begin{itemize}
\item $\dom(a)$ is essentially dense, and
\item $a^{-1}(L)$ is essentially dense if $L$ is essentially dense.\hfill\qed
\end{itemize}
\enddemo

If $P$ is a finitely generated projective $\UU\G$-module,  we may choose a
$\UU\G$-module $Q$ such that $P \oplus Q \cong (\UU\G)^n$. If, moreover, 
$$ h:P\ra P^*:=\Hom_{\UU\G}(P,\UU\G)
$$ is a hermitian form in the sense that $h=h^*:P \cong P^{**} \to P^*$, then we
can extend it by three blocks of zeroes to an element in $\Herm_n(\UU\G)$. Using
Lemma~\ref{lem:independence}, it is easily verified that the
$L^2$-signature $
\sigma^{(2)}_\G(h)$ is well defined, i.e.\ independent of the choice of $Q$. The
following result follows easily now  from the observation that metabolic forms on
$P\oplus P^*$, which represent zero in Witt groups, have trivial
$L^2$-signature.

\proclaim{{C}orollary}\label{cor:L2-signature} The $L^2$\/{\rm -}\/signature is a well\/{\rm -}\/defined real
valued homomorphism on the Witt group of hermitian forms on finitely generated
projective
$\UU\G$\/{\rm -}\/modules. Restricting this homomorphism to nonsingular forms on free
modules gives 
$$
\sigma^{(2)}_\G:L^0(\UU\G) \ra \R.
$$
\endproclaim

As a great example, consider a finite CW-complex $X$ together with a homomorphism
$
\varphi:\pi_1X\to\G$. Then all the twisted homology groups
$H_q(X;\UU\G)$ are finitely generated projective $\UU\G$-modules.  This follows
from the fact that $\UU\G$ is a {\it  von Neumann regular} ring. In particular,
every finitely presented $\UU\G$-module is projective and these modules form an
abelian category \cite{G} (see also \cite{Re}).

If $X$ is, in addition, an oriented Poincar\'{e} complex of dimension~$4k$,
possibly with boundary, then the intersection form
$$ h_X: H_{2k}(X;\UU\G) \ra H_{2k}(X, \partial X; \UU\G) \cong H^{2k}(X;\UU\G)
\cong H_{2k}(X;\UU\G)^*
$$ is a hermitian form. The first isomorphism above is Poincar\'e duality and the
second comes from the universal coefficient spectral sequence 
$$
\Ext^p_{\UU\G}(H_q(X;\UU\G),\UU\G) \Longrightarrow H^{p+q}(X;\UU\G)
$$ which degenerates at $p=0$ by the projectivity of the twisted homology groups.

\numbereddemo{Definition} The $L^2$-{\it signature} $ \sigma^{(2)}(X, \varphi)$ is defined to be the
real number
$\sigma^{(2)}_\G(h_X)$.
\enddemo

\proclaim{Lemma}\label{lem:properties} The $L^2$\/{\rm -}\/signature has the following
properties\/{\rm :}\/
\begin{itemize}
\ritem{1.} If $(X^{4k}, \varphi)$ is the boundary of a $(4k+1)$\/{\rm -}\/dimensional
Poincar{\rm \'{\it e}} complex {\rm (}\/with the homomorphism to $\G$ extending\/{\rm )}  then $
\sigma^{(2)}(X, \varphi)=0$.
\ritem{2.} The resulting homomorphism from the bordism group of oriented Poincar{\rm \'{\it e}}
complexes
$$
\sigma^{(2)}: \Omega_{4k}^{PC}(B\G) \ra\R
$$ is equal to the ordinary signature $ \sigma_0$ on the image of the bordism
group
$\Omega_{4k}^{\rm TOP}(B\G)$ of oriented topological manifolds. 
\ritem{3.} If $(X, \varphi)$ and $(X', \varphi')$ have the same boundary {\rm (}\/and the
homomorphisms to $\G$ agree on it\/{\rm )} then
$$
\sigma^{(2)}(X \cup_{ \partial  X} X', \varphi\cup \varphi')= \sigma^{(2)}(X,
\varphi) + \sigma^{(2)}(X', \varphi').
$$
\ritem{4.} The {\rm  reduced} $L^2$\/{\rm -}\/signature $\sigma^{(2)}(X,
\varphi)- \sigma_0(X)$ of a topological $4k$\/{\rm -}\/manifold only depends on the boundary
$(\partial X, \varphi_ \partial )$.
\end{itemize}

\endproclaim
\demo{Proof}  The proofs of 1.\ and 3.\ are exactly as for the ordinary
signature. One uses the homological properties of $\UU\G$ mentioned above as well
as the usual additivity properties of the $\G$- dimension.  Property 2  for
smooth manifolds is exactly Atiyah's $L^2$-Index theorem \cite{A} applied to the
$L^2$-signature operator~$S$. One needs to check that the definition of $\sigma^{(2)}_\G(X, \varphi)$
given above agrees with Atiyah's definition involving the $L^2$-index of $S$.
This follows from the
$L^2$-Hodge theorem together with the fact that the von Neumann dimension can be
read off after tensoring with $\UU\G$. A detailed argument will be given in a
forthcoming paper by L\"uck and   Schick.

The statement for topological manifolds follows from the fact that the cokernel
of the map from smooth to topological bordism is a torsion group and we are
mapping into the torsion-free group $\R$. Finally, it is clear that Property~4 
is a direct consequence of 2  and 3. 
\enddemo

\numbereddemo{{R}emark}\label{rem:eta} If $ \partial X$ has a smooth structure then one can
pick a Riemannian metric
$g$ and define the $\eta$-invariant $\eta( \partial X,g)$ of the signature
operator. By lifting $g$ and the operator to the $\G$-cover, Cheeger and Gromov
\cite{ChG} also define the von Neumann $\eta$-invariant $\eta^{(2)}( \partial
X,\varphi_ \partial ,g)$. They show that the difference $\eta^{(2)}-\eta$ is
independent of the metric $g$. This difference is referred to as the von Neumann
$\rho$-invariant.

Moreover, if $X$ is smooth then the Index and $L^2$-Index theorems for manifolds
with boundary imply that the reduced $L^2$-signature of $(X,
\varphi)$ equals the von Neumann $\rho$-invariant of $( \partial X,
\varphi_\partial )$.  By an argument similar to that in Lemma~\ref{lem:properties},
Part 2, it follows that this equality holds true if only $ \partial X$ is smooth.
\enddemo

As an example, consider a knot $K:S^{4k-3}\hra S^{4k-1}, k > 1$. Then surgery on $K$
leads to a closed $(4k-1)$-dimensional manifold $M$, together with the
abelianization map
$\varphi:\pi_1M\to\Z$. It follows from Remark~\ref{rem:center valued} that the
corresponding center-valued von Neumann
$\rho$-invariant detects the concordance group modulo torsion.

\proclaim{Lemma}\label{lem:Ranicki} Let $R \sigma$ be Ranicki\/{\rm '}\/s symmetric signature
map. Then the composition
$$
{\Omega^{PC}_{4k}(B\G)}\stackrel{R \sigma}{\longrightarrow} {L^0(\Z\G)}\longrightarrow 
{L^0(\UU\G)} \stackrel{\sigma^{(2)}_\G}{\longrightarrow} {\R}
$$ 
is equal to the $L^2$-signature from Lemma~{\rm \ref{lem:properties},} Part {\rm 3.}
\endproclaim

\demo{Proof}  The result follows from the fact that for von Neumann regular rings
$R$ the bordism class of a chain complex gives the same element in
$L^0(R)$ as the corresponding intersection form on the middle homology. This
follows from the chain homotopy invariance of algebraic Poincar\'e cobordism
(Chapter~1 of \cite{Ra1}) 
 and algebraic surgery below the middle dimension (Chapter~4 of \cite{Ra1}).
\enddemo 

 The next result is not strictly necessary for the definition of our
invariants but it seems appropriate to mention it at this point. 
 In addition to reproving that the $L^2$-signature gives a well-defined slice
obstruction via Theorem~\ref{invariant thm} it also shows that in order to define
our obstructions one can equally well work with
$(n)$-solutions $W$ which are finite Poincar\'{e} complexes (rather than
topological
$4$-manifolds).

Note that by a theorem of Higson-Kasparov \cite{HK} the following assumption is
satisfied for torsion-free amenable groups, and hence in particular for PTFA
groups.

\proclaim{Proposition}\label{prop:Baum-Connes} If $ \G $ is torsion\/{\rm -}\/free and the analytic
assembly map
$$
A_ \G:K_0 ({B} \G) \to K_0(C^{\ast} \G)
$$
 is onto then
$\sigma^{(2)}_ \G(h) =
\sigma_0(\varepsilon_{\ast} h)$ for all $h\in L^0(C^*\G)$. In particular{\rm ,} the
$L^2$\/{\rm -}\/signature from Lemma~{\rm \ref{lem:Ranicki}} equals the ordinary signature on all
of
$\Omega^{PC}_{4k}(B\G)$.
\endproclaim 

\demo{Proof}   Since $h$ is invertible, it follows that $0\not\in \spec(h)$ and
since
$\spec(h)$ is compact, it actually has a gap around $0$. Therefore, the
characteristic functions $p_+$ and $p_-$ are {\it  continuous} functions on
$\spec(h)$ and $p_+(h), p_-(h) \in M_n(C^{\ast} \G)$.  By taking the difference
we get the signature map
$\sigma_ \G:L^0(C^{\ast} \G)\to K_0(C^{\ast} \G)$ which is an isomorphism for any
$C^*$-algebra.  Recall that
$\sigma^{(2)}_
\G(h)$ is given by composing with
$$  K_0(C^{\ast} \G) \longrightarrow K_0( \NN \G)\stackrel{\tr_ \G}{\lrar} \R.
$$
  Moreover, by Atiyah's $L^2$-Index theorem applied to all twisted Dirac
operators, the two compositions
$$ \begin{array}{l}
{K_0({B} \G)}
\stackrel{A_ \G}{\lrar} {K_0(C^{\ast} \G)}
\stackrel{\tr_ \G}{\lrar}{\R}, \\[8pt]
{K_0({B} \G)}{\lrar}
{K_0({\ast})=K_0(\C)}{\displaystyle\mathop{\stackrel{\tr_1}{\lrar}}_{\cong}} {\Z}
\end{array}
$$ are the same (compare \cite[7.15]{BCH}).  The claim now follows from
surjectivity of the analytic assembly map and the naturality of the signature and
assembly maps in the following commutative diagram.
$$
\begin{array}{ccl}
{L^0(C^{\ast} \G)}  &\lrar&{L^0(\C)}\\[5pt]
 {\sigma_ \G}  \!\Big\downarrow\!&&
 { \sigma_1}\!\Big\downarrow\!{ \cong} \\[5pt]
{K_0(C^{\ast} \G)} &\lrar& {K_0(\C)} \\[5pt]
{A_ \G}\!\Big\uparrow \!{ \cong}&&{A_1}\!\Big\uparrow\!{ \cong}\\[5pt]
{K_0({B} \G)} &\lrar& {K_0({\ast})}.\\
\noalign{\vskip-16pt}
\end{array} 
$$ 
\enddemo
 \eject

We next show how to define the  
$L^2$-signature for forms over the quotient field $\KK$ of $\Z\G$, if it exists. 
More generally, for any group $ \G$ we can consider the following diagram of
inclusions of rings with involution:
$$
\begin{array}{ccc} 
{\C\G}&\lrar& {\NN \G}\phantom{.}\\[5pt]
\Big\downarrow&&\Big\downarrow \\[5pt]
{\DD\G} &\lrar& {\UU\G} .
\end{array}
$$ Here the {\it  division closure} ${\cal D} \G$ of $\C \G$ in
$ \UU \G$ is the smallest intermediate ring which is division closed.  This means
that if $r\in {\cal D} \G$ is invertible in $ \UU \G$ then the inverse
$r^{-1}$ already lies in ${\cal D} \G$.  For the case $ \G = \Z$ we obtain
${\cal D} \G = \C(t)$, the quotient field of rational functions on $S^1$.  In
fact, if $ \C\G$ satisfies the Ore condition, then $\DD\G$ is the Ore
localization of $\C\G$
\cite{Re} and we have constructed the $L^2$-signature
$$
\sigma^{(2)}_ \G: L^0(\KK) \ra  L^0(\DD\G)
\longrightarrow L^0( \UU \G) \longrightarrow \R.
$$
 This applies in particular to PTFA groups for which this $L^2$-signature
equals the ordinary signature $\sigma_0$ on the image of $L^0(\Z \G)$ in
$L^0(\KK)$ by Proposition~\ref{prop:Baum-Connes}.

We conclude this section with an innocent looking but extremely useful property.  

\proclaim{Proposition} \label{prop:naturality of trace} For a subgroup $ \G_1
\subseteq 
\G_2${\rm ,} there are commutative diagrams
$$
\BoxedEPSF{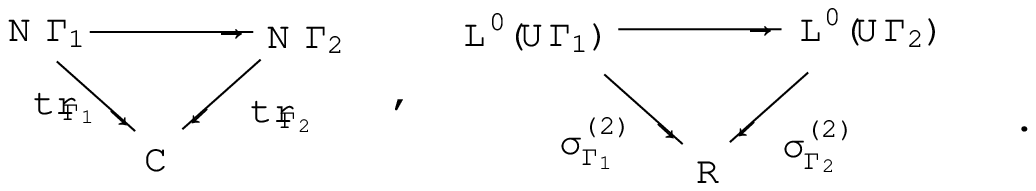 scaled 1000}
$$
\endproclaim 

\demo{Proof}   The most difficult part is to construct the homomorphism $\UU\G_1
\longrightarrow \UU\G_2$. A homomorphism $\NN\G_1 \to \NN\G_2$ is given by
completing 
$$
 a\otimes \id\in \End(\l^2  \G_1  
\otimes_{ \C \G_1} \C  \G_2)
$$ 
to a bounded operator on $\l^2  
 \G_2$ for any $a\in \mathcal N  \G_1$. Since
$$
\langle (e_1\otimes e_2)(a\otimes\id),e_1\otimes e_2 \rangle =
\langle (e_1)a,e_1 \rangle \cdot \langle e_2,e_2 \rangle=\langle (e_1)a,e_1\rangle
$$ 
it follows that the first diagram commutes. For details see
\cite[Thm.~3.3]{Lu1}. This reference also contains the statement that tensoring an
$\NN\G_1$-module with $\NN\G_2$ is a flat functor. In particular, the map
$\NN\G_1 \to \NN\G_2$ sends nonzero divisors to nonzero divisors and thus
induces a homomorphism $\UU\G_1 \to
\UU\G_2$. To see that the second \pagebreak  diagram above commutes just observe that
``diagonalizing'' a hermitian matrix over $\UU\G_1$ and then tensoring up the
$\pm 1$-eigenspaces to $\UU\G_2$ diagonalizes the induced matrix over
$\UU\G_2$. Thus the commutativity of the first diagram proves the claim.
\enddemo 

We should warn the reader that the von Neumann algebra
$ \NN \G$ is not functorial in $ \G$. This has to do with the specific choice of
the Hilbert-space $\l^2 \G$. Proposition~\ref{prop:naturality of trace} gives the
best possible functoriality which is valid for all groups.  If
$ \G$ is amenable, then the equality of the reduced and maximal\break
$C^{\ast}$-algebras (which are functorial!) implies that the   projection $\G
\sra \{1\}$ induces a homomorphism of $C^{\ast}$-algebras
$\varepsilon:C^{\ast} \G \to \C$.  For example, if $ \G = \Z$ then this is given
by evaluating a continuous function at $1\in S^1$.  This clearly does {\it  not}
extend to
$ \NN\Z =L^\infty(S^1)$.

\numbereddemo{{R}emark}\label{rem:Z}
The above proposition is fundamental to all our calculations! We will construct our
knots in such a way that the relevant intersection form over $\Z\G$ will contain as
its entries only linear combinations of powers of a single nontrivial group element
$\eta\in\G$. Since our groups $\G$ are torsion-free, this gives an inclusion of
groups 
$$
\Z\cong \langle \eta \rangle  \subset \G
$$
to which  Proposition~\ref{prop:naturality of trace} can be applied. Thus the
$L^2$-signature for
$\G$ can be calculated for this particular hermitian form as an integral over the
circle of certain twisted signatures. The concrete example   to be
used can be found after Definition~\ref{def:signature}.
\enddemo

\vglue-14pt
\section{Nonslice knots with vanishing Casson-Gordon
invariants}\label{sec:example}
\vglue-6pt 
In this chapter we give examples of knots which are
$(2)$-solvable but not $(2.5)$-solvable. In particular, by Theorem~\ref{vanish}
these are the first examples of knots that have vanishing metabelian concordance
invariants, including Casson-Gordon invariants, but are not slice knots even in
the topological category. Only the highly technical Proposition~\ref{uniqueP}
prevents us from exhibiting $K_n$, for each $n\geq 2$, which is
$(n)$-solvable but not
$(n.5)$-solvable. We focus on one example and indicate how this may be modified
to produce an infinite family of such.

Consider the knot $K$ in Figure~6.1.
The rectangles containing integers symbolize full twists between the two strands
which pass vertically through the rectangles. Thus the rectangle labeled $-2$
symbolizes two left-handed full twists (see \cite{Ki}). The rectangle labeled by
$J^*$ symbolizes the four component string link obtained by taking four
untwisted parallel copies of the knotted arc $J*$, which is shown at the bottom of
Figure~6.1  . We shall show that $K$ is $(2)$-solvable but that it fails to
satisfy Theorem~\ref{invariant thm} for $(2.5)$-solvability,  as detected by the
reduced
$L^2$-signature \pagebreak of Section~\ref{sec:L2}. The same proof  will show\break

\vglue-18pt
\centerline{\BoxedEPSF{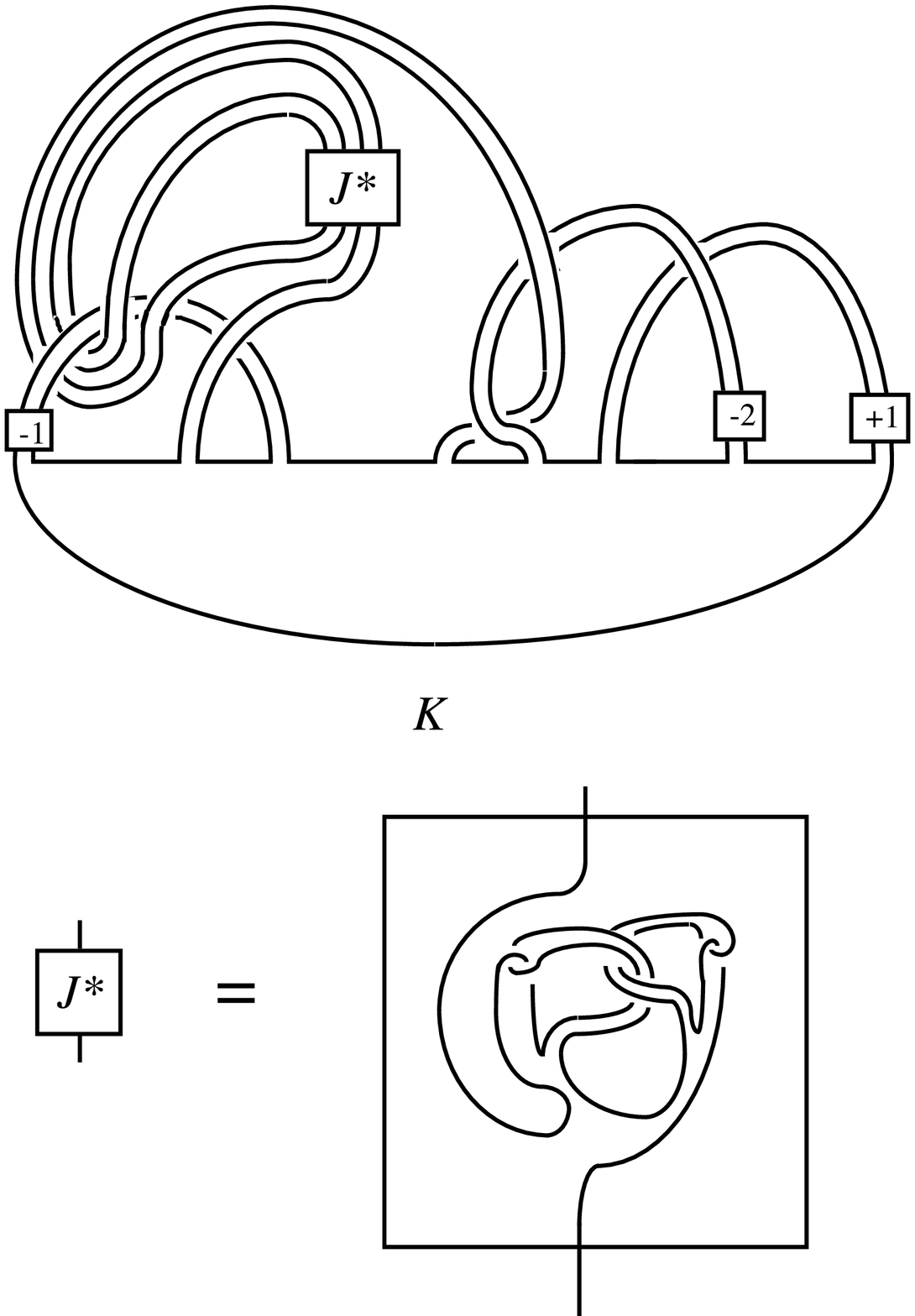 scaled 330}}
   \vglue2pt
  \centerline{Figure 6.1}
 
 \vglue-36pt
\centerline{\BoxedEPSF{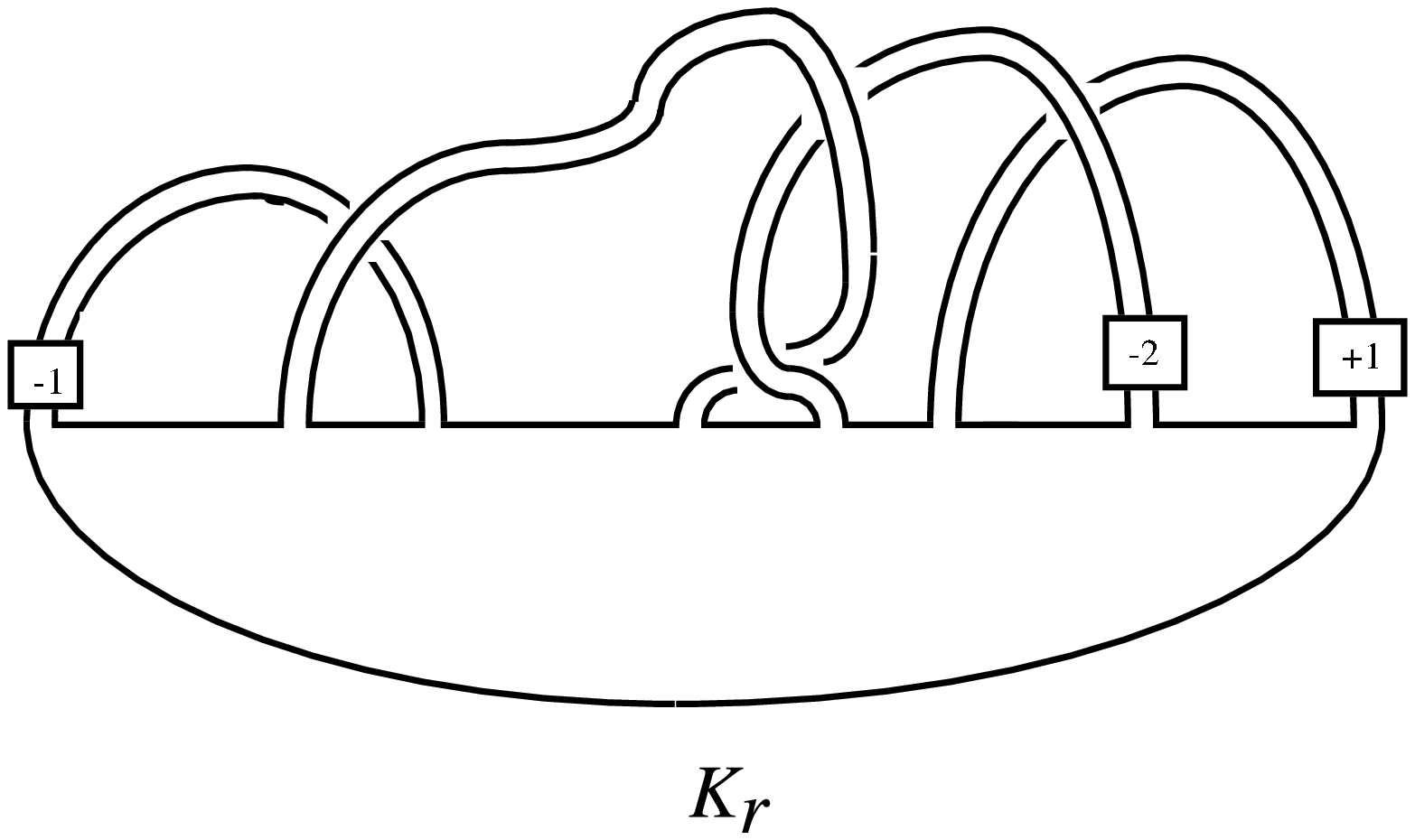 scaled 450}}
\vglue3pt
 \centerline{Figure 6.2}
 \vglue12pt

\noindent  
 that $M$, the
zero surgery on $K$, is not rationally $(2.5)$-solvable with multiplicity $1$
(see the definition above Theorem~\ref{invariant thm}). Consequently,  $K$ is not
slice in any rational homology ball wherein the meridian generates the free part
of $H_1$ of the slice disk complement.

First we sketch the argument that $K$ is $(2)$-solvable but not
$(2.5)$-solvable. The bulk of the work is to find a fibered ribbon knot $K_r$
(see Figure~6.2) which is $(2)$-solvable ``in only one way''; i.e.,  for
which $\AA_0$ and $\AA_1$ have unique self-annihilating submodules. For example
we ask the reader to check that the ordinary Alexander module of $K_r$ is cyclic
of order $p(t)^2$ where
$p(t)=t^{-1}-3+t$, the Alexander polynomial of the figure~8 knot. Since $p(t)$ is
irreducible it follows that this module contains a {\it unique proper
submodule\/}. Thus the rational Alexander module $\AA_0$ certainly contains a
unique submodule $P_0$ which is self-annihilating with respect to
$B\l_0$. Then we form the knot $K$ by modifying the ribbon knot by surgeries on
two circles, in such a subtle way that $\AA_0$ and $\AA_1$ are unaffected and $K$
remains $(2)$-solvable, but such that the $(2)$-solution for $K$ has some
nontrivial second homology.  Our proof will proceed by contradiction. Suppose
$K$ were $(2.5)$-solvable via $W^{\prime}$. Let $W_r$ denote the complement of
the ribbon disk $K_r$ and $Y$ denote the cobordism from $M_r$, the zero surgery
on $K_r$, to $M$, which consists of two relative 2-handles. Let $W$ denote the
union of $W_r$ and $Y$ along $M_r$. We will show that $W$ is a
$(2)$-solution for $M$. Of course $W^{\prime}$ is also a $(2)$-solution for $M$.
By Theorem~\ref{invariant thm} applied to each, there exist representations of
$\pi_1(M)$ into $\G_2$ which extend to $\pi_1(W)$ and $\pi_1(W^{\prime})$
respectively. We use the fact that there are   unique self-annihilating 
submodules
to show that these representations coincide. Let this common map be denoted
$\phi_2$. Since
$W^{\prime}$ is a $(2.5)$-solution, $B(M,\phi_2)=0$ by Theorem~\ref{B=0}.
 But $W$
may also be used to calculate $B(M,\phi_2)$.  The intersection form of $W$ is
represented by a simple $(2 \times 2)$-matrix whose 
$L^2$-signature was calculated to be nonzero in Section~\ref{sec:L2}. This
contradiction will then finish the proof.

As stated in the introduction, if there exists a fibered genus two ribbon knot
for which $\AA_0,\dots,\AA_{n-1}$ all have unique self-annihilating submodules,
then the same procedure we discuss herein creates a knot which is $(n)$-solvable
but not $(n.5)$-solvable.

Now we describe in detail the construction of $K$. Consider
$K'=J\#(-J)$ where $J$ is the figure-eight knot (as shown in Figure~6.3).
This is a well-known fibered ribbon knot. We summarize the argument. Consider a
knotted ball pair
$J'=(B^3,B^1)$ representing the figure-eight knot, and cross this with $[0,1]$.
The result is a knotted $2$-disk $ \Delta=J'\times [0,1]$ in
$B^4$ whose boundary is the ribbon knot $K'$.  Now, $S^3-J=B^3-J'$ is known to be
fibered with fiber the standard Seifert surface (a punctured torus $T$). It
follows that
$B^3\times [0,1]-J'\times [0,1]$ is fibered with fiber
$T\times [0,1]$, a genus $2$ handlebody $H$, and hence that
$K'$ is a genus 2~fibered knot. 
 
\centerline{\BoxedEPSF{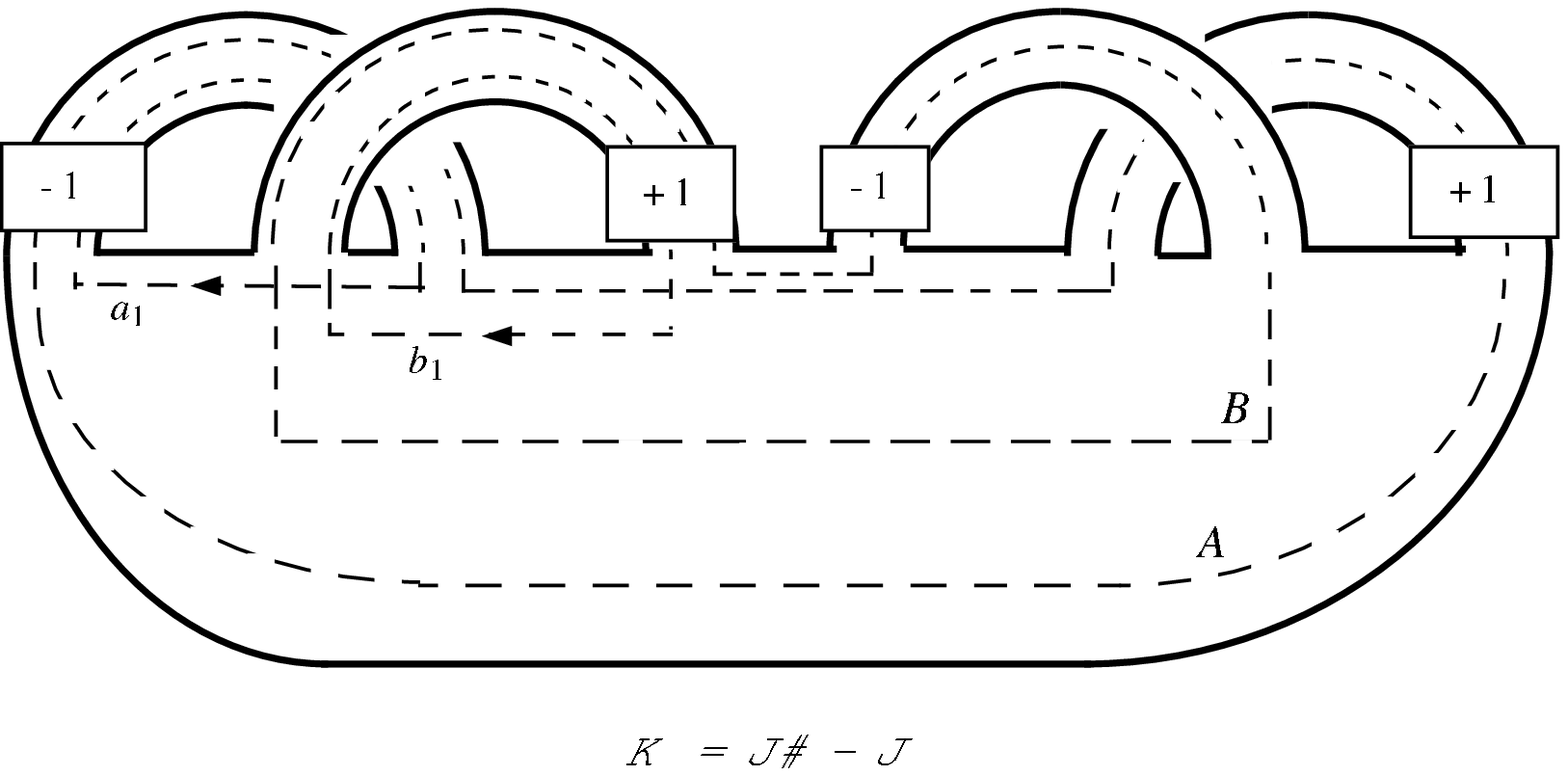 scaled 550}}
 \vfil 
\centerline{Figure 6.3}
 \eject

\centerline{\BoxedEPSF{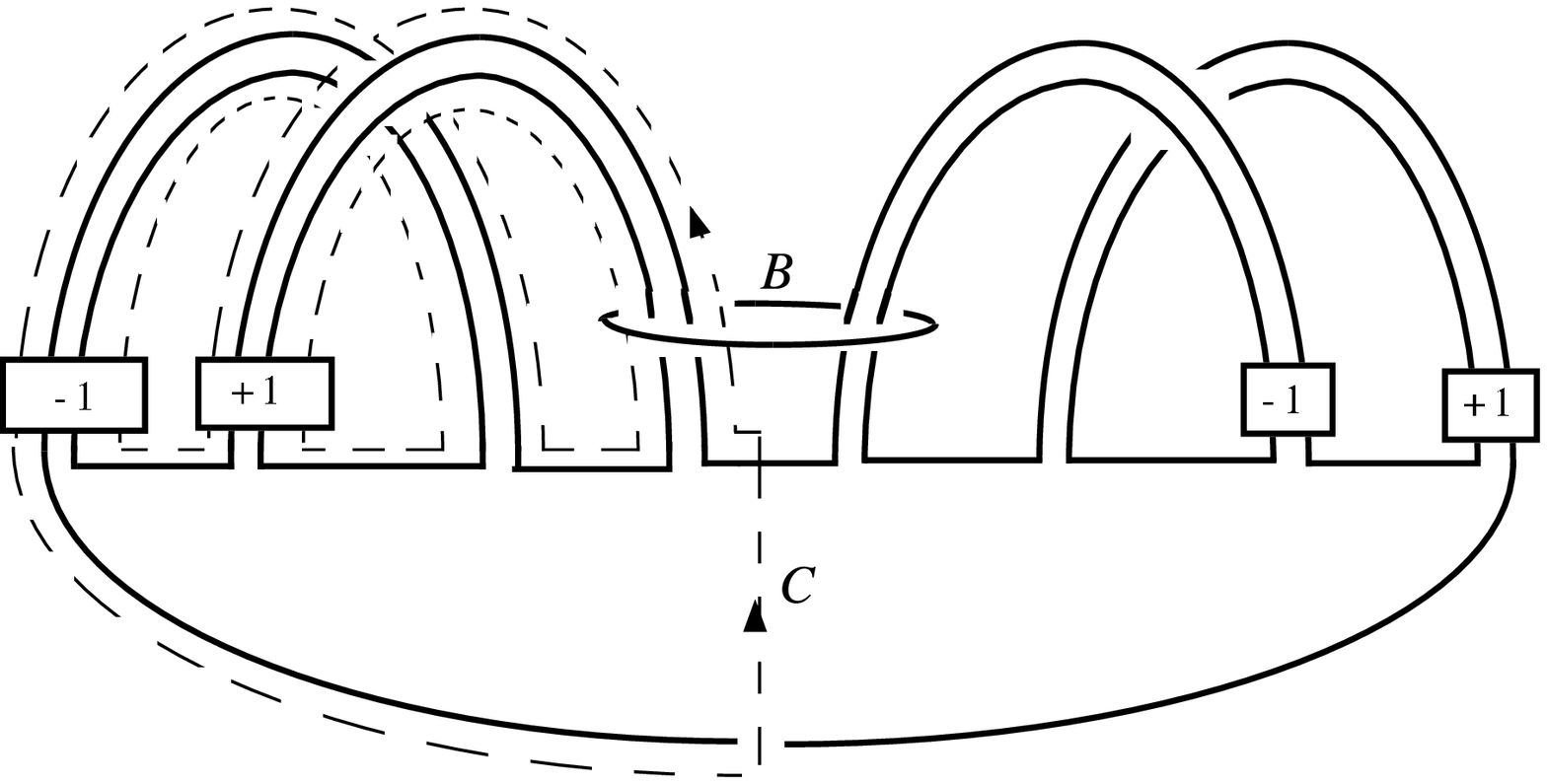 scaled 450}}

    \centerline{Figure 6.4}
\vglue6pt

Moreover from this point of view it is easy to see that the loops labeled $A$ and
$B$ in Figure~6.3 bound embedded disks in $H$, and that the loops labeled
$a_1$ and
$b_1$ map to generators $x$ and $y$ of the fundamental group
$F\left<x,y\right>$ of $H$. Let $f:T \to T$ be the monodromy homeomorphism for
$J$,
$f'$ that of $K'$ and $\tilde f':H \ra H$ that of $  \Delta$. Then we may assume
that
$f$ preserves
$\partial T$. It follows that $f'$ preserves the ``sub-longitude'' loop labeled
$C$ in Figure~6.4.

Note that $B$ is unknotted in $S^3$ and has self-linking
zero on the obvious Seifert surface. Therefore if we perform $+1$ surgery on
$B$, $K'$ will be transformed to a new knot
$K_r$ as shown in Figure~6.2. This may be seen by pushing $B$ off the
Seifert surface, as shown in Figure~6.4, and ``blowing down'' $B$ by one
application of Kirby's calculus \cite{Ki}.
It is known that the result, $K_r$, of such a modification is again a fibered
knot whose monodromy $f_r$ equals $D_B\circ f'$ where $D_B$ is a Dehn-Twist along 
$B$
(\cite{H}, \cite{Sta}). Since $B$ bounds a disk in $H$, $D_B$ extends to $\widetilde
D$ on
$H$ and $f_r^{\phantom{|}}$ extends to $\tilde f_r:H\to H$. Therefore
$K_r$ is also a fibered ribbon knot with fibered ribbon disk
$  \Delta_r$ and fiber $H$. Moreover $\widetilde D$ is homotopic to the identity
so
$\tilde f_r$ is homotopic to
$\tilde f':H\to H$. Thus $B^4-  \Delta_r$ is homotopy equivalent to $B^4- 
\Delta$. In particular note that the element $[x^{-1},y]$ of $\pi_1(H)$ is the
image of the sub-longitude $C=[b^{-1}_1,a_1]$ and thus 
$$ (\tilde f_r)_{\ast}([x^{-1},y])=(\tilde f')_{\ast}([x^{-1},y])=[x^{-1},y]
$$ in $\pi_1(H)$. Moreover $(f_r)_{\ast}(C)$ is represented by the image of $C$
under
$D_B$.

Finally we will modify $K_r$ by two surgeries, resulting in
$K$. The effect of these surgeries is subtle enough that
$\AA_0$ and $\AA_1$ as well as the Casson-Gordon invariants are unchanged (as we
shall see). Consider an embedded circle $\eta$ in the complement of the obvious
Seifert surface for $K_r$. The specific example we wish to consider is shown in
Figure~6.5, but to find examples of knots which are
$(n)$-solvable but not
$(n.5)$-solvable one would choose $\eta$ to represent a nontrivial class in the
$n^{\rm th}$ derived group of
$\pi_1(W_r)$. This $\eta$ was also chosen so that $j_{\ast}(\eta)=C$,
 which will later be shown to generate $\AA_1(W_r)$. Note that $\{A,B,\eta\}$ is the Borromean ring.

\centerline{\BoxedEPSF{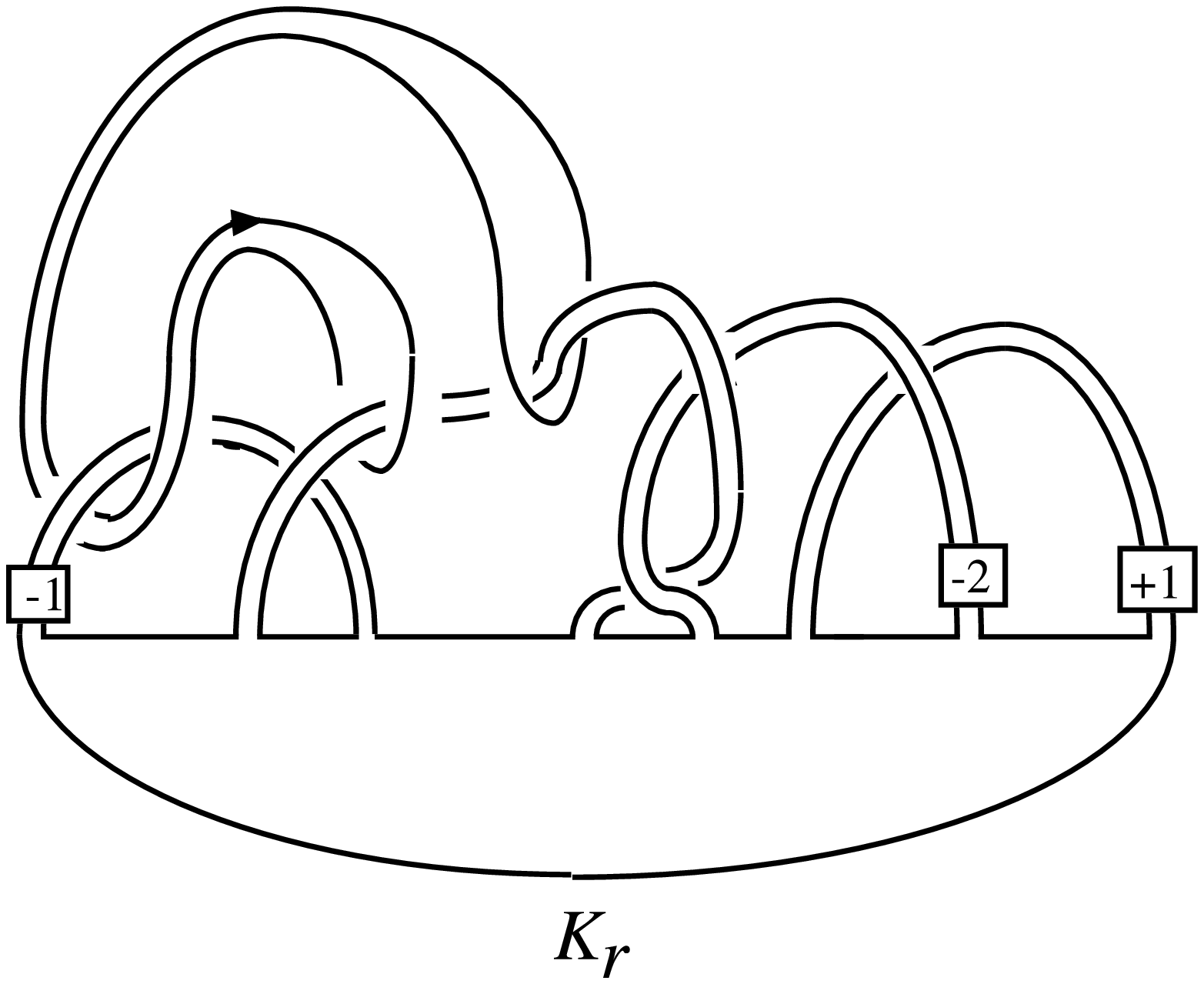 scaled 450}}
\vglue12pt
\centerline{Figure 6.5}
\vglue24pt

\figin{Fig66}{450}
\vglue12pt
\centerline{Figure 6.6}
\vglue12pt

 Now replace a
solid torus neighborhood of $\eta$ by the
$3$-manifold shown in Figure~6.6. This manifold is the result of two Dehn
surgeries on $S^1\times D^2$. Since
$\{\gamma_1,\gamma_2\}$ forms a zero-framed Hopf link in
$S^3$ (ignoring $K_r$), and the result of such a surgery is known to be
homeomorphic to
$S^3$, the image of $K_r$ under this homeomorphism is a new knot $K$ in $S^3$
which is shown in Figure~6.1  . The reader proficient in Kirby's calculus may
confirm the accuracy of Figure~6.1   by first isotoping
$\{\gamma_1,\gamma_2\}$ until it looks more like a Hopf link, next ``sliding''
all strands of $K_r$ which ``pass through
$\gamma_1$'' over $\gamma_2$, then ``sliding'' strands of
$K_r$ over $\gamma_2$ until the Hopf link becomes split off from the knot
\cite{Ki}. We remark that the $3$-manifold of Figure~6.6 is a knot complement,
namely the complement of the knot $J^{\ast}$ obtained by closing up the knotted
arc shown at the bottom of Figure~6.1  . Using the trefoil knot in place of
$J^{\ast}$ would lead to a simpler example of a knot with vanishing Casson-Gordon
invariants which is not
$(2.5)$-solvable, but it is difficult to see if it is
$(2)$-solvable. However, in place of $J^{\ast}$, we could use any Arf invariant
zero knot, such that the integral of the Levine signature function is nonzero,
and reach an identical conclusion.  This will be detailed in a forthcoming paper. 

We will now show that $K$ is $(2)$-solvable, using the fact that $K$ is obtained
from the ribbon knot $K_r$ by performing two surgeries. Let $W_r$ denote the
exterior in
$B^4$ of the fibered ribbon disk for $K_r$. Let $W$ denote the
$4$-manifold obtained from $W_r$ by adding $2$-handles along the zero-framed
circles
$\{\gamma_1,\gamma_2\}$. Then $W$ is an $H_1$-bordism for $M$, the zero surgery
on $K$, and we will show that it is in fact a $(2)$-solution. Note that since
$\{\gamma_1,\gamma_2\}$ are null-homotopic in $M$,
$W\backsimeq W_r\vee S^2\vee S^2$ and so
$\pi_1(W)\cong\pi_1(W_r)$ and $H_2(W; \Z)\cong\Z^2$ are generated by $\{\l_i\}$. The
integer-valued intersection form with respect to this basis is the standard
hyperbolic form. The circles $\gamma_i$ bound obvious immersed disks $D_i$ in a
neighborhood of
$\eta$. It is desirable to introduce a local kink in each, as shown in
Figure~6.6, so that the push-off of $\gamma_i$ into $D_i$ has linking
number zero with
$\gamma_i$. These disks, together with the cores of the
$2$-handles form immersed $2$-spheres $\l_1$ and $\l_2$. Being
$1$-connected, these surfaces lift to any cover. If
$[\eta]\in\pi_1(W_r)^{(n)}$, we shall show that $L=\{\l_1\}$ generates an
$n$-Lagrangian and $\l_2$ is its $n$-dual. Since the particular $\eta$ of
Figure~6.5 lies in
$\pi_1(M_r)^{(2)}$ (it bounds a surface in the complement of a Seifert surface
for $K_r$), this will show that $K$ is
$(2)$-solvable. But we show the more general fact to illustrate how easy it is to
generalize $K$ to an
$(n)$-solvable example. It suffices to show that, with
$\pi_1(W)\cong\pi_1(W_r)$ coefficients,
$\mu(\l_1)=j_{\ast}(\eta)-1$ and $\lambda(\l_1,\l_2)=1$ where
$j_{\ast}:\pi_1(M)\ra\pi_1(W_r)$ because then, with
$\pi_1(W_r)/\pi_1(W_r)^{(n)}$ coefficients,
$\mu(\l_1)=\lambda(\l_1,\l_1)=0$. But this is clear from an analysis of the two
points of self-intersection of $D_1$ and the one point of intersection of $D_1$
with $D_2$. Section~\ref{sec:surfaces} contains a detailed explanation of how to
calculate $\mu$ and $\lambda$. Thus $M$ is
$(n)$-solvable via $W$ ($(2)$-solvable in our special case).

Now suppose that $K$ is $(2.5)$-solvable via a $4$-manifold $W'$. By
Theorem~\ref{invariant thm} we can  find a representation
$\phi_2:\pi_1(M)\ra\G_2^U$  such that
$\phi_2$ extends to
$W'$ and $B(M,\phi_2)=0 = \sigma(W')$. (Actually we shall use Remark~\ref{post inv rem}.3 to
restrict the image of $\phi_2$ to a certain subgroup $\G_2$ of our universal
group $\G^U_2$.) If we can show that $\phi_2$ also extends to $W$ and that
$\phi_2(\eta)\neq1$, we can quickly reach a contradiction as follows. Since
$\phi_2$ extends to W, we can calculate $B(M,\phi_2)$ using W. Note that
$H_2(W_r; \KK_2)\cong H_1(W_r;
\KK_2)=0$ by Proposition~\ref{rank coef} and Proposition~\ref{torsion}, so that $H_2(W; \KK_2)\cong H_2(W,W_r;
\KK_2)$ is free on $\{\l_1,\l_2\}$. (In fact, Lemma~\ref{2comp} implies
$H_2(W_r;\Z\G_2)=0$ so that even $H_2(W;\Z\G_2)$ is free.) The intersection and
self-intersection forms with $\Z\G_2$ coefficients may be computed, by
naturality, from the intersection and self-intersection forms with $\Z\pi_1(W)$
coefficients derived above. Let $\psi_2$ denote the extension to
$\pi_1(W)$ and let $t=\phi_2(\eta)=\psi_2\circ j_{\ast}(\eta)$. Note that
$\mu(\l_1)=\psi_2\circ j_{\ast}(\eta)-1=t-1$ determines that
$\lambda(\l_1,\l_1)=t+t^{-1}-2$ by Property 5 of Definition~\ref{def:form} below
(alternatively this may be computed directly as above). Thus $B(M,\phi_2)\in L^0(
\KK_2)$ is represented by the matrix 
$$
\left(\begin{array}{cc} t+t^{-1}-2  &1\\ 1  &t+t^{-1}-2\end{array}\right).
$$ We claim that the reduced
$L^2$-signature of  $B(M,\phi_2)$ (see Lemma~\ref{lem:properties} and the discussion above
Definition~\ref{Qsolv}) is nonzero. For, since
$t\neq1$, the subgroup of $\G_2$ generated by $t$ is infinite cyclic. Since the
matrix above also represents an element of $L^0(\C (t))$, applying
Proposition~\ref{prop:naturality of trace}, we see that $\sigma^{(2)}_{\G_2}$
agrees with  $\sigma^{(2)}_{\Z}$. But the latter was calculated to be nonzero
below Definition~\ref{def:signature}. Finally, note that the ordinary signature of the
above matrix is zero so the reduced and unreduced $L^2$-signatures agree. Thus
$B(M,\phi_2)\neq0$. This contradiction will complete the proof that $K$ is not
$(2.5)$-solvable. The remainder of this section is devoted to verifying that if
$\phi_2$ is the (essentially unique) map guaranteed by Theorem~\ref{invariant
thm} (applied to $W^\prime$), then $\phi_2(\eta)\neq1$ and  $\phi_2$ extends to
$W$. Note that since $W$ is also a $(2)$-solution, Theorem~\ref{invariant thm}
applied to $W$ implies that {\it  certain} maps extend to $W$. From this point of
view we must show that one of these can be chosen to coincide with $\phi_2$.

Let $\phi_0:\pi_1(M)\ra\Z\equiv\G_0$ be the unique homomorphism sending a
meridian to 1. Since both $W'$ and $W$ are rational
$H_1$-bordisms with multiplicity $1$,
$\phi_0$ extends uniquely to $\psi'_0$ and
$\psi_0$ respectively. By Theorem~\ref{self annihil} with
$n=1$ and $\RR_0=\Q[t^{\pm 1}]$, since $M$ is $(1)$-solvable via $W'$ and $W$,
$B\l_0$ is hyperbolic and the kernels of the inclusion maps
$$ H_1\left(M;\Q[t^{\pm 1}]\right)\xrightarrow{i_{\ast}}  H_1\left(W';\Q[t^{\pm
1}]\right)
$$ and 
$$ H_1\left(M;\Q[t^{\pm 1}]\right)\xrightarrow{j_{\ast}} H_1\left(W;\Q[t^{\pm
1}]\right)
$$ are self-annihilating. But, as mentioned earlier, $\AA_0(K_r)$ has a unique
self-annihilating submodule $P_0$. Now, it is easy to see that $\AA_0(K) \cong \AA_0(K_r)$ by observation that, since any
loop on the obvious Seifert surface for $K_r$ has zero linking number with the $\gamma_i$, the Seifert matrix is
unaffected by the surgeries. Thus $\Ker i_{\ast}=\Ker j_{\ast}$. Choose a nonzero element
$p_0\in P_0$ inducing $\phi_1:\pi_1(M)\ra\G^U_1$ by Theorem~\ref{characterizing
char}. By Theorem~\ref{extend}
$(n=1)$, $\phi_1$ extends to $\psi'_1$ and $\psi_1$ on
$\pi_1(W')$ and $\pi_1(W)$ respectively. Before proceeding, we want to replace
$\G^U_1$ by a much smaller group in order to simplify a subsequent calculation
(\ref{uniqueP}). The basic point is that we can replace $\G^U_1$ by a subgroup
containing the images of $\psi'_1$ and $\psi_1$ and proceed with the argument. In
fact we can do slightly better. Let $S$ be the smallest direct summand of
$H_1\left(W';\Q[t^{\pm 1}]\right)$ which contains the image of
$i_{\ast}$. Since
$$ H_1\left(M;\Q[t^{\pm 1}]\right) \cong \Q[t^{\pm 1}]/p(t)^2
$$ and $P_0=\left<p(t)\right>$, and the kernel of $i_{\ast}$ is
$P_0$, the image of $i_{\ast}$ is cyclic of order $p(t)$. Since $p(t)$ is
irreducible, it follows that
$S\cong\Q[t^{\pm 1}]/p(t)^m$ for some positive integer $m$, and we can choose the
isomorphism so that
$i_{\ast}(1)=p(t)^{m-1}$. Since
$$ (\phi_1)_{\ast}:H_1\left(M;\Q[t^{\pm 1}]\right)\ra\Q(t)/\Q[t^{\pm 1}]
$$ has kernel {\it precisely\/} $P_0$ (by Theorem~\ref{characterizing char} and
since $p_0$ generates $P_0$), 
$$ (\psi'_1)_{\ast}:S\ra\Q(t)/\Q[t^{\pm 1}]
$$ is an embedding. Let $S$ also denote the image of this map. Therefore, if 
$$ H_1\left(W';\Q[t^{\pm 1}]\right)=S \oplus T
$$ we can replace
$\Q(t)/\Q[t^{\pm 1}]$ by the subgroup $S$, replace $\G^U_1$ by
$\G_1=S\rtimes\G_0$, and replace $(\psi'_1)_{\ast}$ by the projection onto $S$.
The map
$(\phi_1)_{\ast}$ is replaced by
$i_{\ast}$ followed by projection, still denoted
$(\phi_1)_{\ast}$. It remains to note that the image of
$(\psi_1)_{\ast}$ also lies in $S$. This is clear because
$\pi_1(W)\cong\pi_1(W_r)$ and $W_r$ is a ribbon disk exterior so that $j_{\ast}$ is
surjective on $\pi_1$ and hence on first homology. Moreover we can identify
$H_1\left(W;\Q[t^{\pm 1}]\right)$ with
$\Q[t^{\pm 1}]/p(t)$ in such a way that $(\psi_1)_{\ast}$ is the standard
embedding.

By Theorem~\ref{characterizing char}, these compatible ``characters'' induce
actual compatible homomorphisms
$\phi_1$, $\psi'_1$, $\psi_1$ from $\pi_1(M)$,
$\pi_1(W')$, $\pi_1(W)$ respectively to $\G_1$. Set
$$
\RR_1=\left(\Q[S]-\{0\}\right)^{-1}\Q\G_1=\K_1 [\mu^{\pm 1}]
$$ as in Definition~\ref{def:univ gr} and Corollary~\ref{Knmu}. The first of these
coefficient systems defines $\AA_1=H_1(M;\K_1[\mu^{\pm 1}])$. It will now suffice
to prove the following facts about $\AA_1$. Note that since $\eta\in \pi_1(M)^{(2)}$, it lifts to the $\G_1-$cover and hence represents a class in $\AA_1$.

\proclaim{Proposition}\label{uniqueP} $\AA_1$ contains a unique proper submodule and hence
a unique submodule $P_1$ such that
$P^\bot_1=P_1$. Moreover there exists $p_1\in P_1$ such that
$B\l_1(p_1,\eta)\neq0$.
\endproclaim 

Before embarking on the proof of Proposition~\ref{uniqueP}, we will show how it
completes the proof that $K$ is not $(2.5)$-solvable. Since $M$ is $(2)$-solvable
via $W$ and via $W'$, Theorem~\ref{self annihil} applies with $n=2$, $\G=\G_1$,
$\RR=\RR_1$ to show that the kernels of the maps
$$ j_{\ast}:H_1(M;\RR_1)\ra H_1(W;\RR_1)
$$ and
$$ i_{\ast}:H_1(M;\RR_1)\ra H_1(W';\RR_1)
$$ are self-annihilating and hence each equal to $P_1$ by
Proposition~\ref{uniqueP}. Choose $p_1\in P_1$ as guaranteed by
Proposition~\ref{uniqueP}. This induces $(\phi_2)_{\ast}:\AA_1\ra \KK_1/\RR_1$ and
$$ 
\phi_2:\pi_1(M)\ra\G_2\equiv (\KK_1/\RR_1)\rtimes\G_1
$$ where $\KK_1$ is the quotient field of $\RR_1$. Apply Theorem~\ref{extend}.1
to both $W$ and $W'$ for $n=2$ and $x=p_1$ to conclude that
$\phi_2$ extends to both
$\pi_1(W)$ and $\pi_1(W')$. Since $W^{\prime}$ is assumed to be a
$(2.5)$-solution, Theorem~\ref{B=0} with $n=2$ implies that $B(M,\phi_2)=0$.
Moreover we claim that
$\phi_2(\eta)\neq1$ since $(\phi_2)_{\ast}(\eta)= B\l_1(p_1,\eta)\neq0$. Thus
$\phi_2$ is the desired coefficient system which leads to a contradiction as
explained above. 

To prove Proposition~\ref{uniqueP}, we must compute $\AA_1$. We shall do this in
two independent ways --- first by using a few general principles and second by
explicitly computing the monodromy of the fibered knot $K_r$. Recall that
$\G_1=S\rtimes\Z$ and $\K_1$ is the (commutative) quotient field of the group
ring $\Q S$. If
$z\in\K_1$, let
$z^{\mu}$ denote the image of $z$ under the automorphism $\mu$ of $\K_1$ (see
Section~\ref{sec:Alexander}).

Upon first glance at the form of $\AA_1$ in part $b$ below, one might conclude
that it had at least two proper submodules if $k\neq1$. However remember that,
although $\K_1$ is commutative, the ring $\K_1[\mu^{\pm1}]$ is not. 

\proclaim{Proposition}\label{H1kmu}
There are isomorphisms of right $\K_1 [\mu^{\pm 1}]$\/{\rm -}\/modules as follows\/{\rm :}\/
\begin{description}
\item[\/{\rm a:}\/] 
$$ 
H_1(W;\K_1 [\mu^{\pm 1}])\cong H_1(W_r;\K_1 [\mu^{\pm 1}])\cong\frac{\K_1
[\mu^{\pm 1}]}{(\mu-1)\K_1 [\mu^{\pm 1}]}
$$
and $j_{\ast}(\eta)$ is sent to the generator $1$.
\item[{\rm b:}] 
$$
\AA_1=H_1(M;\K_1 [\mu^{\pm 1}])\cong H_1(M_r;\K_1 [\mu^{\pm 1}])\cong\frac{\K_1
[\mu^{\pm 1}]}{(\mu- 1)(\mu - k)\K_1 [\mu^{\pm 1}]}
$$ 
for a certain $k=z^\mu z^{-1}$
such that $z \in\Q S$ and{\rm ,} in the expression for $z$ the coefficient of the
additive identity of $S$ is nonzero.
\end{description} 
Moreover{\rm ,} under these identifications{\rm ,} the inclusion induced map
$j_{\ast}$ sends $1$ to~$1$.
\endproclaim 

\demo{{P}roof  of Proposition~{\rm \ref{uniqueP},} assuming
Proposition~{\rm \ref{H1kmu}}} 
The kernel $P_1$ of $j_*$ is of
rank~1 over $\K_1$ since any module of the form $\K_1 [\mu^{\pm 1}]/g\K_1
[\mu^{\pm 1}]$  has $\K_1$-rank equal to the degree of $g$. So there exists a
nonzero {\it  generator} $p_1\in P_1$. Since $j_*(\eta)$ is not zero, $\eta$ does
not lie in $P_1=P_1^{\perp}$ and hence $B\l_1(p_1,\eta)\neq0$. This establishes
one claim of  Proposition~\ref{uniqueP}.

Now we will show that the submodule generated
by $\mu-1$ is the unique proper submodule $P$ of $\AA_1$. Such a submodule
$P$ would have rank~1 over $\K_1$ and thus would be isomorphic to
$\K_1[\mu^{\pm 1}]/(\mu-b)\K_1[\mu^{\pm 1}]$ for some
$b\in\K_1$. Since the degree of $(\mu -1)(\mu -k)$ is two, we may assume that
$i:P\ra\AA_1$ sends 1 to a degree~1 polynomial. Although $i(1)$ need not be
monic, there is some $p\in P$ such that $i(p)$ is monic. But $P$ is cyclic,
generated by any nonzero element, so that
$P\cong\K_1[\mu^{\pm 1}]/(\mu-b)\K_1[\mu^{\pm 1}]$ (for a different $b$) where
$p\mapsto 1$. Thus we may assume $i(1)$ is monic, say $\mu-d$ for some
$d\in\K_1$. This necessitates
$(\mu-d)(\mu-b)=(\mu-1)(\mu-k)$ for some $d$,
$b\in \K_1$. 
The uniqueness of $P$ is therefore implied by the following lemma which is a
purely algebraic statement about the skew polynomial ring $\K_1 [\mu^{\pm 1}]$.
This lemma completes the proof of Proposition~\ref{uniqueP} and that $K$ is not
($2.5$)-solvable, modulo the proof of Proposition~\ref{H1kmu}. \hfill\qed
\enddemo

\proclaim{Lemma} \label{lem:algebra}
If $k\in\K_1$ satisfies the algebraic properties from Proposition~{\rm \ref{H1kmu}.b}
and $d,b\in\K_1$ are arbitrary then the equation in $\K_1 [\mu^{\pm 1}]$
$$
(\mu-d)(\mu-b)=(\mu-1)(\mu-k)
$$
implies that $d=1$.
\endproclaim

\demo{Proof} 
Equating coefficients, using $d\mu=\mu d^\mu$, eliminating the
variable $b$, using $k=z^{\mu}z^{-1}$ and setting $\gamma=d-1$, we are led to:
\begin{equation}\label{gammak}
\gamma z^{\mu} = (\gamma  + 1)\gamma^\mu z\qquad\hbox{ for some }\ \gamma \neq 0\
\hbox{  in }\
\K_1.
\end{equation} The solution $\gamma =0$ corresponds to the known solution
$d=1$. We will show there are no other solutions. Recall the polynomial
$p(t)=t^{-1}-3+t$ and the abelian group
$$ 
S= \Q[t^{\pm 1}]/p(t)^m
$$  
introduced earlier, where $\K_1$ is the quotient field of the group ring $\Q S$.
Suppose there is a nonzero solution $\gamma =p/q$ to Equation~\ref{gammak} where
$p$,
$q\in\Q S$,
$pq\neq0$ and $p$ and $q$ are relatively prime. We may assume that, for $p$, the
coefficient of $e$, the identity element in the group $S$, is nontrivial, by
absorbing a unit into $q$. Note that $S$ is {\it locally\/} free abelian since it
is torsion-free. Thus $p$, $q$, $p^{\mu}$, $q^{\mu}$, $z$ and $z^{\mu}$ lie in a
subring isomorphic to $\Q[\Z^n]$ for some $n$. In particular this ring is a
unique factorization domain and has only trivial units of the form $rs$ where
$r\in \Q$ and $s\in \Z^n$. Equation~\ref{gammak} then becomes
\begin{equation}\label{pqmu} pq^{\mu}z^{\mu} = (p + q)p^{\mu} z,
\end{equation}
 an equation in $ \Q[ \Z^n]$. 

\vglue8pt {\it Case} I.  $p$ and $z$ are relatively prime. 
\vglue8pt

Then any factor of $p$ (on the left-hand side of Equation~\eqref{pqmu}) must
divide $p^{\mu}$ (on the right-hand side). Thus $p=rsp^\mu$ for some unit $rs$ ($r\in
\Q$, $s\in S$). Suppose
$p=\sum r_is_i$ for  nonzero rationals $r_i$, $s_i\in S$ and
$i\in\CC$, a finite index set. Then 
$$
\sum r_is_i=\sum(rr_i)ss^{\mu}_i
$$  so that, for each
$i\in\CC$, $ss^\mu_i=s_{f(i)}$ for some
$f(i)\in\CC$. The permutation $f:\CC\to\CC$ is of finite order since $\CC$ is
finite, so there exists a positive integer $\l$ such that 
$$ 
ss^\mu s^{\mu_2}\dots s^{\mu^{\l-1}}s^{\mu^{\l}}_i=s_i
$$  
for each $i$. Note that this is a statement entirely in $S$ (not $ \Q S$):
the group operation here is from the abelian group structure on
$S$ and the action of $\mu$ comes from the group automorphism $\mu$. Recall that
$S$ is actually the {\it  additive} group of the ring
$\Q[t^{\pm 1}]/p(t)^m$, and $\mu$ acts by multiplication by
$t$. Switching to additive notation and setting
$$ 
-s'=s+s^\mu+\dots+s^{\mu^{\l-1}}
$$  
we have $(t^\l-1)s_i=s'$ for each $i\in\CC$. If
$\CC$ contains two distinct elements $s_0$ and $s_1$, say, then
$s_0-s_1$ is annihilated by $t^\l-1$. This is impossible since $t^\l-1$ and
$t^{-1}-3+t$ are relatively prime. Therefore $\CC$ contains only one element and
$p=r_0s_0$ for some $r_0\in \Q$ and $s_0\in S$. Hence $p$ is a unit and can be
assumed to be 1 by absorbing the unit into
$q$. Now Equation~\eqref{pqmu} reduces to:
\begin{equation}\label{qmu} q^\mu z^{\mu} = (1 + q)z.
\end{equation}
 Let $w=zq$. Then Equation~\eqref{qmu} becomes
$w^\mu=z+w$, in $\Q S$. Let $r_0$ be the coefficient of the additive identity
$e\in S$ in the expression for $w$ and similarly let $c_0$ be the coefficient of
$e$ for $z$. Note that $\mu$ is a {\it  group} automorphism of $S$ and as such
preserves the identity. By equating coefficients of $e$ one sees that
$a_0=c_0+a_0$, implying $c_0=0$, an obvious contradiction to
Proposition~\ref{H1kmu}.b. Thus Case~I is not possible.  

\vglue8pt {\it Case} II. $p$ and $z$ have greatest common factor $f$ in
$\Q[\Z^n]$. \vglue8pt

Suppose $p=f\tilde p$, and $z=f\tilde z$. Then
after dividing out $ff^{\mu}$ from Equation~\eqref{pqmu}, repeat the argument of
Case~I to conclude $\tilde p$ is a unit which may be assumed to be $1$. Setting
$w=\tilde z q$, we reach the same equation $w^\mu=z+w$ and arrive at the same
contradiction.
\enddemo  

\demo{Proof  of Proposition~{\rm \ref{H1kmu}}}  Let $Y$ be the cobordism between
$M_r$ and $M$. It will be convenient to refer to the following commutative
diagram. Here $\psi^r_1$ is defined using $\psi_1$ to make the diagram commute
and $\phi_1^r$ is induced \pagebreak by $\psi_1^r$. 

\centerline{\BoxedEPSF{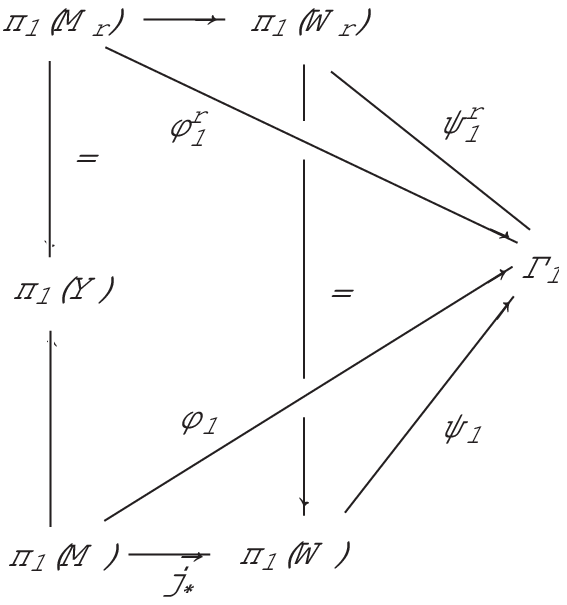 scaled 1000}}
\vglue12pt \noindent  Given any such compatible $\G_i$-coefficient systems (we need the case
$i=1$ pictured above and also the case $i=0$ which has a corresponding diagram), we
claim that one can compare $H_1(M_r;\Z\G_i)$ and $H_1(M;\Z\G_i)$ by considering
the linking matrix with $\Z\G_i$ coefficients. To see this, consider the
commutative diagram below with $\Z\G_i$ coefficients. 
$$
\BoxedEPSF{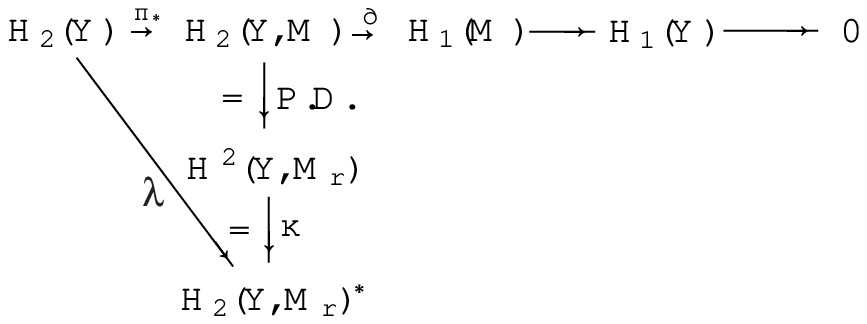 scaled 1000}
$$ 
Since $Y\simeq M_r\vee S^2\vee S^2$, the $\Z\G_i$-modules $H_2(Y,M_r)$ and
$H^2(Y,M_r)$ are free of rank two, and 
$H_1(M_r)\cong H_1(Y)$. Just as it is with untwisted coefficients, the map
$\lambda$ is given by the linking matrix of the attaching circles of the two
$2$-handles. This linking matrix has been calculated earlier to be 
$$
\left(\begin{array}{cc} t+t^{-1}-2  &1\\ 1  &t+t^{-1}-2\end{array}\right)
$$ where $t=\phi^r_i(\eta)$ and $\phi^r_i:\pi_1(M_r)\ra\G_i$. But since
$\eta\in\pi_1(M_r)^{(2)}$, and $\G_i$ is $i$-solvable, for the cases $i=0,1$ 
this is the standard hyperbolic matrix. Alternatively, note that each
$\gamma_i$ admits a Seifert surface $S_i$ which is a punctured torus lying inside
the tubular neighborhood of
$\eta$ and which lifts to the $\G_1$-cover (since
$\phi^r_1(\eta)=1$), and use these to compute the linking matrix. Since this
matrix is invertible, the above sequence implies that
 $M$ and $M_r$ have isomorphic integral Alexander modules and 
$\AA_1\cong H_1(M_r;\K_1[\mu^{\pm 1}])$.
Since
$(\G_1)^{(2)}=\{e\}$, the maps $\phi_1$ and $\phi^r_1$ factor through
$\pi_1(M)/\pi_1(M)^{(2)}$ and
$\pi_1(M_r)/\pi_1(M_r)^{(2)}$ respectively. Therefore the images of $\phi_1$ and
$\phi^r_1$ are completely dictated by the images of the induced maps 
$$ (\phi_1)_{\ast}:H_1(M;\Z[t^{\pm 1}])\ra S \mbox{  and  }
(\phi^r_1)_{\ast}:H_1(M_r;\Z[t^{\pm 1}])\ra S
$$ on the {\it integral\/} Alexander modules. But the integral Alexander modules
of $M$ and $M_r$ are isomorphic and thus images of the two maps above are
identical. We have noted previously that the kernel of
$\phi_1$ on the {\it rational\/} $\AA_0$ is $P_0$ so that the image of $\pi_1(M_r)$ and
$\pi_1(M)$ in $\G_1$ is
$$
\Z[t^{\pm 1}]/p(t)\rtimes\Z\cong(\Z \times \Z)
\rtimes\Z.
$$ 
But the latter is precisely $\pi_1(W_r)/\pi_1(W_r)^{(2)}$, where $W_r$ is the
ribbon disk complement, so that
$\psi_1$ and
$\psi_1^r$ induce a monomorphism modulo the second derived subgroup. Hence
$$ H_1(W_r;\K_1[\mu^{\pm 1}])\cong H_1(W;\K_1[\mu^{\pm 1}])\cong 
H_1(W_r^{(2)};\Z)\otimes\K_1[\mu^{\pm 1}]
$$ where $W_r^{(2)}$ is the universal abelian cover of the infinite cyclic cover
and the tensor product is over $\Z[(\Z
\times \Z)\rtimes\Z]$ (see Remark~\ref{postlore}.1). The infinite cyclic cover of
the {\it  fibered} ribbon disk complement $W_r$ is $H\times \R$ and thus is
homotopy equivalent to a wedge of two circles corresponding to $x$ and $y$. Hence
$W_r^{(2)}$ is homotopy equivalent to the usual planar grid
$$
\{(x_1,x_2)\in\R^2\mid x_1 \mbox{ or } x_2 
\mbox{ is an integer}\}
$$ and so $H_1(W_r^{(2)};\Z)$ is free on one element (we choose
$[x^{-1},y]$, the image of $C$) as a $\Z[\Z \times \Z]$-module. Therefore it is
certainly a cyclic module over
$\Z[(\Z \times \Z)\rtimes\Z]$ and the action of $\mu$ on
$[x^{-1},y]$ is given by the monodromy $\tilde f_r:H\ra H$. But we have
previously observed that $(\tilde f_r)_{\ast}([x^{-1},y])=[x^{-1},y]$.
Consequently if we denote the generator by $\overline C$ then $\overline
C(\mu-1)=0$ so that, as a {\it right\/}
$\K_1[\mu^{\pm 1}]$-module, $H_1(W_r;\K_1[\mu^{\pm 1}])$ is as claimed in
Proposition~\ref{H1kmu}.a. It is not difficult to check (using a presentation as in Figure~6.7) that the loop $\eta$ maps
to $[x^{-1},y]$ under the inclusion $j_*$ (indeed that was how $\eta$ was chosen) and so $j_*(\eta)=j_*(C)$. Thus $j_*(\eta)$ generates.

We can apply Lemma~\ref{Bl iso} to the case at hand where 
$$ 
H_1(M_r;\K_1[\mu^{\pm 1}])/P_1\cong H_1(W_r;\K_1[\mu^{\pm 1}])
$$ 
to conclude that 
$$ 
P_1\cong H_1(W_r;\K_1[\mu^{\pm 1}])^{\#}\cong
\frac{\K_1[\mu^{\pm 1}]}{(\mu^{-1}-1)\K_1[\mu^{\pm 1}]}\cong\frac{\K_1[\mu^{\pm
1}]}{(\mu-1)\K_1[\mu^{\pm 1}]}
$$   
and in particular
$rk_{\K_1}(P_1)=1$. Now we can make use of a theorem that any finitely-generated
$\K_1[\mu^{\pm 1}]$-module is cyclic ~\cite[Prop.~2.2.8 and
Th.~1.5.5]{Co2}. (We will also shortly derive the fact that $\AA_1$ is
cyclic by explicit computation.) Since $rk_{\K_1}\AA_1$ = $2 rk_{\K_1} P_1 = 2$,
we have that
$$
\AA_1\cong H_1(M_r;\K_1[\mu^{\pm 1}])\cong\K_1[\mu^{\pm 1}]/g\K_1[\mu^{\pm 1}]
$$ 

\centerline{\BoxedEPSF{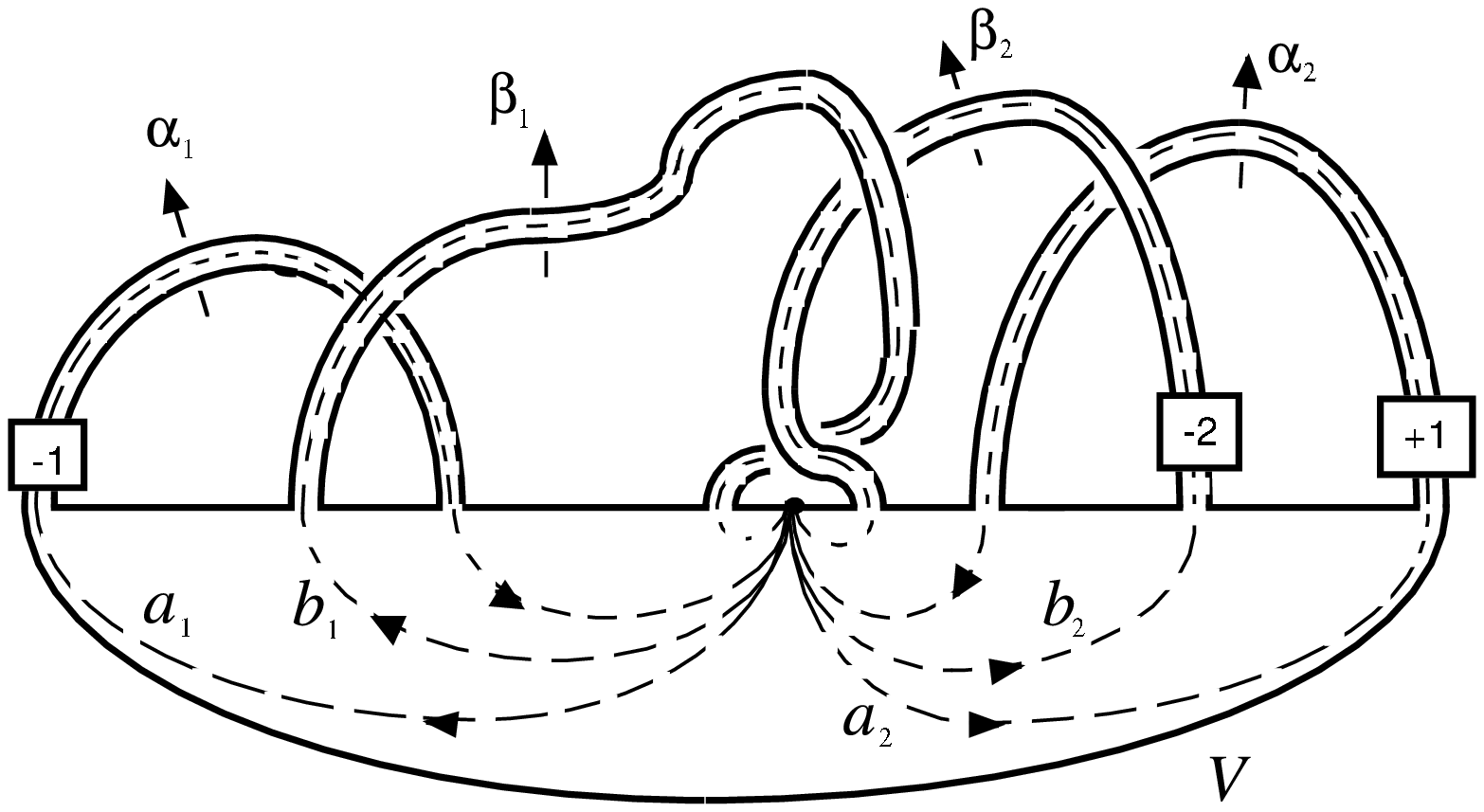 scaled 500}}
\vglue6pt
\centerline{Figure 6.7}
\vglue12pt 

\noindent  where $g$ is a monic degree~2 polynomial in $\mu$. Since this module admits
an epimorphism to
$$
\K_1[\mu^{\pm 1}]/(\mu-1)\K_1[\mu^{\pm 1}] \ \ (\cong\K_1)
$$ it admits such an epimorphism sending $1$ to $1$. The image of
$g$ lies in\break $(\mu-1)\K_1[\mu^{\pm 1}]$  so 
$g=(\mu-1)(\mu-k)$ for some $k\in \K_1$. The kernel of this epimorphism, $P_1$, 
is clearly generated by $\mu -1$ and hence is itself cyclic of order $\mu -k$.
But above we saw that $P_1$ is cyclic of order $\mu -1$.  Hence
$$
\frac{\K_1[\mu^{\pm 1}]}{(\mu-k)\K_1[\mu^{\pm 1}]}\cong\frac{\K_1[\mu^{\pm
1}]}{(\mu-1)\K_1[\mu^{\pm 1}]}.
$$ This does not imply that $k=1$ since $\K_1[\mu^{\pm 1}]$ is noncommutative.
However it can be shown that this is equivalent to the fact that  there exists some nonzero
$z\in\K_1$ such that
$k=z^\mu z^{-1}$ ~\cite[p.~112 Lemma~3.4.2]{Co2}. 

In summary, we have demonstrated Proposition~\ref{H1kmu}.a and b except for
showing that $z$ has the required properties.
 For this purpose we are forced to move to a second,
more explicit, calculation of $\AA_1$ --- including a determination of the
constant
$k$ via a computation of the
$\pi_1$ monodromy of $K_r$. This calculation constitutes the remainder of this
section.

First we show that it is quite easy to determine an explicit presentation of
$\AA_1$ as a $\K_1$-module. We have already identified $\AA_1$ with
$H_1(M_r;\K_1[\mu^{\pm 1}])$. Recall that we argued in the proof of
Proposition~\ref{H1kmu}.a that the image of $\phi_1$ in $\G_1$ is $(\Z \times
\Z)\rtimes \Z$. Let
$\G^\l_1$ denote this subgroup. The infinite cyclic cover of
$S^3-K_r$ is homotopy equivalent to a wedge of four circles whose fundamental group
is free on $\{A,B,a_1,b_1\}$ whereas
$\pi_1$ of the infinite cyclic cover of $W_r$ is free on
$\{j_{\ast}(a_1),j_{\ast}(b_1)\}$. It is then immediate that, as a  $ \Z[\Z \times
\Z]$-module, $H_1(S^3-K_r; \Z\G^\l_1)$ is isomorphic to the direct sum of
$H_1(W_r;
\Z\G^\l_1)$ and a free $ \Z[\Z \times \Z]$-module of rank~2 with basis
$\{A,B\}$. The same holds for coefficients in $\K_1[\mu^{\pm 1}]$ {\it as
$\K_1$-modules}. The inverse of the longitude of $K_r$, $\l^{-1}$, equals
$[b^{-1}_1,a_1][a_2,b^{-1}_2]$ (we use the convention
$[a,b]=aba^{-1}b^{-1})$. Since $A=a_1a^{-1}_2$ and
$B=b_1b^{-1}_2$, the element $\l^{-1}$ represents
$$ A(x-1)+B(xy^{-1}-x) \in H_1(S^3-K_r;\Z\G^\l_1)
$$ where
$\{x,y\}$ are used as the basis of the $\Z \times \Z$ action. This can be
demonstrated by rewriting $\l^{-1}$ as
$$ [b^{-1}_1,a_1]A^{-1}[a_1,b^{-1}_1]\left(b^{-1}_1a_1Ba^{-1}_1b_1\right)
\left(b^{-1}_1Ab_1\right)\left(b^{-1}_1B^{-1}b_1\right),
$$ 
where we recall  that $j_{\ast}(b_1)=x$ and $j_{\ast}(a_1)=y$. Recall also our
convention that if $\widetilde X\xrightarrow{P} X$  is a regular cover, then a
right
$\Z[\pi_1(X)/p_{\ast}(\pi_1(\widetilde X))]$-module structure on
$H_1(\widetilde X)$ (viewed as the abelianization of
$p_{\ast}(\pi_1(\widetilde X))$) is given by
$[\gamma]\mu_{\ast}=[\mu^{-1}\gamma\mu]$ where $[\gamma]\in H_1(\widetilde X)$ and
$\mu=\mu p_{\ast}(\pi_1(\widetilde X))$ is a right coset. Therefore we arrive at
the {\it longitudinal relation\/}:
\begin{equation}\label{longitudinal relation} B = A(x^{-1} - 1)_{\ast}(y^{-1} -
1)^{-1}_{\ast}.
\end{equation} Since $y^{-1}-1$ is a unit in $\K_1$, $\AA_1$, as a $\K_1$-module,
is free on $\{A,C\}$. Since
$H_1(W_r;\K_1[\mu^{\pm 1}])\cong\K_1$, generated by
$j_{\ast}(C)$, and since $A$ bounds a disk in $W_r$, this shows that the kernel of
$j_{\ast}$ in Proposition~\ref{H1kmu}.b is indeed the $\K_1$ subspace spanned by
$A$. Note that this implies that the subspace is invariant under the action of
$\mu$, a fact we shall presently confirm by direct calculation.

To calculate the structure of $\AA_1$ as a $\K_1[\mu^{\pm 1}]$-module, we must
derive the action of $\mu$ on $A$ and $C$. Since $\AA_1\cong
H_1(M_r;\K_1[\mu^{\pm 1}])$, it is certainly sufficient to know the monodromy (on
$\pi_1$) of $K_r$ (indeed it would suffice to know it for $K'$). Let
$V$ be the Seifert surface for $K_r$ and $F=\pi_1(V,{\ast})$ the free group on
$\{a_1,b_1,a_2,b_2\}$ as in Figure~6.7. We will use the dual basis
$\{\a_1,\b_1,\a_2,\b_2\}$ of $\pi_1(S^3-V)$ to help calculate. The basepoint is
on the boundary of a tubular neighborhood of $K_r$, on the negative side of $V$
as shown in Figure~6.7.

Refer to Figure~6.8, which is a schematic picture of the infinite
cyclic cover of $S^3-K_r$. If $x$ is a based loop then
$x\mu_{\ast}=\mu^{-1}x\mu$ is obtained by traveling from the basepoint halfway
around the meridian (in the negative direction) until reaching the positive side
of $V$, traversing $x^+$, then returning to the basepoint along the same path.
This must then be written in terms of the chosen basis, which will be
$\{a_1,b_1,a_2,b_2\}$ or
$\{a^-_1,b^-_1,a^-_2,b^-_2\}$ which are identified. The initial calculations
follow.
\figin{Fig68}{550}
\vglue6pt
\centerline{Figure 6.8}
\eject

These were accomplished using a 25-foot extension cord on a living room
floor to simulate $K_r$. Here $\delta =[\b_1,\b^{-1}_2]$. These are essentially
the positive push-offs of $a_i$ and $b_i$, but based at the basepoint on the
negative side of $V$ as described above.
\begin{eqnarray} 
a_1\mu_{\ast}  &=& \delta\a^{-1}_1\delta^{-1},
\label{initial calc}\\
 a_2\mu_{\ast}  &=& \delta\a^{-1}_2\delta^{-1} ,\nonumber\\
b_1\mu_{\ast}  &=& \delta\a^{-1}_1\b_2\delta^{-1}, \nonumber\\ b_2\mu_{\ast}  &=&
\delta\a^{-1}_2\b_2\b^{-1}_1\b_2\delta^{-1}. \nonumber
\end{eqnarray}
 Now we must translate $\{\a_1,\b_1,\a_2,\b_2\}$ into
$\{a_1,b_1,a_2,b_2\}$ using the negative push-offs:
\begin{eqnarray} a_1  &= &a^-_1 = \a^{-1}_1\b^{-1}_1,   \label{translate}\\ a_2  &=&
a^-_2 =
\a^{-1}_2\b^{-1}_2 ,\nonumber\\ b_1  &= &b^-_1 = \b_1\b_2\b^{-1}_1, \nonumber\\ b_2  &=&
b^-_2 =
\b^2_2\b^{-1}_1. \nonumber
\end{eqnarray}
 These enable us to solve: $\b_2=b_2b_1b^{-1}_2$,
$\b_1=b^2_1b^{-1}_2$, $\a_1=b_2b^{-2}_1a^{-1}_1$,
$\a_2=b_2b^{-1}_1b^{-1}_2a^{-1}_2$. The ``monodromy'' equations~\ref{initial
calc} then become:
\begin{eqnarray} a_1\mu_{\ast}  &=& \delta a_1b^2_1b^{-1}_2\delta^{-1}
\quad\hbox{ where }\ \delta=[b^{-1}_1,b_2],\label{monodromy}\\ a_2\mu_{\ast}  &=&
\delta a_2b_2b_1b^{-1}_2\delta^{-1}, \nonumber\\ b_1\mu_{\ast}  &=&\delta
a_1b^3_1b^{-1}_2\delta^{-1}, \nonumber\\ b_2\mu_{\ast}  &= &\delta a_2b_2b_1. \nonumber
\end{eqnarray}
Using these one may calculate the $\K_1[\mu^{\pm 1}]$-module relations:
\begin{eqnarray} A\mu  &= &A + B(y^{-1}) ,\label{kmu1}\\ B\mu  &= &A + B(y^{-1} +
x^{-1}y^{-1})
,\label{kmu2}\\ C\mu  &=& C + B(x-1).\label{kmu3}
\end{eqnarray}

Let us clarify the meaning of these relations. It is easy to think of taking the
infinite cyclic cover of $M_r$ and then the $\Z \times \Z$ cover of that (this is
the boundary of the universal abelian cover of the infinite cyclic cover
$W_\infty$ of $W_r$.) The group of deck translations of this regular cover is
$\G^\l_1=(\Z\times\Z)\rtimes \Z$ where $\mu$ generates $\G_0= \Z$ and $x$ and
$y$ generate $\Z \times \Z$ (corresponding to the chosen generators of
$H_1(W_\infty)$). The above module relations do constitute a presentation of
$H_1(M_r)$ with coefficients in the ring $ \Z\G^\l_1$ with the subset $
\Z[x^{\pm1},y^{\pm1}]-\{0\}$ inverted. However, this ring is a subring of
$\K_1[\mu^{\pm 1}]$ and it is this larger ring which is our intended coefficient
ring. Recall that
$\G_1= S\rtimes\G_0$ where $S=\Q[t^{\pm 1}]/p(t)^m$ and where 
$$ 
H_1(W_r; \Q[t^{\pm 1}])\cong \Q[t^{\pm 1}]/p(t)
$$ 
embeds in $S$ in the standard fashion, and the map $H_1(M_r;
\Q[t^{\pm 1}])\ra S$ factors through the inclusion to
$H_1(W_r; \Q[t^{\pm 1}])$. The latter fact justifies our replacing things like
$a^{-1}_1\gamma a_1$ with $\gamma y$ in the above calculations, since
$j_{\ast}(a_1)=y$ and
$j_{\ast}(b_1)=x$. Therefore the reader sees that $\G^\l_1$ is naturally the
subgroup
$H_1(W_r; \Z[t^{\pm 1}])\rtimes \Z$ of
$\G_1$, and the elements $x$, $y$ in the above relations are to be viewed in this
way as elements of $\K_1[\mu^{\pm 1}]$, which, as you recall, is $\Z\G_1$ with $
\Z S-\{0\}$ inverted. Finally, if $s\in S$ then $s^\mu$ will denote the element
$\mu^{-1}s\mu$ of $\G_1$ or more precisely the image of $s$ under the action of
$\G_0$ on $S$. What is this action? Note that under our identifications it agrees
with the ordinary action of $\mu$ on the Alexander module
$$ H_1(W_r; \Z[t^{\pm 1}])\cong\Z \times \Z \cong\left<x,y\right>.
$$ The element $x^\mu$ may be calculated from the formula above for
$\mu^{-1}b_1\mu=b_1\mu_{\ast}$ by setting $b_1=b_2=x$ and
$a_1=a_2=y$ and abelianizing. Thus we get:
\begin{equation}\label{xmu} x^\mu = yx^2,\ \ \ y^\mu = yx.
\end{equation}

The relation~\eqref{kmu2} is not needed since $B$ was eliminated using the
longitudinal relation; however it can be used as a ``consistency check''. The
relation~\eqref{kmu3} could have been derived more easily. Recall that we argued
that
$C$ was fixed under the monodromy of $K'$ (connected sum of figure eight with
itself) and so the monodromy of $K_r$ sent $C$ to the image of $C$ under a Dehn
twist along
$B$. Since $C$ intersects $B$ in two points, relation~\eqref{kmu3} can be easily
deduced. From relation~\eqref{kmu3} and the longitudinal
relation~\eqref{longitudinal relation} we obtain
\begin{equation}\label{Cmu} C(\mu - 1)(x - 1)^{-1}w = A
\end{equation} where $w=(y^{-1}-1)(x^{-1}-1)^{-1}$, showing that $\AA_1$ is
cyclic, generated by $C$ and thus completing the proof of
Proposition~\ref{H1kmu}.b. Using relation~\eqref{kmu1} and $B=Aw^{-1}$ we get
$A(\mu-1-w^{-1}y^{-1})=0$. Combining this with~\eqref{Cmu} yields that $C$ is
annihilated by
$$ (\mu-1)(x-1)^{-1}w(\mu-1-w^{-1}y^{-1}).
$$ Let $s=(x-1)^{-1}w$ and let $r=1+w^{-1}y^{-1}$. Note that 
$$ s(\mu-r)=\mu s^\mu-sr= (\mu-sr(s^\mu)^{-1})s^\mu.
$$  Hence $C$ is annihilated by $(\mu-1)(\mu-k)$ where
$k=sr(s^\mu)^{-1}$. This simplifies to
\begin{equation}\label{simpler} k = \frac{(x-1)^\mu}{x-1}\
\frac{(x^{-1}-1)^\mu}{(x^{-1}-1)}.
\end{equation} This provides the specific $k$ of Proposition~\ref{H1kmu}.b. Note
that
$k=z^\mu z^{-1}$ where $z=(x-1)(x^{-1}-1)$ and the coefficient of the identity is
2. This concludes the proof of Proposition~\ref{H1kmu}.b. The only extra
information we have obtained is the specific value of $z$, which we needed in the
proof of Proposition~\ref{uniqueP}. Probably there is a way to deduce this
property of $z$ without explicit calculation, in which case the 25-foot extension
cord is not needed!
\enddemo

\proclaim{Theorem}\label{Q slice} \hskip-8pt The zero surgery on the knot $K$ of Figure~{\rm 6.1}   is
$(2)$\/{\rm -}\/solvable but not rationally $(2.5)$\/{\rm -}\/solvable with multiplicity $1$ \/{\rm (}\/\/see the
discussion above Theorem~\/{\rm \ref{invariant thm}).}\/ In particular\/{\rm ,}\/ the knot $K$ of
Figure~\/{\rm 6.1}\/   is not slice in any rational homology ball wherein the meridian of
$K$ generates the free part of $H_1$ of the exterior of the slice disk.
\endproclaim

\demo{Proof}  Repeat the proof that $K$ is not $(2.5)$-solvable. The only place
where we used solvability was in claiming that $\phi_0$ was the usual epimorphism. The
point is that when the multiplicity is $1$,
$\AA_0$ is the  usual Alexander module of $K$ and our calculations are correct.
If it is not $1$ then $\AA_0$ is larger  and new calculations would need to be
made. This has not been attempted.
\enddemo

\vglue-6pt
\section{$(n)$-surfaces, gropes and Whitney towers}
\label{sec:surfaces} 
\vglue-6pt 

Consider a regular covering $X_N\to X$ of smooth connected
oriented $4$-manifolds, where $\pi_1(X_N)
\cong N$ is a normal subgroup of $\pi_1(X)$. To be precise, all the spaces are
equipped with a base point (which   will be suppressed from the notation).
\numbereddemo{Definition} \label{def:$N$-surface} Let $F$ be a closed oriented surface. An {\it  $N$-surface} in
$X$   is a generic immersion $f:F\imra X$ such that
$f_{\ast}(\pi_1(F)) \leq N$. In addition, the surface is equipped with a {\it  whisker}, i.e.\  an arc in $X$ from the base
point of $X$ to the image of the base point of $F$.
\enddemo 

By covering space theory, an  $N$-surface lifts uniquely to a generic
immersion $f_N:F\imra X_N$ leading to the induced homology class
$[f]:=(f_N)_{\ast}[F]
\in H_2(X_N)$. Clearly any class in $H_2(X_N)$ is represented in this way. The  
group of deck transformations $\pi_1(X)/N$ acts on $X_N$ and thus on
$H_2(X_N)$. On lifts of 
$N$-surfaces and their homology classes, this action is given by pre-composing
the whisker with a loop in $\pi_1(X)$. Moreover, addition in
$H_2(X_N)$ corresponds to a connected sum of $N$-surfaces along their whiskers.

\proclaim{Lemma} \label{lem:mu} An  $N$\/{\rm -}\/surface $f$ has a Wall self\/{\rm -}\/intersection
invariant 
$$
\mu (f) \in \Z [\pi_1(X)/N]/\{a-\bar a\}.
$$ Here there is the usual involution $a \mapsto \bar a$ on elements $a$ in the
group ring $\Z[\pi_1(X)/N]$ which is induced by $\bar g  :=g^{-1}$ for group
elements.

If $\mu (f) = 0$ then the induced homology class $[f] \in   H_2(X_N)$ is
represented by an {\rm  embedded} $N$\/{\rm -}\/surface whose image can be chosen to be
arbitrarily close to the image of $f$.
\endproclaim

\demo{Proof} To define $\mu(f)$ recall that by definition $f:F\imra X$ only has
transverse double points $p$. Choose two arcs in
$f(F)$ leaving $p$ on different sheets and ending at the base point of $F$,
missing all other double points on their way. Together with the whisker for $f$
this gives an element  
$g_p \in \pi_1 (X)$, the so-called {\it  double point loop} at
$p$. Since there is no preferred order of the two sheets, we have to identify
$g_p$ with $g_p^{-1}=\bar g_p $.

Moreover, $f_{\ast}(\pi_1(F)) \leq N$ implies that the choice of the two arcs
only changes $g_p$ by elements in $N$. This implies that
$$
\mu(f) :=\sum_p {\e_p \cdot g_p}
$$ is well-defined in the quotient of the group ring above. Here
$\e_p\in \{\pm 1\}$ is the usual sign of the double point
$p$, coming from the orientations of $F$ and $X$.

Now assume that $\mu(f)=0$. Then the double points of $f$ can be paired up with
signs and double point loops. Consider such a pair $p$, $p'$ with
$\e_{p'}=-\e_p$ and
$g_{p'}=g_p^{\pm 1} \in\pi_1 (X)/N$. Consider a loop $\alpha $ in $f(F)$ which
leaves
$p$ on one sheet, changes sheets at $p'$ and returns on the other sheet to $p$.
There are two ways of making $\alpha $ into a loop at the base point of $X$, and
by assumption (at least) one of them lies in the subgroup $N\leq \pi_1(X)$.
Consider the corresponding subarc $\alpha_0$ of $\alpha$ leading from
$p$ to $p'$. Let $T$ be the normal bundle of $f$ restricted to $\alpha_0$. By
construction, $f_{\ast}(\pi_1(F\cup T))
\leq N$ and we may use $T$ to do surgery on~$f$: Remove two small disks around
$p$ and
$p'$ (on the correct sheets) and replace them by the annulus which is the normal
circle bundle corresponding to $T$. This procedure of {\it  adding a tube} along
$\alpha_0$ removes the pair of double points
$p$,~$p'$, stays in the same homology class in $H_2(X_N)$, and can be done
arbitrarily close to the image of $f$. A finite number of such tube additions
produces the desired   embedded $N$-surface.
\enddemo
 \vglue9pt

Let $w_2(f)\in \Z/2 $ be the second Stiefel-Whitney number of the normal bundle
of an
$N$-surface $f$. Note the identities
$$ w_2(f)= \langle f^{\ast}w_2(X),[F] \rangle = \langle f^{\ast}_N w_2(X_N),[F]  
\rangle = \langle w_2(X_N),[f] \rangle .
$$ In particular, $w_2(f)$ only depends on the induced homology class  $[f] \in
H_2(X_N)$ and vanishes if $X_N$ is a spin manifold.

The self-intersection invariant $\mu(f)$ is clearly unchanged under isotopies,
finger moves and Whitney moves, i.e.,   under a regular homotopy of immersions.
As usual, it is not invariant under an arbitrary homotopy: A local kink changes
$\mu(f)$ by adding $\pm$  the trivial group element. This move also changes the
Euler number $e(f)$ of the normal bundle of $f$ by $\pm 2$. Therefore, if
$w_2(f)=0\in \Z/2 $ then\pagebreak there is a well-defined number of kinks one has to have to
make   the normal bundle of
$f$ trivial. More precisely, we define a {\it  homotopy invariant}
$$
\mu(f):= \sum_p \epsilon_p\cdot g_p - \frac{e(f)}{2}\cdot 1
$$ and we will use $\mu(f)$ in this sense in the rest of the paper. 

\begin{quote} {\it  We will also assume that all surfaces $f$ with
$w_2(f)=0 $ are represented by framed immersions.}
\end{quote}

The same proof as for Lemma~\ref{lem:mu} above shows if that
$\mu(f)=0$ (in the modified definition) then $[f] \in H_2(X_N)$ is represented by
a {\it  framed,   embedded} $N$-surface, i.e.\ an $N$-surface with trivial normal
bundle. One just has to observe that the double point from a local kink can   be
removed by the following procedure (which stays within the class of
$N$-surfaces and does not change   the normal Euler number): A neighborhood
around the double point $p$ can be identified with
$D^2 \times D^2$ with the two sheets of $f$ being $D^2 \times 0$ and $0\times
D^2$. On the boundary of this 4-ball we see a Hopf-link. It bounds a twice-twisted band in
$S^3$. If we cut out the local sheets around $p$ and replace them with this
annulus we get the desired result.
\vglue4pt 
The following result is an exercise in Morse theory and will not be proved here.
We will not need this result and only state it for the sake of completeness.

\proclaim{Lemma}\label{w2} If $w_2(f)=0$ then the homotopy invariant $\mu(f)$  
only depends on $[f] \in H_2(X_N)$. Moreover{\rm ,} $\mu (f)=0$ {\rm  if and only if}
$[f]$ is represented by a framed embedded $N$\/{\rm -}\/surface.
\endproclaim

We next turn from self-intersection to intersection numbers. Let $f,g$ be
$N$-surfaces meeting in general position. Then their (Wall) intersection number
$\lambda(f,g) \in
\Z[\pi_1(X)/N]$ is defined by the formula
$$
\lambda(f,g) =\sum_p \epsilon_p\cdot g_p
$$ where $p$ runs through the intersection points of $f$ and
$g$. Signs $\epsilon_p$ and group elements $g_p$ are defined similarly to the
self-intersection invariant.

It is well-known that $\lambda(f,g)$ only depends on $[f],[g]
\in H_2(X_N)$. In fact, under the isomorphism $H_2(X_N)
\cong H_2 (X;\Z[\pi_1(X)/N])$ the pairing $\lambda$ corresponds to the
composition of the following three obvious maps (where $\Lambda  
:=\Z[\pi_1(X)/N]$):
$$ H_2(X; \Lambda ) \ra H_2(X, \partial ; \Lambda )
\xrightarrow{\cong} H^2(X; \Lambda )    \ra
\Hom_\Lambda (H_2(X; \Lambda ), \Lambda ).
$$
\proclaim{Lemma} \label {lem:l} Let $f,g$ be $N${\rm -}surfaces with
$$
\lambda(f,g)= \sum^r_{i=1} \epsilon_i \cdot g_i.
$$ Then there exists an $N$\/{\rm -}\/surface $f'$ with $[f']=[f] \in H_2(X_N)$ which
intersects
$g$ geometrically in exactly $r$ points $p_i$ with group elements $g_{p_i} = g_i$
and signs
$\epsilon_i$ for $i=1,  
\dots,r$.
\endproclaim

The proof of this lemma proceeds exactly as the proof of   Lemma~\ref{lem:mu}:
Algebraically, canceling pairs of intersection points can be removed by adding
tubes to
$f$. We leave the details to the reader. The algebraic properties of the
intersection invariants $\lambda$, $\mu$ on
$H_2(X_N)$ can be summarized as follows.

\numbereddemo{Definition}\label{def:form} Let $R$ be a ring with involution and $M$ a left
$R$-module. A {\it  quadratic form } on $M$ consists of a $\Z$-bilinear map $
\lambda:   M \times M \ra R$ together with a map $\mu:M
\ra \overline{R}   :=R/\{r- \bar r \mid r\in R\}$ satisfying the following
properties:
\vglue4pt
 1.    $\lambda(a,b)=\overline{\lambda (b,a)}$.
\vglue4pt 2.    $\lambda(r\cdot a,b)=r\cdot \lambda(a,b)$ for all
$r \in  R$.
\vglue4pt 3.  $\mu(r\cdot a)=r \cdot\mu(a)\cdot\bar r \in
\overline{R}$  for all $r \in  R$.
\vglue4pt 4.    $\mu(a+b)= \mu(a) +\mu(b) + \lambda(a,b)  \in \overline{R} $.
\vglue4pt 5.     $\lambda(a,a)=\mu(a)+ \overline{\mu(a)}  \in R$.
\vglue4pt
One also calls $\lambda$ a {\it  hermitian form} on $M$ and
$\mu$  its {\it  quadratic refinement}. The form is {\it  nonsingular} if   the
homomorphism $M\to \Hom_R(M,R)$ induced by $\lambda$ is an isomorphism. It is
{\it  non\/{\rm -}\/degenerate} if this map is a monomorphism.
\enddemo

The algebraic intersection and self-intersection numbers
$\lambda$, $\mu$ above define a quadratic form on the module
$\Ker \{w_2:H_2(X_N)\ra \Z/2 \}$ over the ring
$\Z[\pi_1(X)/N]$. In general, this form may   be degenerate.

A simple example of a nonsingular quadratic form is the hyperbolic form on a free
$R$-module of rank $2g$: On a basis $e_1,\dots, e_g,   f_1, \dots, f_g$ one has
by definition
$$
\lambda (e_i,f_j)= \delta _{i,j} \hbox{  and }  
\lambda(e_i,e_j)=0=\lambda(f_i,f_j) \hbox{  and }   
\mu(e_i)=0=\mu(f_i).
$$ It follows from Lemma~\ref{lem:l} that if one has classes
$e_i,f_j\in \Ker \{w_2:H_2(X_N) \ra \Z/2 \}$ satisfying the above equations, 
then there are disjointly embedded
$N$-surfaces $E_i$ (with trivial normal bundle ) representing $e_i$ and
disjointly embedded  $N$-surfaces
$F_i$ representing $f_i$. Moreover, $E_i$ and $F_i$   intersect (transversely) in
exactly one point, whereas for $i
\neq   j$ the surfaces $E_i$ and $F_j$ are disjoint.

\numbereddemo{{R}emark}    \label {rem:hyperbolic} Let $\lambda,\mu$ be a nonsingular
quadratic form on a free  
$R$-module which has a Lagrangian. More precisely, there is a half-basis $e_1,
\dots, e_g$ satisfying
$$
\lambda(e_i,e_j)=0   \hbox{ and }   \mu(e_i)=0.
$$ Then this quadratic form is hyperbolic. The proof is by induction on $g$ and
proceeds as follows: Since $\lambda$ is nonsingular and $e_i$ are basis vectors,
there   exists a vector $f'_1$ such that
$\lambda(e_i,f'_1)=\delta_{i,1}$. Define
$$ f_1:=   f_1' - \mu(f_1') \cdot e_1 .
$$ Then we still have $ \lambda(e_i,f_1)= \delta_{i,1}$ but in   addition
$\mu(f_1) =0$. This implies that the sub-form $
\langle e_1,f_1 \rangle $ is hyperbolic with $e_2, \dots, e_g$ lying in the
orthogonal complement. By induction, the form $\lambda,\mu$ is hyperbolic, too.
\enddemo

We should also remark that a quadratic form with a Lagrangian
$\langle  e_1, \dots, e_g \rangle$ is nonsingular if and only if the vectors
$e_i$ have {\it  duals}. That is to say, there are vectors $d_1, \dots, d_g$ such
that  
$\lambda (e_i,d_j) = \delta_{i,j} $. The above remark shows how to   improve the
$d_i$ (by summing with linear combinations of $e_j$) to a hyperbolic basis
$f_i$.

We finish this section by explaining Whitney towers and gropes.

\numbereddemo{Definition} \label{def:Whitney tower} Let $\gamma$ be a framed circle in the
boundary
$M$ of a
$4$-manifold $W$. A {\it  Whitney tower of height}~1 is an immersed disk $ 
\Delta$ in
$W$ which bounds $\gamma$ and such that the unique framing on   the normal bundle
of
$  \Delta$ restricts to the given framing on  $\gamma$. If the double points of
$ 
\Delta$ can be paired up (with signs and   double point loops) then the choice of
Whitney circles enables one   to iterate the construction. Recall that a Whitney
circle is framed   by a vector field which is tangent along one sheet and normal
along   the other. By convention, a Whitney disk (which bounds a Whitney circle)
is   allowed to have (transverse) double points but it is always assumed   to be
{\it  framed} in the sense that the above vector field on the   Whitney circle
extends to a nonvanishing normal vector field on the   Whitney disk (see
\cite[p.17]{FQ}).

For $n\in \N$, a {\it  Whitney tower of height~$n$} on $\gamma$ is a sequence
$\CC_j=\{  \Delta_{j,k}\}_k$, $j=1,\dots n$, of   collections of Whitney disks
$ 
\Delta_{j,k}$ in general position   (where $\CC_1$ is the Whitney disk with
boundary
$\gamma$) with the following property:
\begin{itemize}
\item For $j=2,\dots,n$ the collection $\CC_j$ pairs up all
$\CC_{j-1}$-(self)-intersections and has interiors disjoint from
$\CC_1,\dots,\CC_{j-1}$.
\end{itemize} A Whitney tower of height~$(n.5)$ has an additional collection
$\CC_{n+1}$ of framed immersed Whitney disks such that
\begin{itemize}
\item $\CC_{n+1}$ pairs up all $\CC_{n}$-(self)-intersections and has interiors
disjoint from  
$\CC_1,\dots, \CC_{n-1}$ (but $\CC_{n+1}$ is   allowed to intersect the previous
collection $\CC_{n}$).
\end{itemize} Finally, we define the notion of a Whitney tower in a slightly  
different situation:  a {\it  Whitney tower of height~$0$} is a collection
$\CC_0$ of
$2$-spheres $S_i\imra W^4$. For  $n\in \N$, a {\it  Whitney tower of height~$n$ on
$\CC_0$} is a sequence
$\CC_j=\{  \Delta_{j,k}\}_k$, $j=1,\dots n$, of   collections of framed immersed
Whitney disks $  \Delta_{j,k}$ in  general position with the following property:

\begin{itemize}
\item For $j=1,\dots,n$ the collection $\CC_j$ pairs up all 
$\CC_{j-1}$-(self)-intersections and has interiors disjoint   from
$\CC_0,\dots,\CC_{j-1}$.
\end{itemize} A Whitney tower of height~$(n.5)$ has an additional collection
$\CC_{n+1}$ of framed immersed Whitney disks such that

\vglue-28pt
\phantom{why?}

\begin{itemize}
\item $\CC_{n+1}$ pairs up all $\CC_{n}$-(self)-intersections and has interiors
disjoint from  
$\CC_0,\dots, \CC_{n-1}$ (but $\CC_{n+1}$ is   allowed to intersect the previous
collection $\CC_{n}$).
\end{itemize}

\enddemo

\vglue-16pt
{\it {R}emark} 7.8.  By definition, a Whitney tower of height~$(0.5)$ on
$\CC_0$   exists if and only if the algebraic (self)-intersection numbers
$\lambda$ and
$\mu$ vanish on the $2$-spheres $S_i$.
\advance\theoremcount by 1
\vglue6pt

The following definition and lemma are taken from \cite{FT}.

\numbereddemo{Definition}\label{def:grope} A {\it  grope} is a special pair (2-complex, base
circle).  A grope has a {\it  height} $n\in \N$.  For $n=1$ a grope is precisely a
compact oriented surface $\Sigma$ with a single boundary component which is the
base circle.  A grope of height $(n+1)$ is defined inductively as follows:  Let
$\{\alpha_i, i=1, \dots, 2g\}$ be a standard symplectic basis of circles for
$\Sigma$, the {\it  bottom} stage of   the grope.  Then a  grope of height
$(n+1)$ is formed by attaching gropes of height $n$ to each $\alpha_i$ along the
base circles. Finally, a grope of height~$(n.5)$, $n\in\N$, has a   bottom
surface $\Sigma$ which on one half basis of curves bounds   gropes of
height~$(n-1)$  and on the dual half basis of curves   bounds gropes of
height~$n$.
\enddemo

Thus  a grope of height $n$ has $n$ surface stages and its fundamental group is
freely generated by the circles of the symplectic basis for all the surfaces in
the top stage. For example, if all the surfaces in the grope have genus 1 then
there are
$2^{(n-1)}$ top stage surfaces each giving 2 free generators.

\proclaim{Lemma}\label{gropebound} For a space $X${\rm ,} a loop $\gamma$ lies in
$\pi_1(X)^{(n)}$ if and only if $\gamma$ bounds a map of a grope of height $n$
{\rm (}\/i.e.\ $\gamma$ becomes the base circle of that grope\/{\rm ).}  Moreover{\rm ,} the height of a grope
$(g,\gamma)$ is the maximal
$n\in \N$ such that $\gamma\in\pi_1(g)^{(n)}$.
\endproclaim

As one can see from Figure~1.1 every grope
$(g,\gamma)$ embeds properly (i.e.\ boundary goes to boundary) into $\big(\R^3_+,
\R^2\times
\{0\}\big)$ mapping $\gamma$ to the unit circle in
$\R^2$. This determines a framing of the grope or an ``untwisted''  thickening.
Restricted to each surface stage this framing is a nonvanishing normal vector
field which on the boundary restricts to   a vector field tangent to the lower
surface stage. In particular,   the framing does not depend on the embedding into
$\R^3$.

Given a $4$-manifold $W$ with boundary $M$ and a framed circle $\gamma$ in $M$,
we say that $\gamma$ bounds a {\it  grope} in $W$ if $\gamma$ extends to an
embedding of a grope with its untwisted   framing. Knots in $S^3$ always are
equipped with the linking number zero \pagebreak framing.

\section{$H_1$-bordisms} \label {sec:$H_1$-bordism}  
 \vfil
We fix a closed oriented
$3$-manifold $M$ and consider the   following class of $4$-manifolds.
\vfil
\numbereddemo{Definition}\label{bord}  An {\it  $H_1$-bordism} is a
$4$-dimensional spin  manifold $W$ with boundary $M$ such that the inclusion map
induces an isomorphism $ H_1 (M)\cong H_1 (W)$.
\enddemo
\vfil
Note that any spin structure on $M$ extends to a spin structure on an
$H_1$-bordism $W$ because the affine spaces of spin structures are isomorphic via
the isomorphism $H^1(W,\Z/2) \cong H^1(M,\Z/2)$.
\vfil
\numbereddemo{{R}emark} \label {rem:Arf} If $M$ is the 0-surgery on a knot $K$ in $S^3$ then
an $H_1$-bordism exists if and only if the Arf invariant of $K$   vanishes. This
fact is well-known and follows from the computation   of the bordism group
$\Omega^{\rm spin}_3 (S^1)
\cong\Omega^{\rm spin}_2\cong \Z/2 $.
\enddemo
\vfil
Recall that $W^{(n)}$ denotes the regular covering of $W$ which corresponds to
the $n^{\rm th}$ term $\pi_1(W)^{(n)}$ of the derived series of $\pi_1(W)$. An {\it  $(n)$-surface} is by definition a
$\pi_1(W)^{(n)}$-surface in the sense of Definition~\ref{def:$N$-surface}. In Section~\ref{sec:surfaces}, we explained the
quadratic form
$\lambda_n,\mu_n$ on $H_2(W^{(n)})$ in terms of intersection and
self-intersection numbers of $(n)$-surfaces in $W$.

\numbereddemo{Definition} \label {def:Lagrangian} Let $W$ be an $H_1$-bordism such that
$\lambda_0$ is a hyperbolic form.
\begin{itemize}
\item[1.] A {\it  Lagrangian} for $\lambda_0$ is a direct summand of  $H_2(W)$ of half
rank on which $\lambda_0$ vanishes.
\item[2.]    An {\it  $(n)$-Lagrangian} is a submodule $L\subset   H_2(W^{(n)})$ on
which
$\lambda_n$ and $ \mu _n$ vanish and which   maps onto a Lagrangian of the
hyperbolic form $\lambda_0$ on  
$H_2(W)$.
\item[3.]   A {\it  spherical Lagrangian} is a submodule
$L\subset\pi_2(W)$ on which $ \lambda, \mu $ vanish and which maps onto a  
Lagrangian of $ \lambda_0$.
\item[4.] Let $k \leq n$. We say that an $(n)$-Lagrangian $L$ admits   {\it  $(k)$-duals} if
$L$ is generated by
$(n)$-surfaces $\l_1, \dots, \l_g$ and there are
$(k)$-surfaces $d_1, \dots, d_g$ such that $H_2(W)$ has rank
$2g$ and \end{itemize}
$$
\lambda_k(\l_i,d_j) = \delta_{i,j}.
$$ 
\noindent \hangindent = 28pt \hangafter=0
Similarly, {\it  spherical duals} for $L$ are classes $ d_1,
\dots,   d_g \in \pi_2(W) $ satisfying the above equation for
$k=n$.
\enddemo
\eject

\proclaim{Theorem} \label {half} Let $M$ be a closed oriented $3$\/{\rm -}\/manifold and $n\in
\N_0$. Then the following statements are equivalent\/{\rm :}\/ There is an
$H_1$-bordism {\rm ...}
\begin{itemize}
\item[{\rm (i)}]  {\rm ...} which contains an $(n+1)$\/{\rm -}\/Lagrangian with
$(n)$\/{\rm -}\/duals.
\item[{\rm (ii)}] {\rm ...} which contains a spherical Lagrangian with
$(n)$\/{\rm -}\/duals.
\item[{\rm (iii)}] {\rm ...} which contains a spherical Lagrangian admitting a   Whitney
tower of height~$(r.5)$ and with
$(n-r)$\/{\rm -}\/duals for some  $r\in\{0,\dots,n\}$.
\item[{\rm (iv)}] {\rm ...} which contains a spherical Lagrangian admitting a   Whitney tower
of height~$(n.5)$.
\end{itemize}
\endproclaim

\numbereddemo{Definition}\label{def:n.5-solv}  The $3$-manifold $M$ is called {\it  $(n.5)$-solvable} if the   conditions above are satisfied. If $M$ is the
0-surgery on a knot or   a link then the corresponding knot or link is called
$(n.5)$-solvable (and the link has trivial linking numbers, so that $H_1$ is a
free abelian group on the number of components of the   link).
\enddemo

This agrees with the definition given in the introduction.
\numbereddemo{{R}emark} \label{rem:homology cobordism} It is clear that this notion is
invariant under homology   cobordisms. More precisely, assume that $M$ and $M'$
form the   boundary of a $4$-manifold $W$ such that the two inclusions induce  
isomorphisms on
$H_{\ast}$. Then $M$ is
$(n.5)$-solvable if and only if $M'$ is $(n.5)$-solvable. For the proof one glues
together the obvious $4$-manifolds.
\enddemo

\demo{Proof of Theorem~{\rm \ref{half}}}
$ {\rm (ii)} \Rightarrow {\rm (i)}$ is trivially true.

$ {\rm (i)} \Rightarrow {\rm (ii)} $  By Lemma~\ref{lem:l} we may assume that   we have
disjointly embedded framed $(n+1)$-surfaces
$\l_1, \dots,\l_g$. Moreover, the {\it  geometric} intersections with the
$(n)$-duals
$d_1, \dots, d_g$ are
$\delta_{i,j}$ and the duals $d_i$ may be assumed to have trivial normal bundle
since
$W$ is spin. Consider a {\it  standard collection} of simple closed curves
$\alpha_{r,s}$ on $\l_s$. By definition, these are simple closed curves which
represent a basis of $H_1(\l_s)$ such that the algebraic and geometric
(self)-intersections agree. By assumption, there are $(n)$-surfaces $A_{r,s}$
whose boundaries are the curves
$\alpha_{r,s}$. Note that the orientations of the curves
$\alpha_{r,s}$ give a nonvanishing vector field in the normal bundle of
$\alpha_{r,s}$ in $\l_s$. After some boundary twists (see
\cite[p.16]{FQ}) we may assume that this vector field extends to a nonvanishing
vector field for the normal bundle of $A_{r,s}$ in $W$. In this case we may refer
to   the surfaces $A_{r,s}$ as {\it  framed}. By tubing into the
$(n)$-duals $d_i$ we may achieve that the interiors of the
$A_{r,s}$ are disjoint from all $\l_j$. This  preserves the framing on the
$A_{r,s}$.\eject

Now consider tangential push-offs $ \alpha_{r,s}'\subset A_{r,s}$ of 
$\alpha_{r,s}$. These circles have a normal 2-frame on them, one vector field 
pointing into $A_{r,s}$, the other being the nonvanishing normal vector field on
$A_{r,s}$ restricted to $\alpha'_{r,s}$. We do surgery on all
$\alpha_{r,s}'$ such that the 2-frames extend over the new 2-disks $b_{r,s}$.
More precisely, we cut out small neighborhoods of $\alpha'_{r,s}$ homeomorphic to
$S^1
\times D^3$ (disjoint from $d_i$ and $\l_j$) and   add copies of $D^2 \times 
S^2$ using the 2-frames to identify the   boundaries $S^1 \times  S^2$. Denote by
$S_{r,s}$ the disjointly embedded framed $2$-spheres $0\times  S^2$. Every
surgery changes
$H_2(W)$ by the orthogonal sum   with a hyperbolic form on $S_{r,s}$ and
$$ B_{r,s}: = A_{r,s}\cup _{\alpha'_{r,s}} b_{r,s}.
$$ Denoting by $W'$ the result of all these surgeries, we see that it still is an\break
$H_1$-bordism. Moreover, we claim that
$W'$ has a spherical Lagrangian: We may use   two parallels of the disks
$b_{r,s}$ to do symmetric surgery on the  
$(n+1)$-surfaces~$\l_s$. This operation is also called a {\it  contraction} in
\cite[p.34]{FQ}. Call the resulting disjointly   embedded $2$-spheres $L_1,
\dots, L_g$. Then the collection of $2$-spheres $L_j, S_{r,s}$ form a spherical  
Lagrangian because the only geometric intersections among these 
$2$-spheres are two points of intersection between $L_s$ and  $S_{r,s}$ for each
$s=1,
\dots, g$ and each $r$. But these intersections are algebraically trivial because
they can be paired   up by small ribbon Whitney disks (see the figure in 
\cite[p.35]{FQ}). By construction, the $(n)$-surfaces
$d_j,B_{r,s}$ form (geometric)   duals for these $2$-spheres and it is clear that
they together   generate $H_2(W')$.

Note that statement~(ii) is the case $r=0$ and statement~(iv) is the case $r=n$
in statement~(iii). Therefore, to prove the equivalence of (ii), (iii) and (iv),
it suffices to prove two induction steps for statement~(iii), one increasing
$r$, the other decreasing~$r$. 

\vglue4pt {\it The induction step $r\mapsto r-1$}.
 Applying Lemma~\ref{lem:l} to the
$(n-r)$-duals $d_1, \dots, d_g$ we may   assume that their {\it  geometric}
intersection with the framed immersed
$2$-spheres $\l_1, \dots, \l_g$ is $\delta_{i,j}$. In fact, by pushing down
intersections between $d_i$ and Whitney   disks in the tower (introducing many
algebraically canceling pairs   of intersections between $d_i$ and $\l_j$) we may
assume that each $d_i$ intersects the whole tower in a single point. Let
$\alpha_k$ be parallels of the bottom stage Whitney circles such that
$\alpha_k$ lie on the interior of the Whitney disks $  \Delta_k$ of the first
collection $ \CC_1$ in the Whitney tower. Picking one of the two double points
that correspond to $  \Delta_k$,  we get Clifford tori $T_k$ that are disjoint
from all $d_i$ and  
$\l_j$ and intersect the Whitney tower in exactly one point an  $  \Delta_k$.
Both standard circles on $T_k$ are by definition the meridians of the sheets that
are intersecting at that point. By construction, these sheets have
$(n-r)$-duals in
$W\smallsetminus \cup_j\l_j$ and thus the $T_k$ are disjoint $(n-r+1)$-surfaces
in this
$4$-manifold. We now do surgeries on the curves $\alpha_k$. As above these
produce disks $b_k$ which are useful in two respects:  They can be used to do
Whitney moves of
$\l_j$ which make these   disjointly embedded spheres which can be thus surgered
away. Call   the resulting $4$-manifold $W'$. Then the unions $ 
\Delta_k \cup   _{\alpha_k} b_k$ form a spherical Lagrangian in $W'$ which admits
a   Whitney tower of height
$(r-0.5)$ (formed from the upper stages of   the original Whitney tower).
Moreover, these $2$-spheres have  
$(n-r+1)$-duals $T_k$.  Note that the $T_k$ form an $(n-r+1)$-Lagrangian.
\vglue4pt

{\it  The induction step $r\mapsto r+1$}.  By
Remark~\ref{rem:hyperbolic} we may assume that the
$(n-r)$-duals   satisfy $\lambda(d_i,d_j)= \mu(d_i)=0$. More precisely,  this involves summing
the   original $(n-r)$-duals with combinations of the framed $2$-spheres  
$\l_j$. This preserves the property that the $d_i$ are  
$(n-r)$-surfaces (and also the property $\lambda (\l_i,d_j) =  
\delta_{i,j}$). Applying Lemma~\ref{lem:l} to the new
$(n-r)$-duals $d_1, \dots,   d_g$ we may assume that each
$d_i$ intersects the Whitney tower in   exactly one point on
$\l_i$ and that the $d_i$ are represented by   disjointly embedded framed
$(n-r)$-surfaces. Let $\alpha_{r,s}$ be a standard collection of simple closed
curves for $d_s$. By assumption, there are $(n-r-1)$-surfaces $A_{r,s}$ with
boundary
$\alpha_{r,s}$. As in the proof for ${\rm (i)}\Rightarrow {\rm (ii)}$ we can arrange that the
$A_{r,s}$ are framed and have interiors   disjoint from $d_i$. We again do
surgeries on tangential push-offs $  
\alpha'_{r,s} \subset A_{r,s}$ of $\alpha_{r,s}$. Then we do symmetric surgery on
the
$d_i$ to obtain disjointly   embedded framed $2$-spheres $D_1, \dots, D_g$. As
before there are disjoint $2$-spheres $S_{r,s}$ resulting from   each surgery.
They have geometric $(n-r-1)$-duals
$B_{r,s}$ made by   closing off the $A_{r,s}$ with the cores of the $2$-disks
attached.   Recall that the intersections between $D_g$ and $ S_{r,s}$ are  
paired up by ribbon Whitney disks. This time we actually do the corresponding
Whitney moves to make $D_j$ disjoint from $S_{r,s}$ (and keep them disjoint from
$B_{r,s}$). The cost of these last Whitney moves is that the $2$-spheres
$S_{r,s}$ now intersect in pairs, corresponding to the intersections
$\alpha_{r,s} \cap \alpha_{r',s}$. But these intersections again occur in pairs
with disjointly embedded Whitney disks
$  \Delta_{k,s}$ (see Figure~8.1).  
Each Whitney disk $  \Delta_{k,s}$ intersects the contraction
$D_s$ in a single point (on the central square) which we remove by summing   into
the original Whitney tower (which is dual to $d_i$ and hence to $D_i$). Finally, we do surgery on   the
$2$-spheres $D_1,
\dots, D_g$ to obtain our $4$-manifold $W'$. By construction, $H_2(W')$ is
generated by the $2$-spheres $S_{r,s}$, which admit a Whitney tower of height
$(r+1.5)$, and   their
$(n-r-1)$-duals $B_{r,s}$. \enddemo
 \vglue-20pt
\centerline{\BoxedEPSF{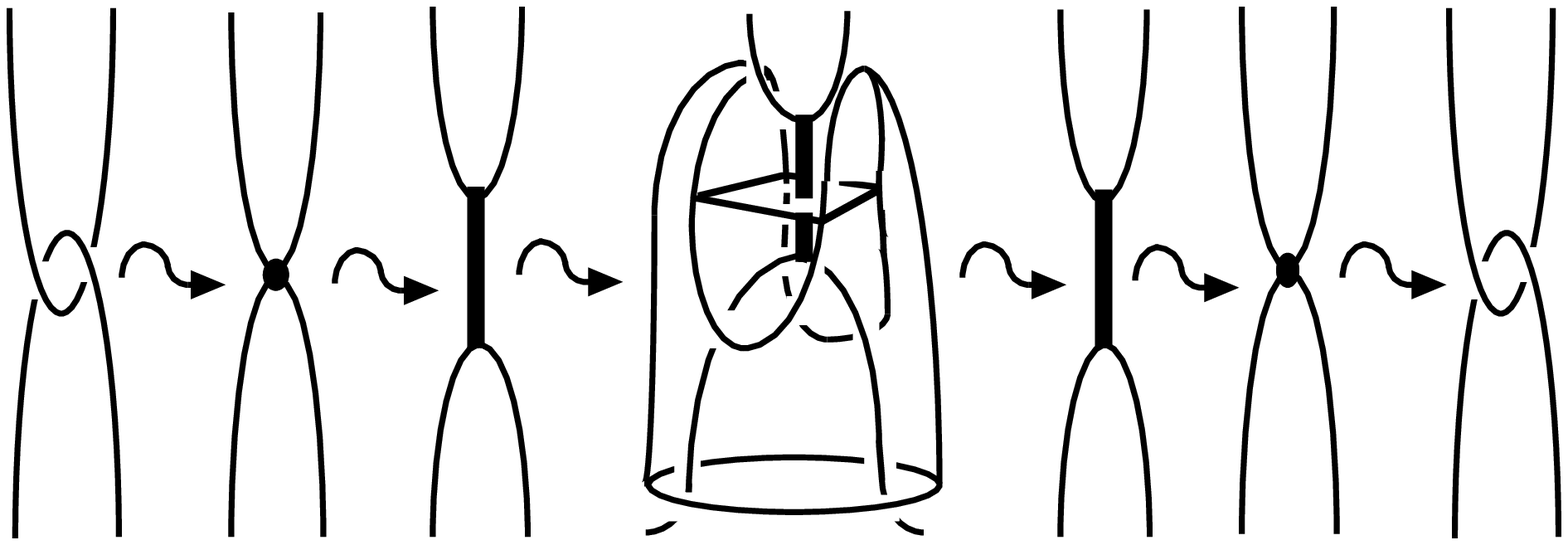 scaled 350}}
\vglue-20pt
\centerline{Figure 8.1. 
       $\Delta_{k,s}$ is the union of the thick arcs in this display.} 
\eject

Looking back at the four statements in the definition of  
$(n.5)$-solvable\break  $3$-manifolds, we see that there is an obvious candidate for what an
$(n)$-solvable $3$-manifold should be.

\numbereddemo{Definition}\label{def:n-solv} A $3$-manifold   $M$ is $(0)$-solvable if it bounds an $H_1$-bordism $W$
such that $(H_2(W),\lambda_0)$ is hyperbolic.  A $3$-manifold $M$ is
$(n)$-solvable, $n>0$,  if any of the conditions of Theorem~\ref{ndashsolv} below are
satisfied.  A link is
$(n)$-solvable if $0$-surgery on the link is an
$(n)$-solvable $3$-manifold.
\enddemo

\proclaim{Theorem} \label{ndashsolv} Let $M$ be a closed oriented $3$\/{\rm -}\/manifold and
$n\in \N$. Then the following statements are equivalent\/{\rm :} There is an
$H_1$\/{\rm -}\/bordism {\rm ...}
\begin{itemize}
\item[{\rm (i)}]  {\rm ...} which contains an $(n)$\/{\rm -}\/Lagrangian with
$(n)$\/{\rm -}\/duals.
\item[{\rm (ii)}] {\rm ...} which contains an $(n)$-Lagrangian with spherical duals.
\item[{\rm (iii)}] {\rm ...} which contains a spherical Lagrangian admitting a   Whitney
tower of height~$(r)$ and with
$(n-r)$\/{\rm -}\/duals for some $r\in\{1,\dots,n\}$.
\item[{\rm (iv)}] {\rm ...} which contains a spherical Lagrangian admitting a   Whitney tower
of height~$(n)$.
\end{itemize}

\endproclaim

\vglue-9pt
{\it Proof}. The arguments that (ii), (iii) and (iv) are equivalent are exactly
as in Theorem~\ref{half}. 
One only needs to make sure that statement (ii) is really equivalent to the $r=1$ case of statement (iii).  In one direction
one uses the $(n)$-Lagrangian from (ii) as the starting point in the induction step $r\Rightarrow r+1$ in the proof of
Theorem 8.4 (there one actually turns $(n)$-duals into an\break $(n)$-Lagrangian first, which is not needed here).  The
output is exactly the $r=1$ case of statement (iii). Conversely, we already observed at the end of the induction step
$r\Rightarrow r-1$ in the proof of Theorem 8.4 that the Clifford tori $T_k$ form an $(n)$-Lagrangian if one begins with
$(n-1)$-duals.  Thus this step gives exactly (ii).

\vglue8pt

(i) {\it implies} (ii).  Let $\l_1, \dots, \l_g$ be an $(n)$-Lagrangian in $W$ with $(n)$-duals\break $d_1,
\dots, d_g$. By Remark~\ref{rem:hyperbolic} and   Lemma~\ref{lem:l} we may assume
that all
$\l_i,d_j$ are represented   by framed embeddings and that the only geometric
intersections among   these $(n)$-surfaces are single points of intersections
$p_i =
\l_i \cap d_i$ for $i=1, \dots, g$. Now this is a perfectly symmetric   setup and
thus we also do a symmetric construction. We do abstract surgery on standard
collections of simple closed curves  $\alpha_{r,s}$ on $\l_i$ {\it  and} $d_j$.
Then we contract to get a geometrically hyperbolic collection of
$2$-spheres $L_1, \dots, L_g,   D_1,
\dots, D_g$. We push the $2$-spheres $S_{r,s}$ off the   contraction, introducing
pairs of double points with Whitney disks   which intersect $L_i$ or $D_j$ in a
single point. We remove this point by summing into the dual
$2$-sphere $D_i$   respectively
$L_j$. This introduces many intersections among the   Whitney disks which will
not be relevant. Finally, we do surgery on the $2$-spheres (say) $L_1, \dots,
L_g$ to obtain a
$4$-manifold $W'$. It contains a spherical Lagrangian $S_{r,s}$ with Whitney  
disks disjoint from these spheres. Therefore, we have actually constructed a
Whitney tower of height~$1$. As discussed in the proof of Theorem~\ref{half} the
$S_{r,s}$ have geometric
$(n-1)$-duals $B_{r,s}$ (using the fact that we started out with
$(n)$-surfaces). We have thus shown that statement (i) implies statement (iii)
with  $r=1$. But this is equivalent to statement (ii).
\hfill\qed\vglue12pt

We next show that there are many $(h)$-solvable knots.

\proclaim{Theorem} \label{surgery curves} If there exists an $(h)$\/{\rm -}\/solvable link $L$ which
forms a standard half basis of untwisted curves on a Seifert surface $F$ for a knot
$K${\rm ,} then $K$ is $(h+1)$\/{\rm -}\/solvable.
\endproclaim

\centerline{\BoxedEPSF{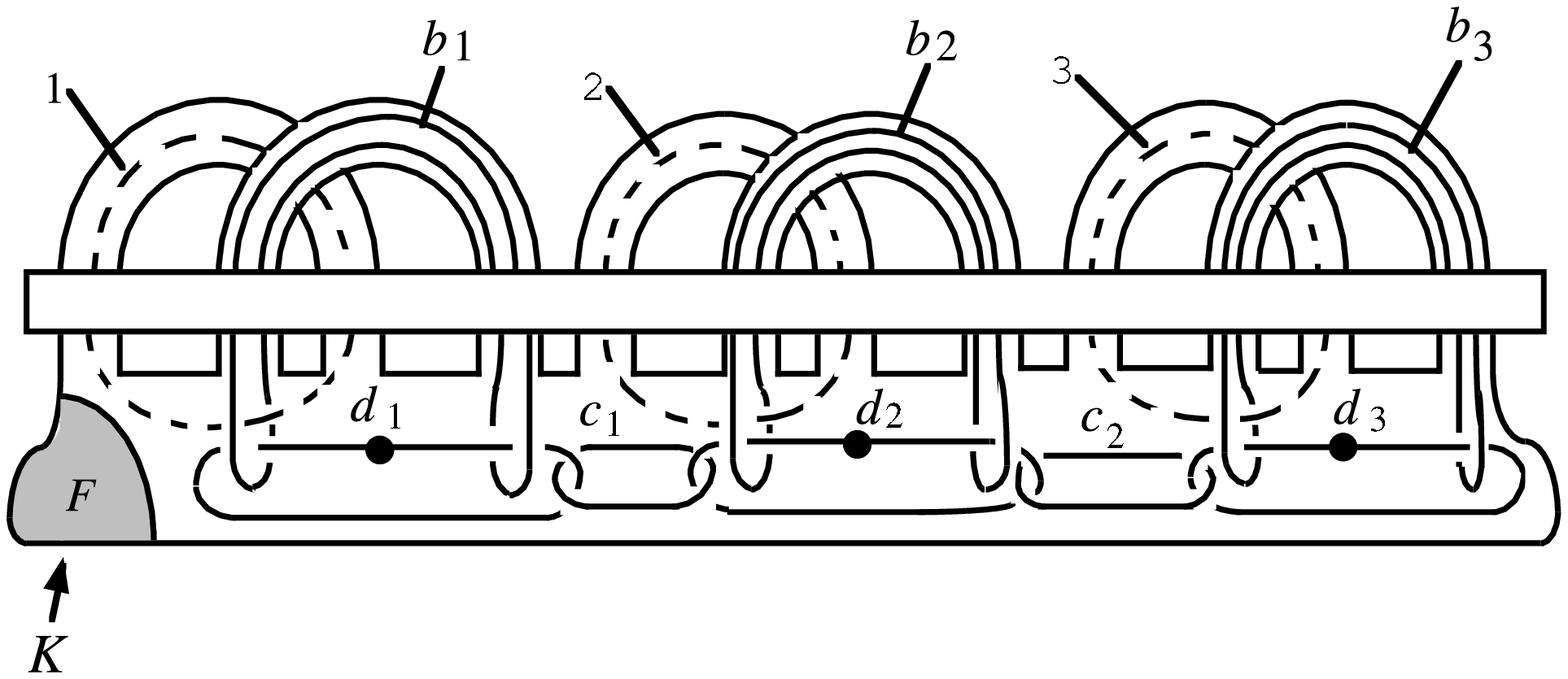 scaled 450}}

 \centerline{Figure 8.2. The cobordism $C$.}

\demo{Proof}  In Figure~8.2  we have drawn the Seifert surface $F$
in the case of genus
$g=3$. It shows a box containing an arbitrary string link of (possibly twisted)
bands for $F$.
The dashed lines $\l_i$ denote the link $L$ whose meridians are   called
$m_i$. In addition, we drew $g$ dotted unlinked circles $d_i$ (whose meridians we
call
$t_i$) and
$g$ solid circles $b_i$ going around the dual bands to $L$. Finally, there are
solid circles $c_1, \dots, c_{g-1}$ connecting pairs of dotted circles. The
figure   determines a
$4$-manifold $C$ in the following way: Start with the   lower boundary
$\partial_-C=\SS^0L$, which is also the   boundary of a 0-handle together with
0-framed 2-handles on the  $\l_i$. To $\partial_-C \times I$ we attach $g$
1-handles corresponding to the dotted circles and $(2g-1)$ 0-framed 2-handles
along $b_i$ and 
$c_j$. Thus the relative chain-complex $C_{\ast}(C,\partial_-C)$ has only   terms
for
${\ast}=1,2$ and the boundary map
$$
\Z^{2g-1} \cong\langle b_i, c_j \rangle =C_2 \xrightarrow{\partial }
          C_1= \langle t_i \rangle \cong \Z^g
$$ satisfies $\partial (b_i)=0, \partial (c_j)=t_{j+1}-t_j$. Therefore,
$H_1(C,\partial_-C)\cong \Z$ is generated by any of the  
$t_j$ and there is an isomorphism
$$ H_2(C,\partial_-C) = H_2(C,\SS^0 L)\xrightarrow{\cong} H_1(\SS^0 L)= \langle
m_i
\rangle \cong \Z^g.
$$
The upper boundary $\partial_+C$ is given by 0-framed surgery on   all the
circles in Figure~8.2, i.~e.\  the  
$\l_i,b_i,c_j,d_i$. In Lemma~\ref{lemcopl} below we show that  
$\partial_+C\cong \SS^0K$. But already from the figure above it follows that the
inclusion   induced map $H_1(\partial_+C) \ra H_1(C,\partial_-C)$ is an  
isomorphism.

Summarizing the above construction, we have a cobordism $C$ between $\SS^0L$ and
$\SS^0K$ which has the following properties
\begin{itemize}
\item[1.] $\partial_-:H_2(C,\SS^0L) \ra H_1(\SS^0L)$ is   an isomorphism.
\item[2.] $i_+:H_1(\SS^0K)\ra H_1(C,\SS^0L)$ is an   isomorphism.
\end{itemize} 
Now recall that $L$ is $(h)$-solvable and let $V$ be the  
$H_1$-bordism for $\SS^0L$ which contains a
$(k)$-Lagrangian  with $(n)$-duals. Here $k=n$ if $h=n
\in\N$ and $k=n+1$ if $h=n.5$ (see Definition~\ref{def:Lagrangian}). Define
$W:=V\cup_{\SS^0L}C$ which is a $4$-manifold with   boundary $\SS^0K$. Consider
the long exact sequence for the pair $(W,V)$, noticing   that by excision
$H_{\ast}(W,V) \cong H_{\ast}(C,\SS^0L)$.
$$ 0 \to H_2V \to H_2W \to H_2(C,\SS^0L) \to H_1V \to H_1W \to   H_1(C,\SS^0L)
\to 0.
$$ By assumption, $H_1(\SS^0L) \cong H_1V$ and therefore the   boundary-map
$H_2(C,\SS^0L)\break \ra H_1V$ is an isomorphism by   1.\ above. This implies that we
have isomorphisms.
$$ H_2V \cong H_2W \hbox{  and } H_1W \cong H_1 (C,\SS^0L).
$$ By 2.\ above this shows that $H_1(\SS^0K) \cong H_1W$. Since $H_1V \ra H_1W$
is the zero map any
$(r)$-surface in $V$ is actually an $(r+1)$-surface when considered in $W$.
Therefore, we actually have a
$(k+1)$-Lagrangian with $(n+1)$-duals in $W$, using   the isomorphism $H_2V
\cong H_2W$. But this shows that $K$ is  
$(h+1)$-solvable.
\enddemo

\proclaim{Lemma} \label{lemcopl} In the above setting{\rm ,} there is a diffeomorphism $\partial_+C
\cong\SS^0K$.
\endproclaim

\demo{Proof}  We first isotope Figure~8.2 into the position of  
Figure~8.3, keeping the dashed, solid and dotted convention   even though
all circles are considered as $0$-framed $2$-handles.
Now we slide each $d_i$ twice over its partner
$\l_i$ leading to the handle diagram
in which $b_i$ are geometric duals for $\l_i$. Thus we may cancel these handles in
pairs, effectively erasing them from the diagram (see Figure~8.4).
But now it is clear that $(g-1)$ more cancellations involving the $c_j$ lead to
$\SS^0K$, by the fact that the $\l_i$  were untwisted.
\enddemo

\figin{Fig83}{500} \vglue4pt
\centerline{Figure 8.3. A $0$-framed handle decomposition for
$\SS^0K$.} 
\vfill
\figin{Fig84}{500} \vglue4pt
       \centerline{Figure 8.4. Before the last $(g-1)$ cancellations.} 
 \vglue12pt

We next show that the slightly abstract notion of $(h)$-solvability is implied by
very concrete geometric conditions in the $4$-ball. The first condition is in
terms of gropes and the second in terms of Whitney towers.

\proclaim{Theorem} \label{gropes} If a knot $K$ bounds a grope of height~$(h+2)$ in
$D^4$ then
$K$ is $(h)$\/{\rm -}\/solvable.
\endproclaim

{\it Proof}.
 Let $\alpha_1, \dots,  \alpha_{2g}$ be a standard collection of   simple closed
curves on the bottom stage $F$ of the grope
$G$ of height~$(n+2)$ which bounds   our knot $K$. (Note that in this proof the
$\alpha_i$ form a {\it  full} basis of   curves rather than just a half-basis as
in the previous proof.) As in Theorem~\ref{half} we do surgery on tangential
push-offs 
$\alpha'_i\subset A_i$ of $\alpha_i$, where $A_i$ are the second surface stages
of our grope $G$. As before this leads to $2$-disks $b_i$, $2$-spheres $S_i$ and 
$(n)$-duals
$B_i=A_i\cup b_i$. Again we use two parallels of the $b_i$ to do symmetric
surgery on
$F$ to obtain a disk
$D$ bounding our knot $K$. Finally, we push the $S_i$ off the contraction $D$ and
remove the   interior of a thickening of
$D$ to obtain a $4$-manifold $W$. By construction, $\partial W$ is 0-surgery on
$K$ and
$H_2(W)$ is freely generated by
$S_i$ and their geometric $(n)$-duals $B_i$.   Therefore, $W$ satisfies condition
(ii) of Theorem~\ref{ndashsolv} and we are done for gropes of integral
height.

If the height of the grope is a half integer $h=n+2.5$ then we may   pick a half
basis
$\alpha_1, \dots, \alpha_g$ which bounds gropes of   height $(n+2)$ and such that
the dual half basis $\alpha_{g+1},\dots, \alpha_{2g}$ bounds gropes of height
$(n+1)$. Then the above   construction gives a
$4$-manifold $W$ with framed embedded  $2$-spheres $S_1,
\dots, S_{2g}$ with geometric $(n+1)$-duals $B_1,\dots, B_g$, respectively
$(n)$-duals
$B_{g+1}, \dots, B_{2g}$. By construction, the surfaces $B_i$ are framed and
embedded   disjointly. The deficiencies in this family of surfaces are pairs of  
intersection points between $S_i$ and $S_{i+g}$ coming from pushing   these
spheres off the contraction $D$. However, we may remove  all of   these intersections
by tubing each $ S_{i+g}$ twice into parallel   copies of $B_i$ for each
$i=1,\dots,g$. Since these
$B_i$ are $(n+1)$-surfaces, we obtain an  $(n+1)$-Lagrangian $S_1,\dots, S_g,
S'_{g+1},\dots, S'_{2g}$ with  $(n)$-duals $B_1,\dots, B_{2g}$. Thus
condition~(i) from Theorem~\ref{half} is satisfied.
\hfill\qed

\proclaim{Theorem} \label{Whitney towers} If a knot $K$ bounds a Whitney tower of
height~$(h+2)$ in
$D^4$ then $K$ is $(h)$\/{\rm -}\/solvable.
\endproclaim

{\it Proof}. 
 Let $T$ be a Whitney tower of height~$(h+2)$ in $D^4$ which bounds   our knot
$K$. Consider the Whitney circles
$\alpha_i$ for the   immersed 2-disk $  \Delta$ that bounds the knot $K$. We do
surgery on tangential push-offs $\alpha'_i
\subset A_i$ of  $\alpha_i$, where $A_i$ denote the next stages of the Whitney
tower 
$T$. This leads to external\break 2-disks $b_i$ which can be used in two ways.  We may
do Whitney moves along $b_i$ to change $  \Delta$ into an   embedded 2-disk
$D$ giving us a $4$-manifold $W$ which is the surgered $D^4$ minus   an open
neighborhood of $D$. Then
$\partial W=\SS^0K$ and $H_2W$ has a Lagrangian which is generated by
$2$-spheres
$S_i:=A_i\cup b_i$ in $W$ which allow a   Whitney tower of height~$h$, formed
from the upper stages of the   Whitney tower $T$.
\hfill\qed\vglue4pt

The last proof   still owed is of the fact that a knot is algebraically slice if and
only if it bounds a grope of height 2.5 in $D^4$ (see Theorem~\ref{algebraically
slice}). Recall from Remark~\ref{rem:$(n)$-solvable} that a knot $K$ is
algebraically slice if and only if it is $(0.5)$-solvable. So by
Theorem~\ref{gropes} it suffices to prove the following result.

\vglue-22pt
\phantom{whow}

\proclaim{Theorem} \label{2.5} If a knot $K$ is algebraically slice then it bounds a
grope of height~$2.5$ in $D^4$.
\endproclaim

\vglue-10pt
{\it Proof}. Using the Levine condition from Theorem~\ref{algebraically slice},
we may start with a Seifert surface $F$ for $K$ in $S^3$ with a half-basis of
curves $
\alpha_1,
\dots, \alpha_g$ with vanishing linking numbers. This implies the existence of
disjointly embedded framed surfaces  $A_i$ in $D^4$ with
$\partial A_i=\alpha_i$. Let $\beta_1, \dots,  \beta_g$ be a dual half-basis of
curves on $F$. By a base change as in Remark~\ref{rem:hyperbolic} we may assume  
that the self-linking numbers of all $b_j$ are even. Then there are   framed
embedded surfaces
$B_j$ in $D^4$ with $\partial B_j=\beta_j$.  

Let $\gamma_{r,s}$ be a full basis of curves on $A_s$. Then, after   some
boundary twists, there are framed embedded surfaces $G_{r,s}$   in the interior
of $D^4$ with
$\partial G_{r,s}=\gamma_{r,s}$. We may assume that $A_i,B_j$ and
$G_{r,s}$ only intersect in isolated points away from their boundaries. Note also
that all these surfaces have interiors disjoint from the   Seifert surface
$F$. We now push the Seifert surface slightly into  $D^4$, or more precisely, we
add a small collar
$S_3 \times  I$ to  $D^4$ with an annulus which connects the original knot $K$ to
the new boundary
$S^3$. In any case, we now see that $K$ bounds a framed grope of height~$2.5$ in
$D^4$ whose bottom stage $F$ is disjoint from all   other surfaces stages. We
shall show next that such a grope can be improved to a framed   embedded grope
of the same height (then thickening this framed grope leads to the desired grope).
The first step is to push down   all intersections among
$G_{r,s}$ and with
$B_j$ into $A_s$ (see  \cite [\S 2.5]{FQ}). This makes the (interiors of)
$G_{r,s}$   disjoint, and also disjoint from $B_j$. Let $T_{\alpha_j}$ denote the
$2$-tori which are the normal circle bundles to $F$ restricted to $\alpha_j$.
They are disjointly embedded, framed and may be assumed to be disjoint from all
other surfaces, except a single point of intersection $T_{\alpha_j} \cap   B_j$.
But this means that we may use tubes into the $T_{\alpha_j}$ to remove all
intersections among the $B_i$. Note that this increases the genus of each
$B_i$ but we do not care.

Finally, consider normal tori $T_{\beta_i}$ to $\beta_i$. Again we may assume
that they are disjointly embedded, framed and   disjoint from everything else
(including the
$T_{\beta_i}$), except   a single intersection point
$T_{\beta_i} \cap A_i$. Thus tubing into $T_{\beta_i}$ removes the last
intersections,   namely those between $B_j$, respectively $G_{r,s}$ and~$A_i$.
Again this procedure increases the genera of the surfaces $B_j$ and  
$G_{r,s}$ but since they form the top stage of the grope, this is irrelevant.
\hfill\qed

\numbereddemo{{R}emark}   \label{rem:2} The above procedure can be used to show that any knot
with trivial   Arf invariant bounds a framed embedded grope of height~$2$ in
$D^4$.   The difficulty in increasing the height by $0.5$ lies in the fact that
if one has to use the tori $T_{\beta_i}$ to remove intersections {\it  among} the
$A_i$ then one cannot find the next stage surfaces $G_{r,s}$ : One of the curves
on each
$T_{\beta_i}$ is by construction the   meridian to the pushed-in Seifert surface
$F$ and is therefore not   null-homologous in $D^4 \smallsetminus F$.
\enddemo

\vglue-16pt
 \section{Casson-Gordon invariants and solvability of knots}\label{sec:CG}
\vglue-6pt

In this section we review the Casson-Gordon invariants, and show they vanish on
$(1.5)$-solvable knots. Throughout this section, all chain and cochain complexes
are cellular, with the cellular structure obtained from lifting a cellular
structure on the base.

\demo{Seifert pairings{\rm ,} linking pairings and $(.5)$\/{\rm -}\/solvability}
We recall the definition of the Seifert pairing and the classical knot slicing
obstructions due to Levine~\cite{L1}. Let $F \subset S^3$ be a Seifert surface for
the knot $K$. The Seifert pairing on $H_1(F)$ is defined by 
$$
\theta(x, y) = \l k(\i_+ x, y)
$$ where $\l k$ is the usual linking in $S^3$ and $\iota_+$ is the positive
push-offs in the normal direction from $F$.  Following Kervaire
\cite{K} and Stoltzfus~\cite{Sto}, we define an isometric structure 
$$ s : H_1(F; \Z) \to H_1(F; \Z)
$$ by the equation $\theta(x, y) = \langle sx, y \rangle_F$ for all $y \in H_1(F;
\Z)$, where $\langle \ , \
\rangle_F$ is the intersection pairing   on $H_1(F; \Z)$. \enddemo

\numbereddemo{Definition}\label{metabolizer}  A metabolizer for the isometric structure on
$H_1(F; \Z)$ is an $s$-invariant direct summand $H \subset H_1(F; \Z)$ such that 
$$ H = H^{\perp} = \{ y \in H_1(F; \Z) \ | \langle x, y \rangle_F = 0
\mbox{ for all x }\in H\}.
$$
\enddemo

Levine shows~\cite{L1} that if $K$ is slice there exists a summand $H \subset
H_1(F;
\Z)$ such that $rk_{\Z}(H) =
\frac{1}{2}rk_{\Z}(H_1(F; \Z))$ and $\theta (H \times H) = 0$.  It follows that
$H$ is a metabolizer for the isometric structure on $H_1(F; \Z)$ defined
above. (See
\cite[p.~95]{K}.)

If $M$ is 0-framed surgery on $K$, $k = p^r \in \Z$ is any prime power, and
$M_k$ is the
$k$-fold cyclic cover of $M$, then
$H_1(M_k) \cong \Z \times TH_1(M_k)$, where the second summand is the
$\Z$-torsion subgroup.  This latter summand has all torsion relatively prime to
$p$ (see~\cite{CG1}, for instance).  The sequence of isomorphisms
$$ TH_1(M_k) \to TH^2(M_k) \to {\rm Ext}_{\Z}^1(TH_1(M_k);\Z) \to
\Hom_{\Z}(TH_1(M_k);\torsionp)
$$ defines a nonsingular $\torsionp$-valued linking pairing on the torsion
subgroup of
$H_1(M_k)$.  The first isomorphism is Poincar\'e duality, the second is Universal
Coefficients, and the last follows from the long exact
$\mbox{Ext}^{\ast}_{\Z}$ sequence associated to the short exact coefficient
sequence
$$ 0 \to \Z \to \Z_{(p)} \to \torsionp \to 0.
$$ This pairing may be computed by the usual formula.

\proclaim{Proposition}\label{.5}  Let $K$ be a $(.5)$\/{\rm -}\/solvable knot via $W$ and
$F$ a closed Seifert surface for $K${\rm ,} i.e.{\rm ,} a surface union a disk at the core of
the added $2$\/{\rm -}\/handle.
\begin{itemize}
\ritem{1.} There exists a choice of oriented $3$\/{\rm -}\/manifold $R
\subset W$ with boundary
$F$ such that 
$$ H = \Ker{(H_1(F) \to H_1(R)/TH_1(R))}
$$  is a metabolizer for the Seifert pairing on $F$.
\ritem{2.} If $k = p^r$ is a prime power{\rm ,} the $\torsionp$\/{\rm -}\/linking pairing on
$M_k${\rm ,} has a self-annihilating subgroup $P \subset H_1(M_k;\Z)${\rm ,} i.e.{\rm ,} a subgroup
$P$ such that \end{itemize}
$$ P = P^{\perp} = \{ y \in H_1(M_k; \Z) \ | \ \l k(x, y) = 0 \mbox{ for all}
\ x
\in P
\}.
$$
\endproclaim 

\demo{Proof}  The second statement follows from the first, and as it is not needed
in this paper the proof is omitted.  Assuming $(1)$-solvability we will prove
Statement~$2$, and more, in Proposition~\ref{QmodZ}.

To prove the first statement, we begin by defining $R$.  By Lemma~\ref{lem:l}, a
basis for the $(.5)$-Lagrangian of $W$ can be represented by disjoint
$(1)$-surfaces $\{ F_i \}$.  We will show $R$ may be chosen to be disjoint from
the surfaces $F_i$. Since $W$ is an $H_1$-bordism, and by transversality,
there is an oriented  $3$-manifold $R' \subset W$ with $\partial R' = F$.  
Now, $R^{\prime}$ intersects $F_1$ in a $1$-manifold. This $1$-manifold is nulhomologous in $F_1$ since $F_1$  lifts to
the universal abelian cover, and since $R'$ is dual to the meridian generating
$H_1(W)$.  A nulhomologous
$1$-manifold bounds a nested collection of subsurfaces on
$F_1$, which, having boundary, have trivial normal bundles in $W$.  Perform ambient
surgeries on
$R^{\prime}$, using these subsurfaces, to remove the intersections of $R'$ with
$F_1$.  More precisely, one successively removes  the trivial regular
neighborhood, in $\RR'$, of each circle of intersection, and replaces it with the circle
bundle, in $W$, over the subsurface.  If we continue in this manner, the resulting
oriented
$3$-manifold
$R$ has the same boundary as $R^{\prime}$ but does not intersect any surface
$F_i$.

We now show the Seifert pairing vanishes on $H$, i.e., $\theta(H\times H) = 0$.
Given $[x], [y] \in H$, there exist $c_x \in C_2(R) \subset C_2(W)$ and
$d_x \in C_2(M)$ such that
$\partial c_x = \partial d_x = \lambda x$ for some $\lambda \in \Z - \{ 0\}$. 
Since  
$(R
\cup M)
\cap F_i = \emptyset$ for all $i$, the homomorphism $H_2(R \cup M ;\Z) \to
H_2(W;\Z)$ factors through the dual of the $(.5)$-Lagrangian $L^{\perp} = L
\subset H_2(W;\Z)$.  Thus, after adding copies of the $F_i$'s to
$c_x \in C_2(W)$, we may assume $(c_x - d_x) = 0 \in H_2(W;\Z)$.

Now choose $c_y \in C_2(R)$ such that $\partial c_y = \mu y$, for some $\mu
\in
\Z - \{ 0\}$.  Since $c_y$ may be pushed off $R$ and is disjoint from all surfaces
$F_i$, the intersection number $c_x \bullet c_y = 0$.  Thus
$$ 0 = c_x \bullet c_y = d_x \bullet c_y = d_x \bullet \lambda y = \lambda
\mu \cdot
\theta ([x], [y]),
$$ 
and therefore $
\theta ([x], [y]) = 0.$

Finally, $H$ is a summand of $H_1(F)$ since $H_1(F)/H$ is torsion-free by
construction. Also, 
$H$ has $\frac{1}{2}$-rank by the usual duality arguments (see, for
instance,~\cite{L1}).  Thus $H$ is a metabolizer for $K$.
\enddemo

{\it Casson\/{\rm -}\/Gordon invariants.}
 Now recall the definitions and fundamental theorems regarding the Casson-Gordon
invariants of knots.
Let $k = p^r$ and $\l = q^s$, where  $p$ and $q$ are  distinct primes. As before, for $K
\subset S^3$ a knot, let $M_k$ denote the $k$-fold cyclic cover of $M$, where
$M$ denotes $0$-framed surgery on $K$.   Suppose we are given a representation 
\begin{equation}\label{rho}
\rho : \pi_1 (M_k) \to \Z_{\l} \times \Z
\end{equation} such that projection to $\Z$ is onto and such that projecting to
$\Z_{\l}$ sends to zero the cycle in $M_k$ whose image in $M$ is $k$ times the
meridian of $K$. Using standard bordism tools, Casson and Gordon observe there is
an oriented $4$-manifold
$W$ and a representation
$$
\psi : \pi_1 (W) \to \Z_{\l'} \times \Z
$$ such that  $\partial W$ is a disjoint union of copies of $M_k$, $\l'$ is a possibly greater power of $q$, and
such that the restriction of $\psi$ to any component of
$\partial W$ is the representation $\rho$.  Now, $W$ can be chosen so that the number
of boundary components is relatively prime to $p$.

Let $\f = \Q(\zeta_{\l})(t)$, where $\zeta_{\l}$ is a primitive $\l^{\rm th}$ root of
unity, and $t$ is an indeterminant. When we use the $\Z_{\l} \times\Z$ cover of
$W$, there is a Hermitian intersection pairing on the middle dimensional homology
$$
\lambda : H_2(W; \f) \otimes H_2(W; \f) \to \f
$$ via the composition of ring homomorphisms
$$
\Z [ \Z_{\l} \times \Z] \to  \Z[\zeta_{\l}][\Z] \to \f .
$$ The first homomorphism is the quotient homomorphism, and the second is
inclusion to the quotient field.  Since $H_2(M_k ; \f) \cong H_1(M_k; \f)
=0$~\cite{CG2} this pairing represents an element
$[\lambda] \in L_0(\f)$. If $\partial W = m M_k$, $(m , p) = 1$, then
$$
\sigma(K, \rho)  = ([\lambda] - [\lambda_0]) \otimes
\frac{1}{m} \in L_0(\f) \otimes \Z_{(p)}
$$ where $[\lambda_0]$ is the class of the intersection pairing on $H_2(W;
\Q)$ viewed as an element of $L_0(\f)$ via the obvious inclusion of rings with
unit
$\Q \to \f$.

\numbereddemo{Definition}\label{def:CG} \cite{CG2} $\sigma (K, \rho) \in L_0(\f) \otimes
\Z_{(p)}$ is called a  {\it  Casson\/{\rm -}\/Gordon invariant of $K$}.
\enddemo

The genius of Casson and Gordon's work is revealed through the concordance
invariance of their obstructions.  We need a definition.

Since $H^1(M_k; \torsionp) \cong \torsionp \oplus \Hom_{\Z}(TH_1(M_k),
\torsionp)$,  the Mayer-Vietoris sequence
\begin{equation}\label{MV} H^1(M_k; \torsionp) \to \oplus^k H^1(M-F;\torsionp)
\to \oplus^k H^1(F;\torsionp)\hskip.5in
\end{equation} together with the Alexander duality isomorphism $$H^1(M-F;
\torsionp) \cong H_1(F;\torsionp)$$ identifies the
$\torsionp$ valued characters on $TH_1(M_k)$ with a subgroup of $$\oplus^k
H_1(F;\torsionp).$$     By a change of basis, Gilmer identifies this subgroup of
$\oplus^k H_1(F;\torsionp)$ with a subgroup $A^k \subset (H_1(F) \otimes
\torsionp)$. (Gilmer prefers to work with the branched cover of $K$ whose homology
is canonically identified with
$TH_1(M_k)$~\cite{G2}.) 

\numbereddemo{Definition} We call the character $\chi_x$ associated to an element $x \in A^k$
the {\it  character Gilmer associated to $x$}. Note that  $\chi_x$  determines a character
we also denote $\chi_x : \pi_1(M_k) \to \torsionp \times \Z$ defined by sending
the meridian to the element $(0,1) \in \torsionp \times \Z$.  We denote the
associated ring homomorphism $\Z \pi_1(M_k) \to \Z [\torsionp
\times \Z] \to \f$ by $\rho_x$.
\enddemo

Theorem~\ref{CG} below, due to P. Gilmer \cite{G2}, is the most general result
about concordance invariance we know, extending Casson and Gordon's original idea
and results.

\vglue6pt {\it Added in proof\/}: Gilmer has informed us that his proof of 9.5 has a serious gap.  For any fixed slice disk,
the theorem remains true for almost all primes~$q$.  These comments also apply to our 9.9.  Theorem 9.11 remains valid in
this same sense and, in any case, is valid for the original Casson-Gordon invariants (using 9.7 and the proof of 9.11 minus
the first and third \pagebreak  paragraphs).

\proclaimtitle{Gilmer \cite{G2}}
\proclaim{Theorem}\label{CG}  If $K$ is slice{\rm ,} then for any Seifert
surface
$F$ for $K$ there is a metabolizer $H
\subset H_1(F)$ for the isometric structure on $H_1(F)$ having the following {\rm  additional} property\/{\rm :}   For
any prime powers $k = p^r$ and $\l = q^s${\rm ,}
$(p,q) = 1${\rm ,} for any $x
\in A^k \cap (H \otimes (\torsionp))$ of order $\l${\rm ,} and for 
$\rho_x : \Z \pi_1(M_k) \to \f,$ the character Gilmer associated to $x${\rm ,} the
Casson\/{\rm -}\/Gordon obstruction

\centerline{$
\sigma (K,\rho_x) \in L_0(\f) \otimes \Z_{(p)}
$}

\noindent vanishes.
\endproclaim

 The aim of this section is to replace the slice hypothesis of
Theorem~\ref{CG} with $(1.5)$-solvability.

We next do an important dimension count.

\proclaim{Lemma}\label{rankeq}  Let $W$ be an $H_1$\/{\rm -}\/bordism for a knot in $S^3${\rm ,} and
let
$W_k$ be its
$k$\/{\rm -}\/fold cyclic cover ($k = p^r$).  Let $\pi_1(M_k) \xrightarrow{\rho} \Z_{\l}
\times \Z
\to
\f$ be a character as in~{\rm \eqref{rho}} with extension $\pi_1(W_k) \to \f$.  Then 
$$ \dim_{\Q} H_2(W_k; \Q) = \dim_{\f} H_2(W_k;\f).
$$
\endproclaim

\demo{Proof}  We show the following equalities where $\chi^{\Bbb F}(W)$ is the
Euler characteristic of $W$ with coefficients in a field $\Bbb  F$: 
\begin{eqnarray*} \dim_{\Q} H_2(W_k; \Q) &\hskip-6pt =\hskip-6pt&  \Sigma (-1)^i
\dim_{\Q} H_i(W_k;\Q)   = 
\chi^{\Q}(W_k)  \\& \hskip-6pt=\hskip-6pt&  \chi^{\f}(W_k)   =\Sigma (-1)^i \dim_{\f} H_i(W_k;\f)   = 
\dim_{\f} H_2(W_k;\f).
\end{eqnarray*}
 The first follows by an easy computation.  The second, third and
fourth equalities are by definition.  The last equality follows from the observation that
$H_{\ast}(M_k, \f) =0$~\cite{CG2}.   In fact, every $4$-manifold with boundary has
the homotopy type of a $3$-dimensional CW-complex, so $H_{\geq 4}(W_k) = 0$ with
any coefficients. Also,
$H_0(W_k;\f)
\cong H_1(W_k;\f) =0$ by Lemma~$4.5$ of~\cite{CG1} and the proof of the corollary
to Lemma~$4$ of~\cite{CG2}.   Also,
$H_3(W_k;\f)
\cong H^1(W_k,M_k;\f)
\cong \Hom_{\f}(H_1(W_k,M_k;\f);\f) =0$, since $H_1(W_k,M_k;\f) =0$.
\enddemo

  $(1)$-{\it solvability and extending characters.}   
Similarly, to the
$\torsionp$ pairing on $TH_1(M_k)$, there are nonsingular relative homology
linking pairings defined for a
$(1)$-solution $W$ as 
follows:
$$ TH_2(W_k,M_k) \cong TH^2(W_k) \cong {\rm Ext}_{\Z}^1(TH_1(W_k);\Z) \cong
\Hom_{\Z}(TH_1(W_k);
\torsionp)
$$
and
\begin{eqnarray*}
 TH_1(W_k) \cong TH^3(W_k, M_k) &\cong& {\rm Ext}_{\Z}^1(TH_2(W_k, M_k);\Z)\\
& \cong&
\Hom_{\Z}(TH_2(W_k,M_k);\torsionp).
\end{eqnarray*}
 Recall from \cite{CG1} that $H_1(W_k;\Z)$ has no $p$-torsion which is used in
the above isomorphisms. By Poincar\'e duality and universal coefficients, 
$H_2(W_k,M_k;\Z)$ has also no
$p$-torsion.

\proclaim{Proposition}\label{QmodZ} Let $K$ be $(1)$\/{\rm -}\/solvable via $W$ and let $$P =
\Ker{( TH_1(M_k)   \to TH_1(W_k))}.$$   Then $P = P^{\perp}${\rm ,} and a character
$TH_1(M_k) \to \torsionp$ given by $x \mapsto \l k(\cdot, x)$  factors through
$H_1(W_k)$ if and only if $x \in P$.
\endproclaim

  Note that this implies the character $\chi_x : \pi_1(M_k) \to
\Z_{\l} \times \Z$ factors through $\pi_1(W_k)$, and similarly for $\rho_x :
\pi_1(M_k) \to \f$.  Also, note that the second statement of the proposition
depends on the particular choice of self-annihilator
$P$ given by the first statement in the proposition.

\demo{Proof}   We only outline the proof, as much of it follows earlier
arguments given in the paper. Consider the following commutative diagram of groups
and homomorphisms:
$$
\BoxedEPSF{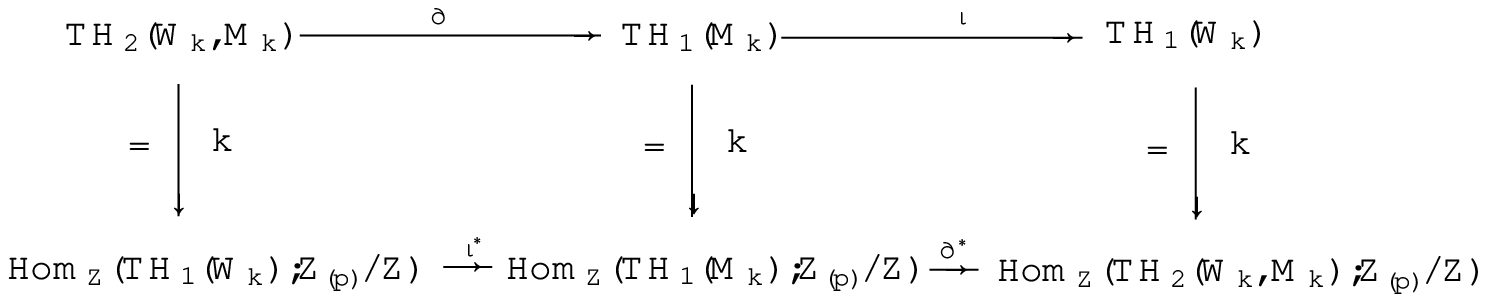 scaled 850}
$$  
That the top horizontal row is exact uses the same argument as in
Lemma~\ref{Nopid} but with $\RR =
\Z[\Z_k]$, and the fact that $rk_\Z H_2(W;\RR) = rk_\Z H_2(W_k;\Z) = 2km$ (see Lemma~\ref{rankeq}).  This
equality follows since the Euler characteristic multiplies in covers.  The
vertical arrows are the nonsingular linking pairings mentioned above.

Assume $x \in P$.  Then $\l k(\cdot,\iota (x)) = \l k(\cdot,0) = 0$. Thus, $\l
k(\cdot,x)$ vanishes on Image$(\partial) = \ker{(\iota)} = P$, and $P \subset
P^{\perp}$. If $x \in P^{\perp}$, then $\partial^{\ast} \circ \l k(\cdot, x) = 0$,
and since the vertical pairings are nonsingular, $\iota (x) = 0$.  Thus
$P^{\perp} \subset P$, and equality follows.

Furthermore, $\l k(\cdot, x)$ extends over 
$H_1(W_k)$ if and only if
\vglue12pt 
\hfill ${\displaystyle (\l k(\cdot, x) \in
\mbox{Image}(\iota^{\ast}) )\Leftrightarrow (\partial^{\ast} \circ \l k(\cdot, x)
= 0)
\Leftrightarrow (\iota (x) = 0 )\Leftrightarrow( x \in P).
}$
\enddemo

\proclaim{{C}orollary}\label{orderext}  Let $\SS$ be the set of characters 
$TH_1(M_k) \to \torsionp$ that extend over
$H_1(W_k)$.  Then the order of $\SS${\rm ,} $|\SS|$ = $\frac{1}{2}\ | T H_1(M_k)|$.
\endproclaim

\demo{Proof}   By Proposition~\ref{QmodZ}, the characters that extend lie in
one-to-one correspondence to elements in a self-annihilator $P = P^{\perp}$.  This
yields an exact sequence
$$ 0 \to P \to TH_1(M_k) \to \Hom_{\Z}(P;\torsionp ) \to 0 
$$ where the latter homomorphism takes an element $x$ to the homomorphism given by
linking with $x$ in $M_k$.  This is onto since
$\torsionp$ is injective for
$p$-torsion-free finite $\Z$-modules.  Since $P \cong \Hom_{\Z}(P;
\torsionp)$ for any finite $p$-torsion free $\Z$-module, $|TH_1(M_k)| = 2|P| =
2|\SS|$.  
\enddemo

\proclaim{Proposition}\label{containsGil}  Let $K$ be $(1)$\/{\rm -}\/solvable via $W$.  Let $H$ be
chosen as in Proposition~{\rm \ref{.5}.} A character $TH_1(M_k) \to
\torsionp$ extends over $TH_1(W_k)$ if and only if it is a character Gilmer
associated to some $x \in A^k \cap (H \otimes \torsionp)$.
\endproclaim 

\demo{Proof}  As in Gilmer~\cite[pp.~5, 6]{G2}, there are isomorphisms (with\break
$\torsionp$-coefficients) where the first is excision and the second is Poincar\'e
duality;
\begin{equation}\label{AD} H^2(W, W-R)  \cong H^2(R \times I, R \times S^0)
\cong H_2(R \times I, F \times I)
\cong H_2(R,F).\enspace
\end{equation}

Now consider the Mayer-Vietoris sequence 
$$ H^1(W_k;\torsionp) \to \oplus^k H^1(W-R;\torsionp) \to \oplus^k
H^1(R;\torsionp).
$$ Pre-composing the isomorphism~\eqref{AD} with the coboundary homomorphism
$$ H^1(W-R;\torsionp) \to H^2(W, W-R;\torsionp)
$$ we get a commutative diagram as in~\cite{G2}, as follows:
$$
\begin{array}{ccccccc}
{\torsionp} &\hskip-4pt\to\hskip-4pt& {H^1(W_k;\torsionp)}&\hskip-4pt\to\hskip-4pt&{\oplus^k
H_2(R,F;\torsionp)}\\[5pt]
 \big\downarrow&& \big\downarrow{\iota^{\ast}} && \big\downarrow {\oplus
\partial}\\[5pt]
{\torsionp} &\hskip-4pt\to\hskip-4pt&  {H^1(M_k;\torsionp)}
&\hskip-4pt\stackrel{j}{\longrightarrow}\hskip-4pt&\hskip-8pt {\oplus^k H_1(F; \torsionp)}
&\hskip-12pt\stackrel{\alpha}{\longrightarrow}\hskip-4pt& \hskip-2pt{\oplus^k H_1(F;
\torsionp).}\\[5pt]
 &&&& \big\downarrow{\oplus \mu}\\[5pt]
&&&&{\oplus^k H_1(R;\torsionp)}
\end{array}
$$ 
Note that the bottom horizontal sequence is exact, but the top may not be.  The
characters on $TH_1(M_k)$ that extend over $TH_1(W_k)$ lie in one-to-one
correspondence with
$\mbox{image}(j\circ \iota^{\ast})$. As in~\cite{G2}, a diagram chase reveals
that 
$$
\mbox{Image} (j \circ \iota^{\ast}) \subset
\oplus^k (H \otimes \torsionp) \cap \ker{(\a)}.
$$  As mentioned following the exact sequence~\eqref{MV}, these are precisely the
characters Gilmer associated to elements in $A^k \cap (H \otimes
\torsionp)$.
\enddemo

 $(1.5)$-{\it solvability and vanishing Casson\/{\rm -}\/Gordon invariants.}
The following lemma is a straightforward application of Proposition~$3.3$
of~\cite{Let}.

\proclaim{Lemma}\label{contraction}  Let $p$ and $q$ be distinct primes. Let $\psi : C
\to D$ be a homomorphism of finitely generated free 
$\Z G$\/{\rm -}\/modules  where $G =\QQ \rtimes \Z${\rm ,}
$\QQ$ a finite abelian $q$\/{\rm -}\/group{\rm ,}   such that the projection homomorphism
$G = \QQ \rtimes \Z \to \Z$ is abelianization.  Let
$N$ be the index $p^r$ subgroup of $G$ {\rm (}\/unique since
$H_1(G)\break \cong \Z${\rm )}  and let  $N \to \Z_{\l} \times \Z${\rm ,} be a group homomorphism
with finite kernel.  Consider the composition
$$
\Z  N \to \Z [\Z_{\l} \times \Z] \to \f.
$$ If $\psi \otimes_{\Z G} {\rm id}_{\Z}$ is a split monomorphism{\rm ,} then $\psi
\otimes_{\Z N} {\rm id}_{\f}$ is a split monomorphism.
\endproclaim

\demo{Proof}   By hypothesis, $K = \mbox{coker} (\psi \otimes {\rm id}_{\Z})$ is
isomorphic to a summand of $D \otimes_{\Z G} \Z$ having $\Z$-rank equal to $d =
{\rm rank}_{\Z G} (D) - {\rm rank}_{\Z G} (C)$.  We can extend $\psi$ to $\psi' : C
\oplus (\Z G)^d \to D$ so that $\psi' \otimes_{\Z G} {\rm id}_{\Z}$ is an isomorphism.

Let $K' = \mbox{coker} (\psi')$. By the right exactness of the tensor product,
$K' \otimes_{\Z G} \Z =0$.  By~\cite[Prop.~3.3]{Let}
$\Q K' = K'
\otimes_{\Z N} \Q N$ is a finite dimensional rational vector space.  Since the
rational vector space $\mbox{image}(\Q  N \to \f )$ is infinite dimensional, and
since $\Q K' = K'
\otimes_{\Z N} \Q N$ is finite dimensional over
$\Q$, it follows that 
$$
\Q K' \otimes_{\Q N} \f =0.
$$  In fact, given $x \in \Q K'$ there is an element $\hat
\gamma \in \Q N$ with nonzero image $\gamma \in \f$ such that
$x \hat \gamma = 0$.  Thus
$$  x \otimes 1 = x \otimes \gamma \cdot \gamma^{-1} = x \hat \gamma  \otimes
\gamma^{-1} = 0.
$$ Again by right exactness of the tensor product,
$$
\psi' \otimes_{\Z N} {\rm id}_{\f} : (C \oplus (\Z G)^d) \otimes_{\Z N} \f \to D
\otimes_{\Z N} \f
$$  is an epimorphism of finite dimensional $\f$-vector spaces of the same rank
($N$ has finite index in $G$) and so is an isomorphism.  In particular, the
restriction to $C \otimes_{\Z N} \f$ given by $\psi \otimes_{\Z N} {\rm id}_{\f}$ is a
monomorphism.
\enddemo

\proclaim{Theorem}\label{vanish}  Let $K \subset S^3$ be a
$(1.5)$\/{\rm -}\/solvable knot.  Then all Casson\/{\rm -}\/Gordon invariants of
$K$ vanish.  That is{\rm ,} the conclusions of Gilmer\/{\rm '}\/s Theorem~{\rm \ref{CG}} hold.
\endproclaim

\demo{Proof}  Let $W$ be a $(1.5)$-solution for $K$.  By
Proposition~\ref{containsGil}  and the fact that $W$ is a $(1)$-solution, the
character Gilmer associated to an element character $x \in A^k \cap (H
\otimes \torsionp)$
$$
\pi_1(M_k) \to \torsionp \times \Z \to \f
$$ factors through a homomorphism $\pi_1(W_k) \to \f$.

Consider $H_1(W_k;\Z)$ as a $\Z[\Z]$-module, where $\Z$ acts through its quotient
group $\Z_k$ as the group of deck transformations of $W_k$.
$TH_1(W_k)$ is a $\Z[\Z]$-submodule.  The maximal abelian $q$-group $\QQ
\subset TH_1(W_k)$ is a split submodule, so there is an epimorphism
$H_1(W_k;\Z) \to \QQ \times \Z$.

Now suppose $x \in A^k \cap (H \otimes (\torsionp) )$ is an element of order
$\l = q^s$.  Then, since $\QQ$ is a split $\Z[\Z]$-module summand of
$TH_1(W_k)$, $\chi_x :
\pi_1(W_k) \to \Z_{q^t} \times \Z \subset (\torsionp) \times
\Z$ factors through $\QQ \times \Z$, where $t$ is an integer possibly bigger than
$s$.  Hence we have a commutative diagram of ring homomorphisms as in the
following diagram.  Here the vertical homomorphisms are induced by the inclusions
of the index
$p^r$ normal subgroups, and the top horizontal homomorphism extends the Gilmer
associated character $\rho_x : \pi_1(M_k) \to
\f$.  Of course, the bottom horizontal homomorphism factors as described since
$\QQ \rtimes\Z$ is $1$-solvable.
$$
\begin{array}{l} 
{\Z\pi_1(W_k)} \vlrar {\Z H_1(W_k)} \vlrar
{\Z[\QQ \times \Z]} \lrar
{\Z[\Z_{\l} \times \Z]} \lrar {\f} \\[6pt]
\qquad \big\downarrow \phantom{{(W_k)}\times \lrar {\Z H_1(W_k)} \lrar
{\Z[\QQ}} \big\downarrow\\[8pt]
{\Z \pi_1(W)} \to {H_1(W; \2)} \lrar {\Z[\QQ \rtimes
\Z]}.
\end{array}
$$ 
 The Casson-Gordon obstruction is the  difference of
$L$-theory classes of the intersection forms on the second homology of the
following chain complexes
$$  C_{\ast}(W_k; \f) \ \ \ \mbox{  and }\ \ \  C_{\ast}(W_k; \Q).
$$  The intersection form $\lambda_0 \otimes {\rm id}_{\Q}$ on $H_2(W_k;
\Q)$ is trivial in $L_0(\f)$.  Indeed, one easily checks that the lifts into $W_k$
of a basis for the image of the Lagrangian $L \subset H_2(W)$ forms a basis for a
Lagrangian of the intersection form $\lambda_0 \otimes {\rm id}_{\Q}$ on
$H_2(W_k;\Q)$.  Thus it remains to show the form on $H_2(W_k; \f)$ is trivial in
$L_0(\f)$ to show the  Casson-Gordon invariant
$\sigma (K,
\rho_x) = 0$.

Let $\{\l_1, \ldots, \l_m \}$ be a set of immersed spheres in $W^{(2)}$ spanning a
$(2)$-Lagrangian $L \subset H_2(W^{(2)})$ and whose projection to
$W$ forms a basis for a Lagrangian in
$H_2(W; \Z)$.    The $\QQ \times \Z$ cover of $W_k$ is a metabelian cover of
$W$, and so is a quotient space of $W^{(2)}$. Hence the intersection pairing with
$\f$ and $\Z [\QQ \times \Z]$-coefficients vanishes on $L$.

By Lemma~\ref{rankeq}, $\dim_{\f} H_2(W_k;\f) = \dim_{\Q} H_2(W_k;\Q) = 2km$.  The
last equality follows from a dimension count, and the fact that the Euler
characteristic multiplies in covers.  Thus, it suffices to show that the image of
the Lagrangian in $H_2(W_k; \f)$ has dimension $km$.

Let $\LL \to W$ be an immersion of $\vee^m S^2$ obtained by basing the
$\l_i$.  Let $\LL_k$ be the induced $k$-fold cover,  $\LL_k = (\bigvee^{m}
S^2)^k$.   Since
$H_2(\LL_k;\f) \cong \f^{km}$, it suffices to show that $H_2(\LL_k;\f) \to
H_2(W_k;\f)$ is one-to-one.   Since $H_3(W_k;\f) =0$, we must show
$H_3(W_k,
\LL_k;\f) =0$.

But $W_k$ has the homotopy type of a $3$-dimensional CW complex, so this is
equivalent to showing the boundary homomorphism $\partial \otimes {\rm id}_{\f}$ below
is one-to-one, where $N = \QQ \times \Z \subset \QQ \rtimes \Z$.
$$ C_3(W_k,\LL_k) \otimes_{\Z N}\f \stackrel{\partial \otimes {\rm id}_{\f}}{\vlrar} C_2(W_k,
\LL_k) \otimes_{\Z  N}\f.
$$ Since $H_3(W, \LL;\Z) =0$, this follows from Lemma~\ref{contraction}.
\enddemo

{\it Letsche obstructions.}
Recall the recently defined Letsche obstructions to slicing a knot.  Our treatment
is brief, and we refer the reader to~\cite{Let} and~\cite{L3} for more details
on the $\eta$-invariant and the Letsche obstructions. Letsche constructs a
homomorphism 
$$
\eta_K : H_1(M; \Z[\Z])
\times R_{\ast}( \G) \to \R
$$  where $R_{\ast}( \G)$ is the representation ring of $ \G = \left( S^{-1}\Z [\Z
]/\Z [\Z]
\right) \rtimes \Z$ by 
$$
\eta_K(x, \theta) = \tilde \eta_{\theta \circ B\l_x}(M) \in \R
$$ for any $\theta \in R_k( \G)$ and for any $k$.  Here
$B\l_x : H_1(M; \Z [\Z]) \to S^{\-1}\Z[\Z]/\Z[\Z]$ is the homomorphism defined by
$B\l_x(y) = B\l (x, y)$, $B\l$ the Blanchfield pairing for the knot $K$.  Now,
$\tilde \eta_{\theta
\circ B\l_x}(M)$ is the reduced $\eta$-invariant associated to the representation 
$$
\theta \circ B\l_x : \pi_1(M) \to U_k  ,
$$  and $U_k$ is the space of $k$-dimensional unitary representations of the
group~$\G$.   Letsche defines a special subclass of representation
$\PP_k(\pi_1(M))$ as those representations $\theta :
\pi_1(M) \to U_k$ that factor through a nonabelian group of the form $\QQ
\rtimes \Z$ where $\QQ$ is a finite abelian $p$-group and such that the image of
the meridian of the knot group, $\theta (\mu)$ has eigenvalues that are
transcendental over $\Q$.  He proves the following theorem, predating our methods.

\nonumproclaim{Letsche's Theorem} If $K$ is slice{\rm ,} then there is a $P
\subset H_1(M; \Z[\Z])$ such that $P = P^{\perp}$  with respect to the Blanchfield
pairing{\rm ,} and such that for all $x \in P$ and $\theta \in R_k(
\G)$ such that $\theta \circ \a \in \PP_k(\pi_1(M))${\rm ,}
$\eta_K(x, \theta) = 0$.
\endproclaim

\proclaim{Theorem}\label{L=0}  If $K\subset S^3$ is $(1.5)$\/{\rm -}\/solvable{\rm ,} then the
conclusions from Letsche\/{\rm '}\/s theorem above also hold.
\endproclaim

The proof is omitted.

\input Cochran.refs
\end{document}

\end{document}

%% file: Cochran.refs
\AuthorRefNames [MvM]